\documentstyle[12pt, fullpage]{amsart}
\begin{document}

\newcommand{\scroll}[1]{\scrollmode#1\errorstopmode}
\newcommand{\comment}[1]{}
\newif\ifproofmode
\newcommand{\proofmodeQ}{\ifproofmode\let\lebal=\label\renewcommand{\label}[1]{\lebal{##1}\quad\fbox{\smaller\smaller\tt ##1}}\hfuzz=50pt\fi}

\newcommand{\createtitle}[2]{\title{#1}\author{Greg Martin}\address{Department of Mathematics\\University of Toronto\\Canada M5S 3G3}\email{gerg@@math.toronto.edu}\subjclass{#2}\maketitle}
\newcommand{\make}[2]{\section{#2}\setcounter{equation}{0}\label{#1sec}\noindent\input{#1}}

\newtheorem{theorem}{Theorem}[section]
\newtheorem{lemma}[theorem]{Lemma}
\newtheorem{corollary}[theorem]{Corollary}
\newtheorem{proposition}[theorem]{Proposition}
\renewcommand{\theequation}{\thesection.\arabic{equation}}

\newenvironment{pflike}[1]{\medskip\noindent{\bf #1}}{\medskip}
\newenvironment{proof}{\noindent{\bf Proof:}}{\qed\medskip}

\newcommand{\3}[1]{\#\{#1\}}
\renewcommand{\mod}[1]{{\ifmmode\text{\rm\ (mod~$#1$)}
 \else\discretionary{}{}{\hbox{ }}\rm(mod~$#1$)\fi}}
\newcommand{\ep}{\varepsilon}
\renewcommand{\implies}{\Longrightarrow}
\newcommand{\rmif}{\hbox{if\ }}
\newcommand{\floor}[1]{\lfloor#1\rfloor}
\newsymbol\dnd 232D % see The Joy of TeX (2nd ed.), Appendix G
\newcommand{\exdiv}{\mathrel{\|}}
\newcommand{\abs}[1]{|#1|}

\newcommand{\li}{\mathop{\rm li}\nolimits}

\vfuzz=2pt % prohibits overfull \vbox messages due to page headers
	   % and footers

%%%

\newcommand{\lcm}{\mathop{\rm lcm}}
\newcommand{\ord}{\mathop{\rm ord}\nolimits}
\newcommand{\Res}{\mathop{\rm Res}}
\newcommand{\cont}{\mathop{\rm cont}}

\newcommand{\hb}{_{h,b}}
\renewcommand{\hdots}{{h_1\dots h_k}}
\newcommand{\hsdots}{{h_1,\dots,h_k}}
\newcommand{\hhb}{_{h_1\dots h_k,b}}
\newcommand{\hsb}{_{h_1,\dots,h_k;b}}
\newcommand{\6}[1]{^{(#1)}}
\newcommand{\R}{{\cal R}}
\newcommand{\Rhs}{\R(f;h_1,\dots,h_k)}

\newcommand{\gov}[1]{{g(#1)\over #1}}
\newcommand{\glov}[1]{{g(#1)\log #1\over #1}}
\newcommand{\giov}[2]{{g_{#1}(#2)\over #2}}
\newcommand{\gilov}[2]{{g_{#1}(#2)\log #2\over #2}}

\newcommand{\manysum}{\mathop{\sum_{n_1\le x_1}\dots\sum_{n_k\le
x_k}}_{(n_i,n_j)=1 \; (1\le i<j\le k)}}

\newcommand{\state}[3]{
\expandafter\expandafter\expandafter\newcommand\expandafter{\csname#2type\endcsname}{#1}
\expandafter\expandafter\expandafter\newcommand\expandafter{\csname#2cont\endcsname}{{#3}}
\begin{#1}\csname#2cont\endcsname\label{#2}\end{#1}}

\newcommand{\uc}[1]{\uppercase{#1}}
\newcommand{\restate}[1]{\edef\temp{\csname#1type\endcsname}
\medskip\noindent{\bf\expandafter\uc\temp~\ref{#1}}.
{\renewcommand{\label}[1]{}\it\csname#1cont\endcsname}
\ifproofmode\quad\fbox{\smaller\smaller\tt#1}\fi\medskip}

\newcommand{\hyp}{Hypothesis~UH}
\newcommand{\hypo}{\hyp1}
\newcommand{\hypap}{Hypothesis~AP}
\newcommand{\sie}{Sier\-pi\'n\-ski}
\newcommand{\pri}[1]{\prod_{i=1}^{#1}}
\newcommand{\prlk}{\prod_{i=l+1}^k}

\def\agov#1{{\abs{g(#1)}\over{#1}}}
\def\aglov#1{{\abs{g(#1)}\log{#1}\over{#1}}}
\def\gqov#1{{g_q(#1)\over{#1}}}

\proofmodefalse
\proofmodeQ

\createtitle{An Asymptotic Formula for the Number of Smooth Values
of~a~Polynomial}{11N32, 11N25}

\section{Introduction}\setcounter{equation}{0}\label{introsec}\noindent Integers without large prime factors, dubbed {\it smooth numbers\/},
are by now firmly established as a useful and versatile tool in number
theory. More than being simply a property of numbers that is
conceptually dual to primality, smoothness has played a major role in
the proofs of many results, from multiplicative questions to Waring's
problem to complexity analyses of factorization and primality-testing
algorithms. In these last applications, what is needed is an
understanding of the distribution of smooth numbers among the values
taken by a polynomial, which is the subject of this
paper. Specifically, we show a connection between the asymptotic
number of prime values taken by a polynomial and the asymptotic number
of smooth values so taken, showing another way in which these two
properties are more than abstractly linked.

There are conjectures about the distribution of prime values of
polynomials that by now have become standard. Dickson \cite{Dic} first
conjectured that any $K$ linear polynomials with integer coefficients
forming an ``admissible'' set infinitely often take prime values
simultaneously, where $\{L_1,\dots,L_K\}$ is admissible if for every
prime $p$, there exists an integer $n$ such that none of
$L_1(n),\dots,L_K(n)$ is a multiple of $p$; subsequently, Hardy and Littlewood
\cite{HarLit} proposed an asymptotic formula for how often this
occurs. Schinzel and \sie's ``Hypothesis~H'' \cite{SchSie} asserts
that for an admissible set $\{F_1,\dots,F_K\}$ of irreducible
polynomials (integer-valued, naturally) of any degree, there are
infinitely many integers $n$ such that each of $F_1(n),\dots,F_K(n)$
is prime; a quantitative version of this conjecture was first
published by Bateman and Horn \cite{BatHor}. We must introduce some
notation before we can describe the conjectured asymptotic formula,
which we prefer to recast in terms of a single polynomial $F$ rather
than a set $\{F_1,\dots,F_K\}$ of irreducible polynomials.

Let $F(t)=F_1(t)\dots F_K(t)$ be the product of $K$ distinct
irreducible polynomials with integer coefficients. We say that the
polynomial $F$ is {\it admissible\/} if the set $\{F_1,\dots,F_K\}$ is
admissible, that is, if for every prime $p$ there exists an integer
$n$ such that $F(n)$ is not a multiple of $p$. Let $\pi(F;x)$ denote
the number of positive integers $n$ not exceeding $x$ such that each
$F_i(n)$ is a prime (positive or negative). When $F$ is an admissible
polynomial, the size of $\pi(F;x)$ is heuristically $C(F)\li(F;x)$,
these two quantities being defined as
\begin{equation}
C(F) = \prod_p \big( 1-\frac1p \big)^{-K} \big( 1-{\sigma(F;p)\over p}
\big),
\label{CFdef}
\end{equation}
where $\sigma(F;n)$ denotes the number of solutions of $F(a)\equiv0\mod
n$, and
\begin{equation}
\li(F;x) = \int\limits \begin{Sb}0<t<x \\
\min\{|F_1(t)|,\dots,|F_K(t)|\}\ge2\end{Sb}
{dt\over\log|F_1(t)|\dots\log|F_K(t)|}.
\label{liFxdef}
\end{equation}
The second condition of integration is included only to avoid having
to worry about the singularities of the integrand (though we could
have instead defined $\li(F;x)$ using the Cauchy principal value, for
example); we note that $\li(F;x)$ reduces to the familiar logarithmic
integral $\li(x)=\int_2^x dt/\log t$ when $F(t)=t$. At this time, the
only case where the conjecture has been established (even
non-quantitatively) is when $F$ is a linear polynomial, when the
problem reduces to counting primes in a fixed arithmetic progression.

The prime number theorem $\pi(x)\sim\li(x)$ says that the
``probability'' of an integer of size $X$ being prime is $1/\log X$;
if we pretend that the primality of the various $F_i(n)$ are
independent random events, then $\li(F;x)$ would be the expected
number of positive integers $n\le x$ for which all the $F_i(n)$ are
prime. (Of course, for a fixed polynomial $F$, the quantity $\li(F;x)$
is itself asymptotic to $x/\log^Kx$, as appears in \cite{BatHor}; but
the specific expression~(\ref{liFxdef}) should be more accurate
uniformly over polynomials~$F$.) This first guess needs to be
modified, however, to account for the fact that the values of $F$
might be more or less likely to be divisible by a given small prime
$p$ than a ``randomly chosen'' integer of the same size; taking this
factor into account is the role of the constant $C(F)$ (which is a
convergent infinite product---see the discussion following
equation~(\ref{nagel}) below).

We would like to have a similar understanding of the distribution of
smooth values of a polynomial. Let us define
\begin{equation*}
\Psi(F;x,y) = \3{1\le n\le x\colon p\mid F(n)\implies p\le y},
\end{equation*}
the number of $y$-smooth values of $F$ on arguments up to $x$; this
generalizes the standard counting function $\Psi(x,y)$ of $y$-smooth
numbers up to $x$. Known upper bounds (see for
instance~\cite{NairTen}) and lower bounds (see~\cite{DarMarTen}) for
$\Psi(F;x,y)$ do indicate the order of magnitude of $\Psi(F;x,y)$ in
certain ranges; however, in contrast to $\pi(F;x)$, there seems to be
no consensus concerning the expected asymptotic formula for
$\Psi(F;x,y)$.

We can form a probabilistic heuristic for the behavior of $\Psi(F;x,y)$
from the known asymptotic formula
\begin{equation}
\Psi(x,y) \sim x\rho\big( {\log x\over\log y} \big),
\label{Psiasymform}
\end{equation}
where $\rho(u)$ is the Dickman rho-function, defined as the (unique)
continuous solution of the differential-difference equation
$u\rho'(u)=-\rho(u-1)$ for $u\ge1$, satisfying the initial condition
$\rho(u)=1$ for $0\le u\le1$. (We record for later reference that
$\rho(u)=1-\log u$ for $1\le u\le2$.) An interpretation of the
asymptotic formula~(\ref{Psiasymform}) is that a ``randomly chosen''
integer of size $X$ has probability $\rho(u)$ of being
$X^{1/u}$-smooth. Again pretending that the multiplicative properties
of the various $F_i(n)$ are independent of one another, we are led to
the probabilistic prediction that
\begin{equation*}
\Psi(F;x,y) \sim x\prod_{i=1}^K \rho\big( {\log F_i(x)\over\log y}
\big),
\end{equation*}
or equivalently, if we let $d_i$ denote the degree of $F_i$,
\begin{equation}
\Psi(F;x,x^{1/u}) \sim x\rho(d_1u)\dots\rho(d_Ku).
\label{probpred}
\end{equation}
It might seem unclear whether this heuristic needs to include some
sort of dependence on the local properties of $F$, analogous to the
constant $C(F)$ defined above. The purpose of this paper is to
demonstrate, assuming a suitable quantitative version of the
conjectured asymptotic formula for prime values of polynomials, that
the probabilistic prediction~(\ref{probpred}) is indeed the correct
one.

Let us define
\begin{equation}
E(F;x) = \pi(F;x) - C(F)\li(F;x),  \label{Efxdef}
\end{equation}
so that $E(F;x)$ is conjecturally the error term in the asymptotic
formula for prime values of polynomials. Our results will depend upon
the following hypothesis:

\begin{pflike}{\hyp.}
\it Let $d\ge K\ge1$ be integers and $B$ any large positive
constant. Then
\begin{equation}
E(F;t)\ll_{d,B} {C(F)t\over\log^{K+1}t}+1  \label{hypbd}
\end{equation}
uniformly for all integer-valued polynomials $F$, of degree $d$ with
precisely $K$ distinct irreducible factors, whose coefficients are at
most $t^B$ in absolute value.
\end{pflike}
\label{hyppage}

\hyp\ asserts that the quantitative version of Hypothesis H holds with
a considerable range of Uniformity in the coefficients of the
polynomial $F$. We immediately remark that we do not need the full
strength of \hyp\ for our main theorem, but the precise requirement is
much more complicated to state (see equation~(\ref{Hfxydef}) and
Proposition~\ref{usehypprop} below). We do comment, however, that in
the case of linear polynomials ($d=K=1$), \hyp\ is equivalent to a
well-accepted conjecture on the distribution of primes in short
segments of arithmetic progressions (see~\ref{zultimatesec}). Thus the
uniformity required by \hyp\ is not unrealistic.

Notice that if $F$ is not admissible, then there exists a prime $p$
such that $\sigma(p)=p$; therefore $C(F)=0$ and hence
$E(F;t)=\pi(F;t)$. On the other hand, every value of $F$ is then divisible
by $p$, and so the only possible prime values of $F$ are $\pm
p$. Since each of these values can be taken at most $\deg F=d$ times,
we see that $\pi(F;t)\ll_d1$ when $F$ is not admissible. Consequently,
\hyp\ automatically holds for non-admissible polynomials. We have
chosen the form~(\ref{hypbd}) for the hypothesized bound on $E(F;t)$
so that it can be applied without checking in advance whether the
polynomial $F$ is admissible.

We are now ready to state our main theorem.

\begin{theorem}
Assume \hyp. Let $F$ be an integer-valued polynomial, let $K$ be the
number of distinct irreducible factors of $F$, and let $d_1$, \dots,
$d_K$ be the degrees of these factors. Let $d=\max\{d_1,\dots,d_K\}$,
and let $k$ be the number of distinct irreducible factors of $F$ whose
degree equals $d$. Then for any real number $U<(d-1/k)^{-1}$, the
asymptotic formula
\begin{equation}
\Psi(F;x,x^{1/u}) = x\rho(d_1u)\dots\rho(d_Ku) + O_{F,U}\big( \frac
x{\log x} \big)
\label{mainasymform}
\end{equation}
holds uniformly for $x\ge1$ and $0<u\le U$.
\label{mainthm}
\end{theorem}

In particular, if $F$ is an irreducible polynomial of degree $d$, then
for any real number $U<1/(d-1)$, the asymptotic formula
\begin{equation*}
\Psi(F;x,x^{1/u}) = x\rho(du) + O_{F,U}\big( {x\over\log x} \big)
\end{equation*}
holds with the stated uniformity. We note that in the case $d=K=1$,
where $k=1$ and $F(t)=qt+a$ is a linear polynomial, Theorem~\ref{mainthm}
reduces to a well-known asymptotic formula for smooth numbers in a
fixed arithmetic progression, since
\begin{equation*}
\begin{split}
\Psi(F;x,x^{1/u}) &= \Psi(qx+a,x^{1/u};q,a) - \Psi(a,x^{1/u};q,a) \\
&= \Psi(qx,x^{1/u};q,a) + O(1) \\
&= {qx\rho(u)\over q} \big( 1 + O_{q,a}\big( \frac1{\log x} \big)
\big) \\
&= x\rho(u) + O_{F,U}\big( {x\over\log x} \big)
\end{split}
\end{equation*}
(see the survey article~\cite{HilTen}; here $\Psi(x,y;q,a)$ denotes
the number of $y$-smooth numbers up to $x$ that are congruent to
$a\mod q$). This asymptotic formula is known unconditionally to hold
for arbitrarily large values of $U$, which is consistent with the
interpretation of the condition $U<(d-1/k)^{-1}=\infty$ in this case.

Theorem~\ref{mainthm} shows, given what we believe to be true about
prime values of polynomials, that the probabilistic
prediction~(\ref{probpred}) of the asymptotic formula for
$\Psi(F;x,y)$ is indeed valid for a suitable range of $u$. We remark
that the formula~(\ref{mainasymform}) is trivially true for $u<1/d$,
since in this range the smoothness parameter $x^{1/u}$ will
asymptotically exceed the sizes of the factors $F_i(n)$ of $F(n)$,
each of which is $\ll_F x^d$. Since $(d-1/k)^{-1}>1/d$,
Theorem~\ref{mainthm} applies to a nontrivial range of $u$; even
though this range is limited, the theorem provides the first hard
evidence that the probabilistic prediction is the correct one for
general polynomials. The proof of Theorem~\ref{mainthm} comprises the
bulk of this paper and is outlined in Section~\ref{rgosec}.

If we let $\Phi(F;x,y)$ denote the number of primes $q$ up to $x$ for
which $F(q)$ is $y$-smooth, then a similar probabilistic argument
yields the prediction
\begin{equation}
\Phi(F;x,x^{1/u}) \sim \pi(F;x)\rho(d_1u)\dots\rho(d_Ku),
\label{sevenandahalf}
\end{equation}
and the methods used to establish Theorem~\ref{mainthm} could be
extended to prove an analogous result for $\Phi(F;x,x^{1/u})$. Instead
of proceeding in this generality, we prefer to focus on a special case
for which a stronger theorem can be established. Shifted primes
$q+1$ without large prime factors played an important role in the
recent proof by Alford, Granville, and Pomerance \cite{AGP} that there
are infinitely many Carmichael numbers; the counting function of
these smooth shifted primes is precisely $\Phi(F;x,y)$ where
$F(t)=t+1$. More generally, for any nonzero integer $a$ we define
\begin{equation*}
\Phi_a(x,y) = \3{q\le x,\, q\hbox{ prime}\colon p\mid q-a\implies p\le
y},
\end{equation*}
the number of shifted primes $q-a$ with $q$ not exceeding $x$ that are
$y$-smooth. We prove the following theorem:

\begin{theorem}
Assume \hyp. Let $0<U<3$ be a real number and $a\ne0$ an integer. Then
\begin{equation}
\Phi_a(x,x^{1/u}) = \pi(x)\rho(u) + O_{a,U}\big( {\pi(x)\over\log x}
\big)
\label{primeasymform}
\end{equation}
uniformly for $x\ge1+\max\{a,0\}$ and $0<u\le U$.
\label{primethm}
\end{theorem}

Our method (see Section~\ref{linearsec}) is capable in principle of
establishing Theorem~\ref{primethm} for arbitrarily large values of
$U$; however, the combinatorial complexity of the proof would increase
each time $U$ passed an integer. This is because there is no
Buchstab-type identity for smooth shifted primes like the one commonly
used to establish the asymptotic formula~(\ref{Psiasymform}) for
smooth numbers, so one must explicitly write an inclusion-exclusion
formula to characterize the smooth numbers in question. The proof of
Theorem~\ref{primethm} for the stated range of $U$ is indicative of
how the more general result would be established. The fact that the
Dickman function $\rho(u)$ appears unaltered in Theorem~\ref{primethm}
in the expanded range $0<u<3$, where the behavior of $\rho$ is more
complicated, gives further evidence that the probabilistic
predictions~(\ref{probpred}) and~(\ref{sevenandahalf}) are correct in
general, without any modifications resulting from the local properties
of the polynomial~$F$.

As with Theorem~\ref{mainthm}, the precise hypothesis that we require for
establishing Theorem~\ref{primethm} is much less stringent than \hyp. We
remark that we could show, using no new ideas but with a substantial
amount of bookkeeping, that Theorem~\ref{primethm} holds uniformly for
$a\ll x$ if we insert an additional factor of $\log\log x$ into the
error term in the asymptotic formula~(\ref{primeasymform}).  In
particular, by taking $x=a=N$, such a strengthening would
encompass the Goldbach-like problem of counting the number of
representations of an integer $N$ as $N=p+n$ where $p$ is prime and
$n$ is $N^{1/u}$-smooth.

\section{The Structure of the Proof of Theorem~\ref{mainthm}}\setcounter{equation}{0}\label{rgosec}\noindent The proof of Theorem~\ref{mainthm} is quite long and, unfortunately,
mired in technical details in many places. In the interest of making
the overall structure of the proof more evident, we have broken the
proof down into several propositions, precisely stated in this
section, from which Theorem~\ref{mainthm} will be deduced near the end
of the section. These propositions will themselves be proved in the
sections to come.

In descriptive terms, the proof of Theorem~\ref{mainthm} proceeds
according to the following outline:

\begin{enumerate}
\item We show that we may assume, without loss of generality, that the
polynomial $F$ satisfies several technical conditions on its
factorization and on its local properties. This reduction is described
in Proposition~\ref{salvationprop} below.
\item Using an inclusion-exclusion argument on the factors of $F$, we
write $\Psi(F;x,y)$ as a combination (see Proposition~\ref{introMprop})
of terms of the form $M(f;x,y)$, defined in equation~(\ref{Mfxydef})
below, where $f$ runs through divisors of~$F$.
\item We establish a ``pseudo-asymptotic formula'' (where the behavior
of one of the primary terms $H(f;x,y)$, defined in
equation~(\ref{Hfxydef}) below, is unknown), given in
Theorem~\ref{Masymthm}, for the expressions $M(f;x,y)$. The proof of
this theorem is itself broken into several steps:
\begin{enumerate}
\item We convert $M(f;x,y)$ into a combination of terms of the form
$\pi(f_i;x_i)$, where the $f_i$ are drawn from a family of polynomials
determined by $f$, as described in Proposition~\ref{alterMprop};
\item We apply Brun's sieve method to establish in Proposition~\ref{brunprop}
an upper bound for the smaller terms of this form;
\item By investigating the relationships among the local properties of
the polynomials in this family, we evaluate a certain finite sum (see
Proposition~\ref{cflbprop}) involving constants of the type $C(f_i)$
defined in equation~(\ref{CFdef});
\item We establish an asymptotic formula for a weighted mean value of the
multiplicative function that arises in the previous step, which is the
content of Proposition~\ref{lipsprop}.
\end{enumerate}
\item Next, we show that \hyp\ implies an upper bound for the term
$H(f;x,y)$, as given in Proposition~\ref{usehypprop}, and thus
converts the pseudo-asymptotic formula from Theorem~\ref{Masymthm}
into a true asymptotic formula for the $M(f;x,y)$.
\item Finally, since $\Psi(F;x,y)$ has been expressed as a combination
of terms of the form $M(f;x,y)$ in Proposition~\ref{introMprop}, we
can recast Theorem~\ref{Masymthm} into our final goal, an asymptotic
formula for $\Psi(F;x,y)$.
\end{enumerate}

The rest of this section will be devoted to making this outline
rigorous and establishing Theorems~\ref{Masymthm} and~\ref{mainthm}
assuming the validity of the propositions to be stated. We begin by
defining the technical conditions we would like to impose on the
polynomial $F$. Other than the first two terms, we make no claims that
the following terminology is (or should become) standard---we simply
wish to use concise, reasonably provocative words rather than include
long, awkward phrases at every turn.

\begin{itemize}
\item Recall that a polynomial is {\it squarefree\/} if it is not
divisible by the square of any nonconstant polynomial.
\item Recall also that a polynomial with integer coefficients is {\it
primitive\/} if the greatest common divisor of its coefficients is
1. A polynomial is primitive if and only if each of its irreducible
factors is primitive, and the factorization of a primitive polynomial
over the rationals is the same as its factorization over the integers
(both these statements are consequences of Gauss' Lemma).
\item We say that a polynomial is {\it balanced\/} if it is the
product of distinct irreducible polynomials all of the same degree, so
that in particular, a balanced polynomial is squarefree.
\item We say that an irreducible polynomial $g$ is {\it effective\/}
if $g(0)\ge2$ and $g'(t)\ge1$ for all $t\ge0$, and we say that an
arbitrary polynomial $f$ is effective if all of its irreducible
factors are effective. This certainly implies that $f(0)\ge2$ and
$f'(t)\ge1$ for all $t\ge0$, and hence $f(t)>t$ for all $t\ge0$.
\item In agreement with the use of the term in Section~\ref{introsec},
we call a polynomial with integer coefficients {\it admissible\/} if
it takes at least one nonzero value modulo every prime. In particular,
any admissible polynomial is primitive.
\item Finally, we say that a polynomial $f$ with integer coefficients
is {\it exclusive\/} if no two distinct irreducible factors of $f$
have a common zero modulo any prime. A primitive polynomial $f$ is
exclusive if and only if the resultant $\Res(g,h)$ of $g$ and $h$
equals 1 for every pair $g,h$ of distinct irreducible factors of $f$.
\end{itemize}

Now that this terminology is in place, we can state the proposition
that allows us to place various restrictions on the polynomial $F$ in
Theorem~\ref{mainthm}.

\state{proposition}{salvationprop}
{Suppose that Theorem~\ref{mainthm} holds for any polynomial with
integer coefficients that is balanced, effective, admissible, and
exclusive. Then it holds for any integer-valued polynomial.}

\noindent This proposition is established in Section~\ref{vationsec}
by considering an arbitrary polynomial restricted to suitable
arithmetic progressions.

Next we introduce the quantity $M(f;x,y)$ for which we shall
establish a pseudo-asymptotic formula in Theorem~\ref{Masymthm}. Given a
primitive polynomial $f$, define
\begin{multline}
M(f;x,y) = \#\{ 1\le n\le x\colon \hbox{for each irreducible factor
$g$ of $f$,} \\
\hbox{there exists a prime $p>y$ such that $p\mid g(n)$} \}.
\label{Mfxydef}
\end{multline}
The connection between smooth values of polynomials and these
expressions $M(f;x,y)$ is given in the following proposition:

\state{proposition}{introMprop}
{Let $F$ be a primitive polynomial, and let $F_1,\dots,F_K$ denote the
distinct irreducible factors of $F$. Then for any real numbers $x\ge
y\ge1$,
\begin{equation*}
\Psi(F;x,y) = \floor x + \sum_{1\le k\le K} (-1)^k \sum_{1\le
i_1<\dots<i_k\le K} M(F_{i_1}\dots F_{i_k};x,y).
\end{equation*}}

The statement of the next proposition requires a bit more notation. In
fact, the following notation will be used in the statements of
Propositions~\ref{alterMprop}---\ref{stupidOprop},
Theorem~\ref{Masymthm}, and Proposition~\ref{usehypprop}:
\begin{itemize}
\item Let $f$ be a polynomial with integer coefficients that is
primitive and balanced. Let $k$ denote the number of irreducible
factors of $f$, let $f_1,\dots,f_k$ denote these irreducible factors,
and let $d$ denote their common degree.
\item Let $u$ and $U$ be real numbers in the range
\begin{equation}
1/d\le u\le U<\min\{(d-1/k)^{-1},2/d\}.  \label{Urange}
\end{equation}
\item Given a real number $x\ge1$, set $y=x^{1/u}$ and $\xi_i=f_i(x)$
for $1\le i\le k$. Note that $\xi_i\asymp x^d$ as $x$ goes to
infinity (where $X\asymp Y$ means $Y\ll X\ll Y$), and that the upper
bound~(\ref{Urange}) on $U$ implies that $y\ge\sqrt{\xi_i}$ when $x$ is
sufficiently large in terms of $f$ and $U$.
\item As in Section~\ref{introsec}, let $\pi(f;t)$ denote the
number of integers $1\le n\le t$ for which each $f_i(n)$ is prime.
\item Given positive integers $\hsdots$, define
\begin{equation}
\Rhs = \{ 1\le b\le\hdots\colon h_i|f_i(b)\hbox{ for each }1\le i\le k\}.
\label{Rhsdef}
\end{equation}
For any element $b$ of $\Rhs$, define the polynomial $f\hhb$ by
\begin{equation}
f\hhb(t) = {f(\hdots t+b)\over\hdots}.  \label{fhhbdef}
\end{equation}
We remark that $f\hhb$ actually has integer coefficients (see the
remarks following equation~(\ref{taylorF}) below).
\end{itemize}
\label{defspage}
When we prove Propositions~\ref{alterMprop}---\ref{stupidOprop}
and~\ref{usehypprop} in later sections, we shall adhere to this notation
as well. With the notation in place, we can describe the connection
between the expressions $M(f;x,y)$ defined in equation (\ref{Mfxydef})
and prime values of polynomials.

\state{proposition}{alterMprop}
{Let $f$ be a polynomial that is balanced, effective, primitive, and
exclusive. Then when $x$ is sufficiently large we have
\begin{equation}
M(f;x,y) = \mathop{\sum_{h_1\le\xi_1/y} \dots
\sum_{h_k\le\xi_k/y}}_{(h_i,h_j)=1 \; (1\le i<j\le k)} \sum_{b\in\Rhs}
\big( \pi\big( f\hhb;{x-b\over\hdots} \big) - \pi(f\hhb;\eta_\hsdots)
\big) \label{twopis} \end{equation} for certain quantities
$\eta_\hsdots$ that satisfy
$\eta_\hsdots\asymp(y\max\{\hsdots\})^{1/d}(\hdots)^{-1}$.}

\noindent Together, Propositions~\ref{introMprop} and~\ref{alterMprop}
provide the link between smooth values of polynomials and prime values
of polynomials. These propositions are combinatorial in nature, and
their proofs are given in Section~\ref{trashsec}.

Of the two terms of the form $\pi(f;t)$ in each summand on the
right-hand side of equation~(\ref{twopis}), the first term is the
significant one. The following proposition provides a tidy bound
for the contribution to the sum from the second terms (those
containing the quantities $\eta_\hsdots$).

\state{proposition}{brunprop}
{Let $f$ be a polynomial that is balanced, effective, and
primitive. If the quantities $\eta_\hsdots$ satisfy
$\eta_\hsdots\asymp(y\max\{\hsdots\})^{1/d}(\hdots)^{-1}$, then
\begin{equation*}
\mathop{\sum_{h_1\le\xi_1/y} \dots \sum_{h_k\le\xi_k/y}}_{(h_i,h_j)=1
\; (1\le i<j\le k)} \sum_{b\in\Rhs} \pi(f\hhb;\eta_\hsdots) \ll_f
{x\over\log x}.
\end{equation*}}

\noindent This proposition is established by an
application of Brun's upper bound sieve method to each summand, as we
show in Section~\ref{brunsec}. The only complication is keeping track
of the dependence on $\hsdots,$ and $b$ in the bounds obtained; the
multiple sum over $\hsdots$ that results is then bounded by a
mean-value theorem for multiplicative functions
(Proposition~\ref{mvprop} in Appendix~A).

The next proposition is, in a sense, the most important step in the
proof, as it relates a complicated expression involving the numbers of
local roots of a family of polynomials to a value of an explicit
multiplicative function. To state this proposition, we recall from
Section~\ref{introsec} that $\sigma(f;n)$ denotes the number of
solutions of $f(a)\equiv0\mod n$, and we define two related 
multiplicative functions as follows:
\begin{equation}
G(f;n) = \prod_{p\mid n} \big( 1-{\sigma(f;p)\over p} \big)^{-1}
\quad\hbox{and}\quad \sigma^*(f;n) = \prod_{p^\nu\exdiv n} \big(
\sigma(f;p^\nu) - {\sigma(f;p^{\nu+1})\over p} \big).
\label{Gsigdefs}
\end{equation}
If we assume that $f$ is admissible, then $\sigma(f;p)<p$ for all primes
$p$, and hence $G(f;p)$ is well-defined. We remark that both of these
multiplicative functions are nonnegative (for $\sigma^*$ this is not
transparent---see equation~(\ref{sigstarnonneg}) below) and that
$\sigma^*(f;n)\le\sigma(f;n)$.

We also recall the definition~(\ref{liFxdef}) of $\li(f;x)$, and for
effective polynomials $f$ we define the modified function
\begin{equation}
\li_\hsdots(f;x) = \int\limits \begin{Sb}0<t<x \\
\min\{f_1(t)/h_1,\dots,f_k(t)/h_k\}\ge2\end{Sb}
{dt\over\log(f_1(t)/h_1)\dots\log(f_k(t)/h_k)}
\label{lihsFxdef}
\end{equation}
for any positive integers $\hsdots$. We can now state the following
proposition:

\state{proposition}{cflbprop}
{Let $f$ be a polynomial that is balanced, effective, admissible, and
exclusive, and let $\hsdots$ be pairwise coprime positive
integers. Then
\begin{equation}
\begin{split}
\sum_{b\in\Rhs} C & (f\hhb) \li\big( f\hhb; {x-b\over\hdots} \big)
\\
&= {C(f) G(f;h_1\dots h_k)\sigma^*(f_1;h_1)\dots\sigma^*(f_k;h_k)
\li_\hsdots(f;x) \over \hdots} \\
&\qquad+ O_f(G(f;h_1\dots h_k)\sigma^*(f_1;h_1)\dots\sigma^*(f_k;h_k))
\label{arises}
\end{split}
\end{equation}
uniformly for $x\ge1$.}

\noindent As one might guess, the sum on the left-hand side
equation~(\ref{arises}) arises from the first terms in the summands of
equation~(\ref{twopis}). The $\li(\cdot)$ factor is quite benign, as
it does not depend very much on $b$; the important part of the
evaluation, as we shall see Section~\ref{multappendixsec}, is to
understand the relationship between the constants $C(f\hhb)$ and
$C(f)$ itself.

We have two more propositions to state before coming to our proof of
the pseudo-asymp\-totic formula for $M(f;x,y)$. In that proof, we
shall need to sum the right-hand side of equation~(\ref{arises}) over
the possible values of the $h_i$. The first of these two propositions
gives us an asymptotic formula for the resulting sum of the main term
from~(\ref{arises}), while the second proposition estimates the
resulting error term.

\state{proposition}{lipsprop}
{Let $f$ be a polynomial that is balanced, effective, admissible, and
exclusive. Then
\begin{multline}
\mathop{\sum_{h_1\le\xi_1/y}\dots\sum_{h_k\le\xi_k/y}}_{(h_i,h_j)=1 \;
(1\le i<j\le k)} {G(f;\hdots)\sigma^*(f_1;h_1)\dots
\sigma^*(f_k;h_k) \li_\hsdots(f;x) \over \hdots} \\
= C(f)^{-1} x \log^k(du) + O_{f,U}\big( {x\over\log x} \big)
\label{werenot}
\end{multline}
uniformly for $x\ge1$ and $0<u\le U$.}

\noindent Proposition~\ref{lipsprop} is definitely the longest part
of the proof, taking all of Section~\ref{psappendixsec} to complete and
requiring the machinery of Appendix~A, although the difficulties are
technical rather than conceptual. Even ignoring the presence of the
$\li_\hsdots$ term in the sum on the left-hand side of
equation~(\ref{werenot}), that sum is made more difficult to analyze
because of the coprimality condition of summation---otherwise, the
multiple sum would just be a $k$-fold product of sums of multiplicative
functions of one variable. Proposition~\ref{mv2prop} below provides a
general asymptotic formula for multiple sums of multiplicative functions
where the variables of summation are restricted in this way. Once we
understand the behavior of that sum with the $\li_\hsdots$ term omitted
(see Lemma~\ref{cancellem}), we employ a lengthy partial summation
argument to assess the effect of the $\li_\hsdots$ term on the actual sum.

The second of these two propositions is much easier to establish than
the first, due to the fact that we seek only an upper bound rather
than an asymptotic formula. This means that the coprimality condition
of summation can be omitted and a mean-value theorem for
multiplicative functions (Proposition~\ref{mvprop}) invoked in the proof,
which can be found in Section~\ref{1multsec}.

\state{proposition}{stupidOprop}
{Let $f$ be a polynomial that is balanced and admissible. Then
\begin{equation*}
\mathop{\sum_{h_1\le\xi_1/y}\dots\sum_{h_k\le\xi_k/y}}_{(h_i,h_j)=1 \;
(1\le i<j\le k)} G(f;\hdots)\sigma^*(f_1;h_1) \dots \sigma^*(f_k;h_k)
\ll_{f,U} {x\over\log x}
\end{equation*}
uniformly for $x\ge1$ and $0<u\le U$.}

Assuming the validity of
Propositions~\ref{alterMprop}--\ref{stupidOprop}, we can at last state
and prove our main result concerning the expressions $M(f;x,y)$. We
have repeated from page~\pageref{defspage} the definitions of the
various symbols so as to have the statement of the theorem be
self-contained.

\begin{theorem}
Let $f$ be a polynomial that is balanced, effective, admissible, and
exclusive. Let $k$ denote the number of irreducible factors of $f$,
let $f_1,\dots,f_k$ denote these irreducible factors, and let $d$
denote their common degree. Let $u$ and $U$ be real numbers in the
range~(\ref{Urange}). Given a real number $x\ge1$, set $y=x^{1/u}$
and $\xi_i=f_i(x)$ for $1\le i\le k$. If $M(f;x,y)$ is defined as in
equation~(\ref{Mfxydef}) and $E(f;x)$ is defined as in
equation~(\ref{Efxdef}), then we have
\begin{equation}
M(f;x,y) = x\log^k(du) + H(f;x,y) + O_{f,U}\big( \frac x{\log x} \big)
\label{Masymeq}
\end{equation}
uniformly for $x\ge1$ and $0<u\le U$, where $H(f;x,y)$ is defined by
\begin{equation}
H(f;x,y) = \mathop{\sum_{h_1\le\xi_1/y}\dots
\sum_{h_k\le\xi_k/y}}_{(h_i,h_j)=1\, (1\le i<j\le k)} \sum_{b\in\Rhs}
E\big( f\hhb;{x-b\over\hdots} \big).
\label{Hfxydef}
\end{equation}
\label{Masymthm}
\end{theorem}

\noindent Of course, it is the unknown size of $H(f;x,y)$ that keeps
equation~(\ref{Masymeq}) from being an unconditional asymptotic
formula for $M(f;x,y)$, which would lead to an unconditional
asymptotic formula for $\Psi(F;x,y)$ using
Proposition~\ref{introMprop}; it is to control the size of $H(f;x,y)$
that \hyp\ will be used (see the proof of Proposition~\ref{usehypprop}
below).

\medskip
\begin{proof}
Proposition~\ref{alterMprop} tells us that
\begin{equation*}
M(f;x,y) = \mathop{\sum_{h_1\le\xi_1/y} \dots
\sum_{h_k\le\xi_k/y}}_{(h_i,h_j)=1 \; (1\le i<j\le k)} \sum_{b\in\Rhs}
\big( \pi\big( f\hhb;{x-b\over\hdots} \big) - \pi(f\hhb;\eta_\hsdots)
\big).
\end{equation*}
From the definition~(\ref{Efxdef}) of $E(f;x)$, we can write this as
\begin{equation}
\begin{split}
M(f;x,y) &= \mathop{\sum_{h_1\le\xi_1/y} \dots
\sum_{h_k\le\xi_k/y}}_{(h_i,h_j)=1 \; (1\le i<j\le k)} \sum_{b\in\Rhs}
C(f\hhb) \li\big( f\hhb;{x-b\over\hdots} \big) \\
&\qquad+ \mathop{\sum_{h_1\le\xi_1/y} \dots
\sum_{h_k\le\xi_k/y}}_{(h_i,h_j)=1 \; (1\le i<j\le k)} \sum_{b\in\Rhs}
E\big( f\hhb;{x-b\over\hdots} \big) \\
&\qquad- \mathop{\sum_{h_1\le\xi_1/y} \dots
\sum_{h_k\le\xi_k/y}}_{(h_i,h_j)=1 \; (1\le i<j\le k)} \sum_{b\in\Rhs}
\pi(f\hhb;\eta_\hsdots).
\end{split}
\label{M1rears}
\end{equation}
If we define
\begin{equation}
M_1(f;x,y) = \mathop{\sum_{h_1\le\xi_1/y} \dots
\sum_{h_k\le\xi_k/y}}_{(h_i,h_j)=1 \; (1\le i<j\le k)} \sum_{b\in\Rhs}
C(f\hhb) \li\big( f\hhb;{x-b\over\hdots} \big)
\label{M1def}
\end{equation}
and use the definition~(\ref{Hfxydef}) of $H(f;x,y)$,
equation~(\ref{M1rears}) becomes
\begin{equation*}
M(f;x,y) = M_1(f;x,y) + H(f;x,y) + O_{f,U}\big( {x\over\log x} \big),
\end{equation*}
where we have used Proposition~\ref{brunprop} to estimate the final
sum in~(\ref{M1rears}). It therefore suffices to show that
\begin{equation}
M_1(f;x,y) = x\log^k(du) + O_{f,U}\big( \frac x{\log x} \big).
\label{M1suffices}
\end{equation}

The inner sum in the definition~(\ref{M1def}) can be evaluated by
Proposition~\ref{cflbprop}, yielding
\begin{equation*}
\begin{split}
M_1(f;x,y) &= C(f) \mathop{\sum_{h_1\le\xi_1/y} \dots
\sum_{h_k\le\xi_k/y}}_{(h_i,h_j)=1 \; (1\le i<j\le k)} {G(f;h_1\dots
h_k)\sigma^*(f_1;h_1)\dots\sigma^*(f_k;h_k) \li_\hsdots(f;x) \over
\hdots} \\
&\qquad+ O\bigg( \mathop{\sum_{h_1\le\xi_1/y} \dots
\sum_{h_k\le\xi_k/y}}_{(h_i,h_j)=1 \; (1\le i<j\le k)} G(f;h_1\dots
h_k)\sigma^*(f_1;h_1)\dots\sigma^*(f_k;h_k) \bigg).
\end{split}
\end{equation*}
Using the asymptotic formula given by Proposition~\ref{lipsprop} for
the main term in this last expression and the estimate in
Proposition~\ref{stupidOprop} to bound the error term, we see that
\begin{equation*}
\begin{split}
M_1(f;x,y) &= C(f) \big( C(f)^{-1} x \log^k(du) + O_{f,U}\big(
{x\over\log x} \big) \big) + O_{f,U}\big( {x\over\log x} \big) \\
&= x \log^k(du) + O_{f,U}\big( {x\over\log x} \big),
\end{split}
\end{equation*}
which establishes equation~(\ref{M1suffices}) and hence the theorem.
\end{proof}

As we mentioned before, \hyp\ implies an upper bound for the
expression $H(f;x,y)$ defined in equation~(\ref{Hfxydef}). The
following proposition, proved in Section~\ref{1multsec} using again
the mean-value theorem for multiplicative functions
(Proposition~\ref{mvprop}), provides the needed estimate.

\state{proposition}{usehypprop}
{Assume \hyp. Let $f$ be a polynomial with integer coefficients that
is balanced, effective, and admissible. Then $H(f;x,y)\ll_{f,U} x/\log
x$.}

If we assume for now the validity of Theorem~\ref{Masymthm} and
Propositions~\ref{salvationprop},~\ref{introMprop},
and~\ref{usehypprop}, we have all the tools we need to establish our
main theorem.

\begin{pflike}{Proof of Theorem~\ref{mainthm}:}
As mentioned after the statement of Theorem~\ref{mainthm}, the theorem is
already known unconditionally in the case $d=K=1$, and so we may assume
that $dK\ge2$. First suppose that the polynomial $F$ has integer
coefficients and is balanced, effective, admissible, and exclusive; in
this case the statement of Theorem~\ref{mainthm} reduces to
\begin{equation}
\Psi(F;x,y) = x \rho(du)^K + O_{F,U}\big ( {x\over\log x} \big).
\label{reducesto}
\end{equation}
Now Proposition~\ref{introMprop} tells us that
\begin{equation}
\Psi(F;x,y) = \floor x + \sum_{1\le k\le K} (-1)^k \sum_{1\le
i_1<\dots<i_k\le K} M(F_{i_1}\dots F_{i_k};x,y).
\label{thistime}
\end{equation}
We note that any divisor of a polynomial that is balanced, effective,
admissible, and exclusive will itself have those four properties.
Therefore, setting $f$ to equal any of these polynomials $F_{i_1}\dots
F_{i_k}$ (so that the convention that $k$ denotes the number of
irreducible factors of $f$ is consistent with the use of $k$ in
equation~(\ref{thistime})), we see from Theorem~\ref{Masymthm} that
\begin{equation}
M(f;x,y) = x\log^k(du) + H(f;x,y) + O_{f,U}\big( \frac x{\log x}
\big).
\label{absorbed}
\end{equation}
But from Proposition~\ref{usehypprop}, under \hyp\ we have
$H(f;x,y)\ll_{f,U} x/\log x$, and therefore $H(f;x,y)$ can be absorbed
into the error term in the asymptotic formula~(\ref{absorbed}). With
this formula, equation~(\ref{thistime}) becomes
\begin{equation*}
\Psi(F;x,y) = \floor x + \sum_{1\le k\le K} (-1)^k \sum_{1\le
i_1<\dots<i_k\le K} \big( x\log^k(du) + O_{F,U}\big( \frac x{\log x}
\big) \big).
\end{equation*}
For any $1\le k\le K$ there are $K\choose k$ ways to choose the $k$
indices $i_1,\dots,i_k$, and so this becomes
\begin{equation*}
\begin{split}
\Psi(F;x,y) &= x + \sum_{1\le k\le K} (-1)^k {\textstyle{K\choose k}}
x\log^k(du) + O_{F,U}\big ( {x\over\log x} \big) \\
&= x (1-\log du)^K + O_{F,U}\big( {x\over\log x} \big)
\end{split}
\end{equation*}
by the binomial theorem. But $1\le du\le 2$ when $u$ is in the
range~(\ref{Urange}), and so $\rho(du)=1-\log du$ as remarked
earlier. Therefore this last equation is equivalent to
equation~(\ref{reducesto}), which establishes Theorem~\ref{mainthm}
for polynomials with integer coefficients that are balanced,
effective, admissible, and exclusive. But then by
Proposition~\ref{salvationprop}, Theorem~\ref{mainthm} holds for all
integer-valued polynomials.
\qed
\end{pflike}

This completes the proof of Theorem~\ref{mainthm} modulo the proofs
of Propositions~\ref{salvationprop} through~\ref{stupidOprop}
and~\ref{usehypprop}, which will be given in the next seven
sections. As remarked in the introduction, Theorem~\ref{primethm} will
be addressed in Section~\ref{linearsec}, while Appendix~A contains
mean-value results for multiplicative functions, and~\ref{zultimatesec} is
devoted to showing that \hyp\ in the case of linear polynomials is
equivalent to a certain statement concerning the number of primes in
short segments of arithmetic progressions.

\section{Technical Conditions on the Polynomial $F$}\setcounter{equation}{0}\label{vationsec}\noindent Our goal for this section is to establish Proposition~\ref{salvationprop}.
The idea is to relate the number of smooth values of a given
integer-valued polynomial $F$ to the number of smooth values of certain
polynomials having the properties in the statement of the proposition.
The tricky properties are admissibility and exclusiveness; the other
properties are addressed by the following elementary lemma.

\begin{lemma}
Let $F(t)$ be a nonconstant integer-valued polynomial, let $d$ be the
largest of the degrees of the irreducible factors of $F$, and let $k$
be the number of distinct irreducible factors of $F$ of degree
$d$. Let $\alpha$ be a real number exceeding $d-1$. There exists an
effective polynomial $F_1(t)$ with integer coefficients that is the
product of $k$ distinct irreducible polynomials of degree $d$, such
that
\begin{equation}
\Psi(F_1;x,y) = \Psi(F;x,y) + O_{F,\alpha}(1)
\label{baleffeq}
\end{equation}
uniformly for $x\ge1$ and $y\ge x^\alpha$.
\label{balefflem}
\end{lemma}

\begin{proof}
We remark that it suffices to show that equation~(\ref{baleffeq})
holds when $x$ is sufficiently large, by adjusting the constant
implicit in the $O$-notation if necessary. Recall that the {\it
content\/} of a polynomial $F$ is the greatest common divisor of all
of its coefficients (so that a polynomial is primitive precisely when
its content equals 1). Since $F(t)$ is integer-valued, the
coefficients of $F$ are all rational, so we can choose a positive
integer $m$ such that $mF(t)$ has integer coefficients. Write
\begin{equation}
mF(t) = \pm c_0 G_1(t)^{\mu_1}\dots G_k(t)^{\mu_k} H_1(t)^{\nu_1}\dots
H_l(t)^{\nu_l}
\label{mfteq}
\end{equation}
where $c_0$ is the content of $mF(t)$, the $\mu_i$ and $\nu_i$ are
positive integers, and the $G_i$ and $H_i$ are distinct primitive,
irreducible polynomials with positive leading coefficients satisfying
$\deg G_i=d$ for $1\le i\le k$ and $\deg H_i\le d-1$ for $1\le i\le
l$.

Set $F_0(t) = G_1(t)\dots G_k(t)$, so that $F_0$ is a balanced
polynomial with integer coefficients. If $n\le x$ is an integer such
that $F(n)\ne0$, then it is clear from equation~(\ref{mfteq}) that the
largest prime factor $p$ of $F_0(n)$ is the same as the largest prime
factor of $F(n)$ provided that $p$ exceeds all prime divisors of
$mc_0H_1(n)\dots H_l(n)$. In particular, as long as
\begin{equation}
y>\max\{m,c_0,|H_1(n)|,\dots,|H_l(n)|\},  \label{yineq}
\end{equation}
then $F_0(n)$ is $y$-smooth precisely when $F(n)$ is $y$-smooth. But
$y\ge x^\alpha$ and each $H_i(x)\ll_F x^{d-1}$, so we see that the
inequality~(\ref{yineq}) always holds when $x$ is sufficiently large
(in terms of $F$ and $\alpha$), since $\alpha>d-1$. Therefore
\begin{equation}
\Psi(F_0;x,y)=\Psi(F;x,y)+O_{F,\alpha}(1),  \label{step0}
\end{equation}
the error arising from values of $n$ for which $F(n)=0$, of which
there can be at most $\deg F$.

Now we choose a positive real number $t_0=t_0(F)$ such that for each
$1\le i\le k$, we have $G_i(t_0)\ge2$ and $G_i'(t)\ge1$ for all $t\ge
t_0$. If we set $F_1(t)=F_0(t+t_0)$, then $F_1$ is again a balanced
polynomial with integer coefficients, and moreover $F_1$ is effective
by our choice of $t_0$. Finally, we see that
\begin{equation*}
\Psi(F_1;x,y) = \Psi(F_0;x+t_0,y)-\Psi(F_0;t_0,y) =
\Psi(F_0;x,y)+O_F(1).
\end{equation*}
This together with equation~(\ref{step0}) establishes the lemma.
\end{proof}

To address the properties of admissibility and exclusiveness, we
consider the restriction of a polynomial $F$ to an arithmetic
progression. If $Q$ is a positive integer and $a$ is any integer, we
can use the Taylor expansion of $F$ at $a$ to see that
\begin{equation}
F(Qt+a) = F(a) + QF'(a)t + Q^2{F''(a)\over2}t^2 + \dots +
Q^D{F^{(D)}(a)\over D!}t^D,
\label{taylorF}
\end{equation}
where $D=\deg F$; note that each of the expressions $F^{(j)}(a)/j!$ is an
integer. If $a$ is chosen so that $Q$ divides $F(a)$, we see that
every coefficient on the right-hand side of~(\ref{taylorF}) is
divisible by $Q$, whence $F(Qt+a)/Q$ is a polynomial with integer
coefficients.

Even if $f(a)$ is not a multiple of $Q$, the coefficients of $F(Qt+a)$
might all be divisible by some common factor which we would like to
remove. If we divide all of the coefficients of a polynomial $F$ by
$\cont F$, we call the resulting polynomial the {\it primitivization\/}
of $F$. The following lemma shows that if we consider a polynomial
restricted to an arithmetic progression to a suitable modulus, then its
primitivization has the desired properties of admissibility and
exclusivity.

\begin{lemma}
Let $F_1$ be a squarefree polynomial with integer coefficients, and
let $D$ denote the degree of $F_1$ and $\Delta$ the discriminant of
$F_1$. Let $Q$ be a positive integer satisfying:
\begin{itemize}
\item for all primes $p\le D$, we have $\ord_p(Q) > \ord_p(\Delta)$;
\item for every pair $G_1,H_1$ of distinct nonconstant irreducible
factors of $F_1$, and for all primes $p$ dividing their resultant
$\Res(G_1,H_1)$, we have $\ord_p(Q) > \ord_p(\Res(G_1,H_1))$.
\end{itemize}
Let $a$ be any positive integer, let $F_2$ denote the polynomial
$F_2(t)=F_1(Qt+a)$, and let $F_3$ denote the primitivization of
$F_2$. Then $F_3$ is admissible and exclusive.
\label{salvationlem}
\end{lemma}

\noindent We remark that it is always possible to find an integer $Q$
satisfying the hypotheses of the lemma when $F_1$ is squarefree, since
this implies that $\Delta\ne0$ and $\Res(G_1,H_1)\ne0$ for any two
distinct nonconstant irreducible factors $G_1$ and $H_1$ of $F_1$. In
fact, the smallest such $Q$ is
\begin{equation}
Q = Q(F_1) = \prod_{p\le D} p^{\ord_p(\Delta)+1} \prod
\begin{Sb}G_1\ne H_1{\rm\,irreducible} \\ G_1H_1\mid F_1\end{Sb}
\bigg( \prod_{p\mid\Res(G_1,H_1)} p^{\ord_p(\Res(G_1,H_1))+1} \bigg).
\label{smallQ}
\end{equation}

\medskip
\begin{proof}
We shall strive first for admissibility and then for exclusiveness (as
if we were climbing the social ladder). If $p$ is a prime greater than
$D$, then the primitivization $F_3$ of $F_2$ is not the zero
polynomial\mod p, since any primitive polynomial has at least one
nonzero coefficient modulo every prime; and in fact $F_3$ has degree
at most $D$. Therefore $F_3$ has at most $D$ zeros\mod p, and so takes
at least one nonzero value\mod p.

If $p$ is a prime not exceeding $D$, then from the
identity~(\ref{taylorF}) applied to $F_1$, we see that the content
of $F_2$ is
\begin{equation*}
\cont(F_2)=\gcd\{F_1(a),QF_1'(a),Q^2{F_1''(a)\over2}, \dots,
Q^D{F_1^{(D)}(a)\over D!}\}.
\end{equation*}
Thus for any prime $p$,
\begin{multline}
\ord_p(\cont(F_2)) = \min\big\{ \ord_p(F_1(a)), \ord_p(Q) +
\ord_p(F_1'(a)), \\
2\ord_p(Q) + \ord_p\big( {F_1''(a)\over2} \big), \dots, D\ord_p(Q) +
\ord_p\big( {F_1^{(D)}(a)\over D!} \big) \big\},
\label{contmin}
\end{multline}
where all of the terms of the form $\ord_p(\cdot)$ are nonnegative
since the quantities involved are all integers. By general properties
of the discriminant of a polynomial, we know that if $p^\nu$ divides
both $F_1(a)$ and $F'_1(a)$, then $p^\nu$ divides $\Delta$. Put
another way,
\begin{equation*}
\min\{\ord_p(F_1(a)),\ord_p(F'_1(a))\} \le \ord_p(\Delta) < \ord_p(Q),
\end{equation*}
where the last inequality is one of our hypotheses on $Q$. In particular,
one of the first two terms on the right-hand side of
equation~(\ref{contmin}) is less than $2\ord_p(Q)$, and hence
equation~(\ref{contmin}) can be simplified to
\begin{equation}
\ord_p(\cont(F_2)) = \min\{ \ord_p(F_1(a)), \ord_p(Q) +
\ord_p(F_1'(a))\}.
\label{contmin2}
\end{equation}
{\hfuzz=3pt

If $\ord_p(F_1(a)) \le \ord_p(Q) + \ord_p(F_1'(a))$, then $\cont(F_2)
= bp^{\ord_p(F_1(a))}$ for some integer $b$ that is not divisible by
$p$, and so
\begin{equation*}
F_3(0) = {F_2(0)\over bp^{\ord_p(F_1(a))}} \equiv b^{-1} {F_1(a)\over
p^{\ord_p(F_1(a))}} \not\equiv 0\mod p.
\end{equation*}
On the other hand, if $\ord_p(F_1(a)) > \ord_p(Q) + \ord_p(F_1'(a))$,
then $\cont(F_2) = bp^{\ord_p(Q) + \ord_p(F_1'(a))}$ for some integer
$b$ that is not divisible by $p$, and so
\begin{equation*}
F_3(1) = {F_2(1)\over bp^{\ord_p(Q) + \ord_p(F_1'(a))}} \equiv b^{-1}
{QF_1'(a)\over p^{\ord_p(QF_1'(a))}} \not\equiv0\mod p.
\end{equation*}
In either case we see that $F_3$ takes a nonzero value\mod p, and so
$F_3$ is admissible.

} % end \hfuzz
To show that $F_3$ is exclusive we want to show, given two distinct
irreducible factors $G_3$ and $H_3$ of $F_3$, that $G_3$ and $H_3$
have no common zeros modulo any prime. Note that any irreducible
factor $G_3$ of $F_3$ is the primitivization of an irreducible factor
$G_2$ of $F_2$, by Gauss' lemma on the contents of polynomials with
integer coefficients; and any irreducible factor $G_2$ of $F_2$ has
the form $G_2(t) = G_1(Qt+a)$ for some irreducible factor $G_1$ of
$F_1$. Similarly, there is an irreducible factor $H_1$ of $F_1$ such
that $H_3$ is the primitivization of $H_2(t)=H_1(Qt+a)$. If $p$ is a
prime not dividing $\Res(G_1,H_1)$, then $G_1$ and $H_1$ have no
common zeros\mod p by general properties of the resultant of two
polynomials, whence the same is clearly true for $G_2$ and $H_2$ and
thus for $G_3$ and $H_3$; so it suffices to consider the case where
$p$ divides $\Res(G_1,H_1)$.

Again, by general properties of the resultant of two polynomials, we know
that if $p^\nu$ divides both $G_1(a)$ and $H_1(a)$, then $p^\nu$ divides
$\Res(G_1,H_1)$. Put another way,
\begin{equation*}
\min\{\ord_p(G_1(a)),\ord_p(H_1(a))\} \le \ord_p(\Res(G_1,H_1)).
\end{equation*}
By exchanging the $G$s and $H$s if necessary, we can assume
without loss of generality that
\begin{equation}
\ord_p(G_1(a)) \le \ord_p(\Res(G_1,H_1)) < \ord_p(Q),  \label{analg}
\end{equation}
where the last inequality is another of our hypotheses on $Q$. From
the equation analogous to~(\ref{contmin}) for $G_2$, we see that the
inequality~(\ref{analg}) implies that $\ord_p(\cont(G_2)) =
\ord_p(G_1(a))$. Therefore $\cont(G_2) = bp^{\ord_p(G_1(a))}$ for some
integer $b$ that is not divisible by $p$, and so
\begin{equation*}
G_3(n) = {G_2(n)\over bp^{\ord_p(G_1(a))}} \equiv b^{-1} {G_1(a)\over
p^{\ord_p(G_1(a))}} \not\equiv 0\mod p
\end{equation*}
independent of $n$, so that $G_3$ has no zeros\mod p whatsoever. Hence
certainly $G_3$ and $H_3$ have no common zeros\mod p, which shows that
$F_3$ is exclusive.
\end{proof}

With these two lemmas in hand we can now establish
Proposition~\ref{salvationprop}, which we restate here for the
reader's convenience.

\restate{salvationprop}

\begin{proof}
Let $F$ be an integer-valued polynomial, let $K$ be the number of
distinct irreducible factors of $F$, and let $d_1$, \dots, $d_K$ be
the degrees of these factors. Let $d=\max\{d_1,\dots,d_K\}$, and let
$k$ be the number of distinct irreducible factors of $F$ whose degree
equals $d$. Let $u$ and $U$ be real numbers satisfying $1/d\le u\le
U<(d-1/k)^{-1}$; let $x\ge1$ be a real number and set $y=x^{1/u}$. We
are trying to show that the asymptotic formula~(\ref{mainasymform})
holds. If $d_j\le d-1$ then $d_ju<(d-1)(d-1/k)^{-1}\le1$, in which
case $\rho(d_ju)=1$. Thus we are trying to prove that
\begin{equation}
\Psi(F;x,x^{1/u}) = x\rho(du)^k + O\big( \frac x{\log x} \big)
\label{kasymform}
\end{equation}
holds uniformly for $x\ge1$ and $0<u\le U$, where here and throughout
this proof, all constants implicit in $O$-notation may depend on $F$
and $U$. Notice that $y\ge x^\alpha$ where $\alpha=1/U>d-1/k\ge
d-1$. Therefore, by Lemma~\ref{balefflem} we can find an effective
polynomial $F_1(t)$ with integer coefficients that is the product of
$k$ distinct irreducible polynomials of degree $d$, such that
\begin{equation}
\Psi(F_1;x,y)=\Psi(F;x,y)+O(1).  \label{PsiPsi1}
\end{equation}

Let $Q=Q(F_1)$ be the integer given by equation~(\ref{smallQ}), so
that $Q$ satisfies the hypotheses of Lemma~\ref{salvationlem}. For
each fixed integer $0\le a<Q$, define the polynomial $F_2(a;t) =
F_1(Qt+a)$, and let $F_3(a;t)$ be the primitivization of
$F_2(a;t)$. Each of these polynomials $F_3(a;t)$ has integer
coefficients, and is balanced and effective because $F_1$
is. Moreover, by Lemma~\ref{salvationlem} each $F_3(a;t)$ is
admissible and exclusive as well. By the hypothesis of the proposition
to be proved, we know that Theorem~\ref{mainthm} holds for each
$F_3(a;t)$. Therefore (assuming \hyp) we have
\begin{equation}
\Psi\big( F_3(a;t);x,x^{1/u} \big) = x\rho(du)^k + O\big( {x\over\log x}
\big),
\label{F3asym}
\end{equation}
since each $F_3(a;t)$, like $F_1$, is the product of $k$ distinct
irreducible polynomials of degree $d$.

On the other hand, every integer $1\le n\le x$ is congruent\mod Q to
some integer $0\le a<Q$, and so every value $F_1(n)$ for $1\le n\le x$
corresponds to a value $F_2(a;m)$ for some $0\le m<(x-a)/Q$. Furthermore,
the corresponding value $F_3(a;m)$ simply equals $F_2(a;m)/\!\cont(F_2)$,
a difference which does not affect whether the value is $x^{1/u}$-smooth
as soon as $x$ exceeds $\cont(F_2)^U$. Therefore, when
$x$ is sufficiently large in terms of $F$, we have
\begin{equation}
\begin{split}
\Psi(F_1;x,x^{1/u}) &= \sum_{0\le a<Q} \Psi\big( F_3(a;t);{x-a\over
Q},x^{1/u} \big) \\
&= \sum_{0\le a<Q} \Psi\Big( F_3(a;t);{x-a\over Q},\big( {x-a\over Q}
\big)^{1/u_a} \Big),
\label{PsiPsis}
\end{split}
\end{equation}
where
\begin{equation*}
u_a = {u\log((x-a)/Q)\over\log x} = u + O\big( {1\over\log x} \big).
\end{equation*}
The function $\rho(u)$ satisfies $-1<\rho'(u)<0$ for all $u>1$, and so
$\rho(du_a) = \rho(du) + O(1/\log x)$. Using the asymptotic
formula~(\ref{F3asym}) in equation~(\ref{PsiPsis}), we conclude that
\begin{equation*}
\begin{split}
\Psi(F_1;x,x^{1/u}) &= \sum_{0\le a<Q} \big( {x-a\over Q}\rho(du_a)^k
+ O\big( {x\over\log x} \big) \big) \\
&= \sum_{0\le a<Q} \big( {x-a\over Q}\big ( \rho(du)+O\big(
{1\over\log x} \big) \big)^k \big) + O\big( {x\over\log x} \big) =
x\rho(du)^k + O\big( {x\over\log x} \big).
\end{split}
\end{equation*}
Together with equation~(\ref{PsiPsi1}), this shows that the asymptotic
formula~(\ref{kasymform}) holds for the original polynomial $F$,
which establishes the proposition.
\end{proof}

\section{Two Combinatorial Propositions}\setcounter{equation}{0}\label{trashsec}\noindent In this section we establish the combinatorial
Propositions~\ref{introMprop} and~\ref{alterMprop}, which link the
counting function of smooth values of a polynomial to the counting
functions of prime values of a related family of polynomials, via the
expressions $M(f;x,y)$ defined in equation~(\ref{Mfxydef}). The first
of these two propositions exhibits $\Psi(F;x,y)$ as a combination of
terms of the form $M(f;x,y)$.

\restate{introMprop}

\begin{proof}
This is simply inclusion--exclusion on the factors of $F$ that, for a
given argument $n$, have large prime divisors. For any integer $n$
define
\begin{equation}
X_i(n) = \begin{cases}
1,&\hbox{if there exists a prime $p>y$ dividing $F_i(n)$,}\\
0,&\hbox{otherwise (i.e., if $F_i(n)$ is $y$-smooth).}\end{cases}
\label{Xindef}
\end{equation}
Then for any nonempty subset $S$ of $\{1,\dots,K\}$, the
definition~(\ref{Mfxydef}) of $M(f;x,y)$ for $f=\prod_{i\in S} F_i$ is
equivalent~to
\begin{equation*}
M\bigg( \big( \prod_{i\in S} F_i \big);x,y \bigg) = \sum_{n\le x}
\bigg( \prod_{i\in S} X_i(n) \bigg).
\end{equation*}
But the definition~(\ref{Xindef}) also implies that
\begin{equation*}
\prod_{i=1}^K (1-X_i(n)) = \begin{cases}
1,&\hbox{if $F(n)$ is $y$-smooth,}\\0,&\hbox{otherwise,}
\end{cases}
\end{equation*}
since $F(n)$ is $y$-smooth if and only if each $F_i(n)$ is $y$-smooth.
We therefore find that
\begin{equation*}
\begin{split}
\Psi(F;x,y) = \sum_{n\le x} \prod_{i=1}^K (1-X_i(n)) &= \floor x +
\sum_{n\le x} \sum \begin{Sb}S\subset\{1,\dots,K\} \\
S\ne\emptyset\end{Sb} \prod_{i\in S} (-X_i(n)) \\
&= \floor x + \sum_{1\le k\le K} (-1)^k \sum
\begin{Sb}S\subset\{1,\dots,K\} \\ |S|=k\end{Sb} M\bigg( \big(
\prod_{i\in S} F_i \big);x,y \bigg),
\end{split}
\end{equation*}
which is equivalent to the statement of the proposition.
\end{proof}

The second of these two propositions exhibits $M(f;x,y)$ as a
combination of terms of the form $\pi(f\hhb,\cdot)$. All of the
notation introduced on page~\pageref{defspage} will be assumed for the
remainder of this section. We also need to define, for any positive
integers $\hsdots$ and any element $b$ of $\Rhs$, the
polynomial
\begin{equation}
f\6i\hsb(t) = {f_i(\hdots t+b)\over h_i}  \label{fihsbdef}
\end{equation}
for each $1\le i\le k$. The fact that $h_i$ divides each $f_i(b)$
implies, by the Taylor expansion~(\ref{taylorF}) applied to $f_i$,
that each $f\6i\hsb$ has integer coefficients. Because of this, the
polynomial $f\hhb$ has the natural factorization
$f\hhb=f\61\hsb\dots f\6k\hsb$ into irreducible polynomials.

\restate{alterMprop}

\begin{proof}
We recall the definition~(\ref{Mfxydef}) of $M(f;x,y)$:
\begin{multline*}
M(f;x,y) = \#\{ 1\le n\le x\colon \hbox{for each irreducible factor
$g$ of $f$,} \\
\hbox{there exists a prime $p>y$ such that $p\mid g(n)$} \}.
\end{multline*}
Since $y\ge\sqrt{\xi_i}$ for each $i$ when $x$ is sufficiently large,
there is a unique $k$-tuple $(p_1,\dots,p_k)$ of primes for each
argument~$n$ that is counted by $M(f;x,y)$. If we write the values
$f_i(n)$ as $p_ih_i$ for suitable integers $h_i$, we see that
\begin{equation*}
\begin{split}
M(f;x,y) &= \3{(n,p_1,\dots,p_k,h_1,\dots,h_k)\colon 1\le n\le x,\,
\hbox{each }p_i>y,\, \hbox{each }f_i(n)=p_ih_i} \\
&= \sum_{h_1\le\xi_1/y} \dots \sum_{h_k\le\xi_k/y} \#\{ 1\le n\le
x\colon \\
&\hskip2in \hbox{each }h_i|f_i(n),\, \hbox{each }f_i(n)/h_i\hbox{ is a
prime exceeding }y \}.
\end{split}
\end{equation*}
Because $f$ is exclusive, no prime can divide two different values
$f_i(n)$ at the same argument~$n$; therefore, we may insert the
condition of summation $(h_i,h_j)=1$ $(1\le i<j\le k)$ without
changing the sum. Furthermore, if $h_i$ divides $f_i(n)$ for each $i$,
then $n$ is congruent to some element $b$ of $\Rhs$ by its
definition~(\ref{Rhsdef}), and conversely. Therefore
\begin{multline*}
M(f;x,y) = \mathop{\sum_{h_1\le\xi_1/y} \dots
\sum_{h_k\le\xi_k/y}}_{(h_i,h_j)=1 \; (1\le i<j\le k)} \sum_{b\in\Rhs}
\#\{ 1\le n\le x,\, n\equiv b\mod\hdots\colon \\
\hbox{each }f_i(n)/h_i\hbox{ is a prime exceeding }y \}.
\end{multline*}
Making the change of variables $n=\hdots m+b$, we see that
\begin{equation*}
\begin{split}
M(f;x,y) &= \mathop{\sum_{h_1\le\xi_1/y} \dots
\sum_{h_k\le\xi_k/y}}_{(h_i,h_j)=1 \; (1\le i<j\le k)} \sum_{b\in\Rhs}
\#\big\{ 0\le m\le {x-b\over\hdots}\colon \\
&\hskip2in \hbox{each }f_i(\hdots m+b)/h_i\hbox{ is a prime exceeding
}y \big\} \\
&= \mathop{\sum_{h_1\le\xi_1/y} \dots
\sum_{h_k\le\xi_k/y}}_{(h_i,h_j)=1 \; (1\le i<j\le k)} \sum_{b\in\Rhs}
\#\big\{ 0\le m\le {x-b\over\hdots}\colon \\
&\hskip2in \hbox{each }f\6i\hsb(m)\hbox{ is a prime exceeding }y
\big\}
\end{split}
\end{equation*}
by the definition~(\ref{fihsbdef}) of the $f\6i\hsb$. The polynomial
$f$ is effective and hence strictly increasing for $t\ge0$, so if we
define $\eta_\hsdots$ to be the smallest real number $\eta$ such that
$f\6i\hsb(\eta)\ge y$ for each $i$, then we see that
\begin{equation*}
M(f;x,y) = \mathop{\sum_{h_1\le\xi_1/y} \dots
\sum_{h_k\le\xi_k/y}}_{(h_i,h_j)=1 \; (1\le i<j\le k)} \sum_{b\in\Rhs}
\big( \pi\big( f\hhb;{x-b\over\hdots} \big) - \pi(f\hhb;\eta_\hsdots)
\big).
\end{equation*}
As for the size of $\eta_\hsdots$, we see from the
definition~(\ref{fihsbdef}) of $f\6i\hsb$ that
\begin{equation*}
(f\6i\hsb)^{-1}(y) = {f_i^{-1}(h_iy)-b\over\hdots} \asymp
{(h_iy)^{1/d}\over\hdots}
\end{equation*}
since each $f_i$ has degree $d$, and so
$\eta_\hsdots\asymp(y\max\{\hdots\})^{1/d}(\hdots)^{-1}$. This
establishes the proposition.
\end{proof}

\section{Roots of Polynomials Modulo Integers}\setcounter{equation}{0}\label{multappendixsec}\noindent The object of this section is to establish
Proposition~\ref{cflbprop}. The tools we shall need are relationships
between the number of local roots of polynomials of the form $f\hhb$,
defined in equation~(\ref{fhhbdef}), and the number of local roots of
the polynomial $f$ itself. For the first part of this section, we make
no particular assumption on the polynomial $f$ except that it has
integer coefficients. For any positive integer $h$ we define
\begin{equation}
\R(f;h) = \{ 1\le b\le h\colon h|f(b)\},
\label{Rhdef}
\end{equation}
analogous to the definition~(\ref{Rhsdef}) of $\Rhs$. For any integer
$b$ that is congruent\mod h to some element of $\R(f;h)$, we define
the polynomial $f\hb$ by $f\hb(t) = f(ht+b)/h$, consistent with the
notation $f\hhb$ defined in equation~(\ref{fhhbdef}).

We recall that $\sigma(f;h)$ denotes the number of solutions of
$f(b)\equiv0\mod h$, or simply the cardinality of $\R(f;h)$. Clearly
$\sigma$ is a nonnegative, integer-valued function satisfying
$\sigma(f;h)\le h$ for all positive integers $h$. Also, if $m$ and $n$
are relatively prime positive integers, then by the Chinese remainder
theorem there is a bijection between roots of $f\mod{mn}$ and pairs of
roots of $f\mod m$ and\mod n; this implies that $\sigma(f;h)$ is a
multiplicative function of $h$.

The following four lemmas exhibit simple relationships between values
of $\sigma(f\hb)$ and values of~$\sigma(f)$.

\begin{lemma}
Let $f$ be a polynomial with integer coefficients, let $h$ be a
positive integer, and let $b$ be an element of $\R(f;h)$. If a prime $p$
does not divide $h$, then $\sigma(f\hb;p)=\sigma(f;p)$.
\label{mostcaseslem}
\end{lemma}

\begin{proof}
By definition,
\begin{equation}
\begin{split}
\sigma(f\hb;p) &= \3{a\mod p\colon f\hb(a)\equiv0\mod p} \\
&= \3{a\mod p\colon {f(ha+b)\over h} \equiv0\mod p}.
\end{split}
\label{cantthinkeq}
\end{equation}
Since $p$ does not divide $h$, we may multiply both sides of the
latter congruence by $h$. Then making the bijective change of
variables $a'\equiv ha+b\mod p$, we see that
\begin{equation*}
\sigma(f\hb;p) = \3{a'\mod p\colon f(a')\equiv0\mod p} = \sigma(f;p),
\end{equation*}
as claimed.
\end{proof}

\begin{lemma}
Let $f$ be a polynomial with integer coefficients, let $h$ be a
positive integer, and let $b$ and $b'$ be integers with $f(b)\equiv
f(b')\equiv 0\mod h$. If $b\equiv b'\mod h$ then
$\sigma(f\hb;p)=\sigma(f_{h,b'};p)$.
\label{bblem}
\end{lemma}

\begin{proof}
Write $b=b'+hq$ for some integer $q$. Making the bijective change of
variables $a'\equiv a+q\mod p$ in the latter congruence in
equation~(\ref{cantthinkeq}), we see that
\begin{equation*}
\sigma(f\hb;p) = \3{a'\mod p\colon {f(ha'+b')\over h}\equiv0\mod p} =
\sigma(f_{h,b'};p),
\end{equation*}
as claimed.
\end{proof}

\begin{lemma}
Let $f$ be a polynomial with integer coefficients, let $h$ be a
positive integer, and let $b$ be an element of $\R(f;h)$. If a prime
power $p^\nu$ exactly divides $h$, then
$\sigma(f\hb;p)=\sigma(f_{p^\nu,b};p)$.
\label{pnulem}
\end{lemma}

\begin{proof}
Write $h=p^\nu h'$ for some integer $h'$. Since $p$ does not divide
$h'$, we may multiply both sides of the latter congruence in
equation~(\ref{cantthinkeq}) by $h'$. Then making the bijective change
of variables $a'\equiv h'a\mod p$, we see that
\begin{equation*}
\sigma(f\hb;p) = \3{a'\mod p\colon {f(p^\nu a'+b)\over
p^\nu}\equiv0\mod p} = \sigma(f_{p^\nu,b};p),
\end{equation*}
as claimed.
\end{proof}

\begin{lemma}
Let $f$ be a polynomial with integer coefficients and let $h$ and $n$
be positive integers. Then $\sum_{b\in\R(f;h)} \sigma(f\hb;n) =
\sigma(f;nh)$.
\label{sigsumlem}
\end{lemma}

\begin{proof}
By definition,
\begin{equation}
\begin{split}
\sigma(f\hb;n) &= \3{a\mod n\colon f\hb(a)\equiv0\mod n} \\
&= \3{a\mod n\colon f(h a+b)/h \equiv0\mod n} \\
&= \3{a\mod n\colon f(h a+b)\equiv0\mod{nh}} \\
&= \3{c\mod{nh},\, c\equiv b\mod h\colon f(c)\equiv0\mod{nh}} \\
&= \3{c\in\R(f;nh)\colon c\equiv b\mod h}.
\end{split}
\label{copied}
\end{equation}
Now every root of $f\mod{nh}$ is certainly a root of $f\mod h$, and so
it follows that every $c\in\R(f;nh)$ is congruent to some
$b\in\R(f;h)$. Therefore
\begin{equation*}
\sigma(f;nh) = \3{c\in\R(f;nh)} = \sum_{b\in\R(f;h)} \3{c\in\R(f;nh)\colon
c\equiv b\mod h} = \sum_{b\in\R(f;h)} \sigma(f\hb;n)
\end{equation*}
from the last line of equation~(\ref{copied}), which establishes the
lemma.
\end{proof}

We remark that the last line in equation~(\ref{copied}) shows that the
quantity $\sigma(f\hb;n)$ can be interpreted as the number of
``lifts'' of the root $b$ of $f\mod h$ to roots of $f\mod{nh}$. Also,
since the inequality $\sigma(f;n)\le n$ is trivial for any polynomial
$f$, we see from Lemma~\ref{sigsumlem} that
\begin{equation}
\sigma(f;nh) = \sum_{b\in\R(f;h)} \sigma(f\hb;n) \le
\sum_{b\in\R(f;h)} n = n\sigma(f;h).
\label{usein7}
\end{equation}

Next we recall the definition~(\ref{Gsigdefs}) of the multiplicative
function $\sigma^*(f;n)$,
\begin{equation*}
\sigma^*(f;n) = \prod_{p^\nu\exdiv n} \big( \sigma(f;p^\nu) -
{\sigma(f;p^{\nu+1})\over p} \big),
\end{equation*}
and establish some of its properties. First, we may put $h=p^\nu$ and
$n=p$ in equation~(\ref{usein7}) to see that
\begin{equation}
\sigma^*(f;p^\nu) = \sigma(f;p^\nu) - {\sigma(f;p^{\nu+1})\over p} \ge
\sigma(f;p^\nu) - {p\sigma(f;p^\nu)\over p} = 0.
\label{sigstarnonneg}
\end{equation}
This shows that $\sigma^*$ is always nonnegative. Next we give an
alternate characterization (and in fact the reason for existence) of
the quantity $\sigma^*(f;h)$. The following lemma provides the key
fact needed for the proof of Proposition~\ref{cflbprop}.

\begin{lemma}
For any polynomial $f$ with integer coefficients and any positive
integer $h$,
\begin{equation*}
\sum_{b\in\R(f;h)} \prod_{p\mid h} \big( 1-\frac{\sigma(f_{h;b};p)}p
\big) = \sigma^*(f;h).
\end{equation*}
\label{chineselem}
\end{lemma}

\begin{proof}
Factor $h=p_1^{\nu_1}\dots p_l^{\nu_l}$ into powers of distinct
primes. We first note that by Lemma~\ref{pnulem},
\begin{equation*}
\sum_{b\in\R(f;h)} \prod_{p\mid h} \big( 1-\frac{\sigma(f_{h;b};p)}p \big)
= \sum_{b\in\R(f;h)} \prod_{p^\nu\exdiv h} \big(
1-{\sigma(f_{p^\nu,b};p)\over p} \big) = \sum_{b\in\R(f;h)} \prod_{i=1}^l
\big( 1-{\sigma(f_{p_i^{\nu_i},b};p_i)\over p_i} \big).
\end{equation*}
For every $l$-tuple $(b_1,\dots,b_l)$ such that
each $b_i\in\R(f;p_i^{\nu_i})$, the Chinese Remainder Theorem gives us a
$b\in\R(f;h)$ such that $b\equiv b_i$\mod{p_i^{\nu_i}} for each $1\le
i\le l$, and this correspondence is bijective. Therefore
\begin{equation*}
\sum_{b\in\R(f;h)} \prod_{i=1}^l \big(
1-{\sigma(f_{p_i^{\nu_i},b};p_i)\over p_i} \big) = \sum
\begin{Sb}(b_1,\dots,b_l) \\ b_i\in\R(f;p_i^{\nu_i})\end{Sb} \prod_{i=1}^l
\big( 1-{\sigma(f_{p_i^{\nu_i},b_i};p_i)\over p_i} \big),
\end{equation*}
since $\sigma(f_{p_i^{\nu_i},b};p_i)$ only depends on
$b$\mod{p_i^{\nu_i}} by Lemma~\ref{bblem}. This last sum now factors:
\begin{equation}
\begin{split}
\sum \begin{Sb}(b_1,\dots,b_l) \\ b_i\in\R(f;p_i^{\nu_i})\end{Sb}
\prod_{i=1}^l \big( 1-{\sigma(f_{p_i^{\nu_i},b_i};p_i)\over p_i} \big)
&= \prod_{i=1}^l \sum_{b_i\in\R(f;p_i^{\nu_i})} \big(
1-{\sigma(f_{p_i^{\nu_i},b_i};p_i)\over p_i} \big) \\
&= \prod_{i=1}^l \bigg( \sum_{b_i\in\R(f;p_i^{\nu_i})} 1\; - {1\over p_i}
\sum_{b_i\in\R(f;p_i^{\nu_i})} \sigma(f_{p_i^{\nu_i},b_i};p_i) \bigg) \\
&= \prod_{i=1}^l \big( \sigma(f;p_i^{\nu_i}) -
{\sigma(f;p_i^{\nu_i+1})\over p_i} \big)
\end{split}
\label{yikes}
\end{equation}
by the definition of $\sigma$ and by Lemma~\ref{sigsumlem} with
$h=p_i^{\nu_i}$ and $n=p_i$. Since the last expression of
equation~(\ref{yikes}) is just $\prod_{i=1}^l \sigma^*(f;p_i^{\nu_i})
= \sigma^*(f;h)$, the lemma is established.
\end{proof}

For the remainder of this section, we specialize to the case where $f$
is a primitive, balanced polynomial with integer coefficients, and we
let $k$ denote the number of distinct irreducible factors of $f$. We
also recall the definitions (\ref{liFxdef}), (\ref{Rhsdef}),
and~(\ref{lihsFxdef}) of $\li(F;x)$, $\Rhs$, and $\li_\hsdots(f;x)$,
respectively.

\begin{lemma}
Let $f$ be a polynomial that is effective, primitive, and balanced.
For any positive integers $\hsdots$ and any $b\in\Rhs$, we have
\begin{equation*}
\li\big( f\hhb;{x-b\over h_1\dots h_k} \big) =
{\li_\hsdots(f;x)\over\hdots} + O_f(1)
\end{equation*}
uniformly for $x\ge1$.
\label{lichangelem}
\end{lemma}

\begin{proof}
From the definition~(\ref{liFxdef}) of $\li(F,x)$, we see that
\begin{equation*}
\begin{split}
\li\big( f\hhb;{x-b\over h_1\dots h_k} \big) &=
\int\limits\begin{Sb}0<t\le(x-b)/h_1\dots h_k \\
\min\{f\61\hsb(t),\dots,f\6k\hsb(t)\}\ge2\end{Sb} {dt\over\log
f\61\hsb(t)\dots\log f\6k\hsb(t)} \\
&= {1\over h_1\dots h_k} \int\limits\begin{Sb}b<v\le x \\
\min\{f_1(v)/h_1,\dots,f_k(v)/h_k\}\ge2\end{Sb}
{dv\over\log(f_1(v)/h_1)\dots\log(f_k(v)/h_k)}
\end{split}
\end{equation*}
by the change of variables $v=h_1\dots h_kt+b$. We may change the
lower limit of integration from $b$ to 0, incurring an error that is
$\le b(\log 2)^{-k}\ll_f b\le\hdots$. Therefore
\begin{equation*}
\li\big( f\hhb;{x-b\over h_1\dots h_k} \big) =
{\li_\hsdots(f;x)\over\hdots} + O_f(1)
\end{equation*}
as claimed.
\end{proof}

We have one more lemma to establish before we can prove
Proposition~\ref{cflbprop}, and for this lemma we must define one more
piece of notation. To avoid double subscripts such as $(f_i)\hb$, we
define the polynomial $f\6i\hb$ by $f\6i\hb(t) = f_i(ht+b)/h$, which
as before has integer coefficients if $h$ divides $f_i(b)$; this is to
be distinguished from the polynomial $f\6i\hsb$, which still
corresponds to its definition~(\ref{fihsbdef}).

\begin{lemma}
Let $f$ be a polynomial with integer coefficients that is primitive,
balanced, and exclusive; let $\hsdots$ be pairwise coprime positive
integers and let $b$ be an element of $\Rhs$. If a prime $p$ divides
$h_i$ for some $1\le i\le k$, then
$\sigma(f\hhb;p)=\sigma(f\6i_{h_i,b};p)$.
\label{othercaseslem}
\end{lemma}

\begin{proof}
By the definition of $\sigma$, and using the remark following
equation~(\ref{fihsbdef}) to factor the polynomial $f\hhb$, we have
\begin{equation}
\begin{split}
\sigma(f\hhb;p) &= \3{a\mod p\colon f\hhb(a)\equiv0\mod p} \\
&= \3{a\mod p\colon f\61\hsb(a)\dots f\6k\hsb(a)\equiv0\mod p} \\
&= \3{a\mod p\colon \prod_{j=1}^k {f_j(\hdots a+b)\over
h_j} \equiv0\mod p}.
\end{split}
\label{fromlattereq}
\end{equation}
Since the $h_j$ are pairwise coprime, $p$ does not divide any of the
$h_j$ other than $h_i$. Also, since $b$ is a root of $f_i\mod p$ and
$f$ is exclusive, $b$ cannot be a root of any other $f_j\mod
p$. Therefore, since $h_i\equiv0\mod p$, we see that $f_j(\hdots
a+b)\equiv f_j(b)\not\equiv0\mod p$ for all $j\ne i$.

We may therefore divide the congruence in
equation~(\ref{fromlattereq}) by $f_j(\hdots a+b)/h_j$ for each $j\ne
i$, obtaining
\begin{equation*}
\begin{split}
\sigma(f\hhb;p) &= \3{a\mod p\colon {f(h_ia+b)\over h_i}\equiv0\mod
p} \\
&= \3{a\mod p\colon f\6i_{h_i,b}(a)\equiv0\mod p} =
\sigma(f\6i_{h_i,b};p),
\end{split}
\end{equation*}
which establishes the lemma.
\end{proof}

Armed with these several lemmas, we are now ready to establish:

\restate{cflbprop}

\begin{proof}
By Lemma~\ref{lichangelem} we can write
\begin{multline*}
\sum_{b\in\Rhs} C(f\hhb) \li\big( f\hhb; {x-b\over\hdots} \big) \\
= \big( {\li_\hsdots(f;x)\over\hdots} + O_f(1) \big) \sum_{b\in\Rhs}
C(f\hhb).
\end{multline*}
Thus to establish the lemma, it suffices to show that
\begin{equation}
\sum_{b\in\Rhs} C(f\hhb) = C(f) G(f;h_1\dots
h_k)\sigma^*(f_1;h_1)\dots\sigma^*(f_k;h_k).
\label{withoutli}
\end{equation}

If $p$ is a prime that does not divide $h_1\dots h_k$, then we know
that $\sigma(f\hsb;p)=\sigma(f;p)$ by Lemma~\ref{mostcaseslem}. Thus
\begin{equation*}
\begin{split}
\sum_{b\in\Rhs} C(f\hhb) &= \sum_{b\in\Rhs} \prod_p \big( 1-\frac1p
\big)^{-k} \big( 1-\frac{\sigma(f\hhb;p)}p \big) \\
&= \prod_p \big( 1-\frac1p \big)^{-k} \big( 1-\frac{\sigma(f;p)}p \big)
\prod_{p\mid h_1\dots h_k} \big( 1-\frac{\sigma(f;p)}p \big)^{-1} \\
&\qquad \times \sum_{b\in\Rhs} \prod_{p\mid h_1\dots h_k} \big(
1-\frac{\sigma(f\hhb;p)}p \big) \\
&= C(f)G(f;h_1\dots h_k) \sum_{b\in\Rhs} \prod_{i=1}^k
\prod_{p\mid h_i} \big( 1-\frac{\sigma(f\6i_{h_i,b};p)}p \big),
\end{split}
\end{equation*}
where we have used Lemma \ref{othercaseslem} to change $\sigma(f\hhb)$
to the $\sigma(f\6i_{h_i,b})$ in the last equality. All that remains
to establish equation~(\ref{withoutli}), and hence the proposition, is
to show that
\begin{equation}
\sum_{b\in\Rhs} \prod_{i=1}^k \prod_{p\mid h_i} \big(
1-\frac{\sigma(f\6i_{h_i,b};p)}p \big) =
\sigma^*(f_1;h_1)\dots\sigma^*(f_k;h_k).
\label{withoutlimbo}
\end{equation}

For every $k$-tuple $(b_1,\dots,b_k)$ such that each $b_i\in\R(f_i;h_i)$,
the Chinese remainder theorem gives us a $b\in\Rhs$ such that $b\equiv
b_i$\mod{h_i} for each $1\le i\le k$, and this correspondence is
bijective. Therefore
\begin{equation*}
\sum_{b\in\Rhs} \prod_{i=1}^k \prod_{p\mid h_i} \big(
1-\frac{\sigma(f\6i_{h_i,b};p)}p \big) = \sum
\begin{Sb}(b_1,\dots,b_k) \\ b_i\in\R(f_i;h_i)\end{Sb} \prod_{i=1}^k
\prod_{p\mid h_i} \big( 1-\frac{\sigma(f\6i_{h_i,b_i};p)}p
\big),
\end{equation*}
since for a given prime $p$ the quantity $\sigma(f\6i_{h_i,b_i};p)$
depends only on $b\mod{h_i}$ by Lemma~\ref{bblem}. This last sum now
factors as
\begin{equation*}
\sum \begin{Sb}(b_1,\dots,b_k) \\ b_i\in\R(f_i;h_i)\end{Sb} \prod_{i=1}^k
\prod_{p\mid h_i} \big( 1-\frac{\sigma(f\6i_{h_i,b_i};p)}p \big) =
\prod_{i=1}^k \sum_{b_i\in\R(f_i;h_i)} \prod_{p\mid h_i} \big(
1-\frac{\sigma(f\6i_{h_i,b_i};p)}p \big) = \prod_{i=1}^k
\sigma^*(f_i;h_i)
\end{equation*}
by Lemma~\ref{chineselem}. This establishes
equation~(\ref{withoutlimbo}) and hence the proposition.
\end{proof}

\section{Multiplicative Functions of One Variable}\setcounter{equation}{0}\label{1multsec}\noindent Two propositions from Section~\ref{rgosec} are established in
this section. Proposition~\ref{stupidOprop} is rather easy, as it
follows from crude upper bounds for the multiplicative functions
$G(f;n)$ and $\sigma(f;n)$ defined in equation~(\ref{Gsigdefs}).
Proposition~\ref{usehypprop}, on the other hand, requires an accurate
knowledge of the order of magnitude of the summatory function of
$G(f;n)\sigma(f;n)$. For this we use a mean-value result (see
Proposition~\ref{mvabprop} in Appendix~A) for a general sum $\sum_{n\le
x}g(n)$ of a multiplicative function $g$. To apply such a
mean-value result, some conditions on the values of the multiplicative
function in question must be verified; the next three lemmas provide
simple estimates of this type for the multiplicative functions
$\sigma(f;n)$, $G(f;n)$, and $\sigma^*(f;n)$. In these three lemmas, all
constants implicit in the $\ll$ and $O$-notations may depend on the
polynomial $f$.

\begin{lemma}
If $f$ is a squarefree polynomial with integer coefficients, then
$\sigma(f;p^\nu)\ll1$ uniformly for all prime powers $p^\nu$.
\label{siglem}
\end{lemma}

\begin{proof}
Let $\Delta$ be the discriminant of $f$, which is nonzero since $f$ is
squarefree, and write $\Delta=\prod_p p^{\theta(p)}$ where all but
finitely many of the $\theta(p)$ are zero. Huxley~\cite{Hux:ANoPC}
gives a bound for $\sigma$ that implies
\begin{equation*}
\sigma(f;p^\nu)\le (\deg f)p^{\theta(p)/2}
\end{equation*}
for any squarefree polynomial $f$ and any prime power $p^\nu$ (this
estimate is improved by Stewart \cite{Ste:OtNoSoPCaTE}, though it will
suffice for our purposes as stated). In particular we see that
$\sigma(f;p^\nu)\le(\deg f)\Delta^{1/2}$ for all prime powers $p^\nu$,
which establishes the lemma.
\end{proof}

\begin{lemma}
If $f$ is a polynomial with integer coefficients that is squarefree
and admissible, then $G(f;p^\nu)=1+O(1/p)$ and
$\sigma^*(f;p^\nu)=\sigma(f;p^\nu)+O(1/p)$ uniformly for all prime
powers~$p^\nu$.
\label{Olem}
\end{lemma}

\begin{proof}
From general facts about polynomials over finite fields, for a given
prime $p$ either $\sigma(f;p)=p$ or $\sigma(f;p)\le\deg f$. However,
the fact that $f$ is admissible means $\sigma(f;p)<p$, and so
$\sigma(f;p)\le\deg f$ for all primes $p$. This implies that for all
primes $p\ge2\deg f$,
\begin{equation*}
G(f;p^\nu)-1 = \big( 1-\frac{\sigma(f;p)}p \big)^{-1}-1 =
{\sigma(f;p)\over p-\sigma(f;p)} \le {\deg f\over p/2} \ll \frac1p.
\end{equation*}
This shows that $G(f;p^\nu)=1+O(1/p)$ for all prime powers $p$ (by
adjusting the implicit constant if necessary). Similarly,
\begin{equation*}
\sigma(f;p^\nu) - \sigma^*(f;p^\nu) = {\sigma(f;p^{\nu+1})\over p} \ll
\frac1p
\end{equation*}
by Lemma~\ref{siglem}, showing that
$\sigma^*(f;p^\nu)=\sigma(f;p^\nu)+O(1/p)$.
\end{proof}

\begin{lemma}
If $f$ is a polynomial with integer coefficients that is squarefree
and admissible, then for any $\ep>0$, we have $\sigma(f;n)\ll_{\ep}
n^\ep$ and $G(f;n)\ll_{\ep} n^\ep$.
\label{eplem}
\end{lemma}

\begin{proof}
It is well-known that the number $\omega(n)$ of distinct prime
divisors of $n$ satisfies $\omega(n)\ll\log n/\log\log n$, which
implies that $A^{\omega(n)}\ll_{A,\ep}n^\ep$ for any positive
constants $A$ and~$\ep$. Therefore, any multiplicative function $g(n)$
satisfying $|g(p^\nu)|\le A$ for all prime powers $p^\nu$
automatically satisfies $g(n)\ll_{A,\ep} n^\ep$. By
Lemmas~\ref{siglem} and~\ref{Olem}, respectively, both $\sigma$ and
$G$ have this property for some constant $A$ depending on the
polynomial $f$, and so the lemma is established.
\end{proof}

The crude upper bounds given in Lemma~\ref{eplem} are enough to
establish Proposition~\ref{stupidOprop}. For this proposition, as well
as for Proposition~\ref{usehypprop} which is proved at the end of this
section, we recall from page~\pageref{defspage} the definitions of the
parameters $d$, $k$, $u$, $U$, $\xi_i$, and~$y$.

\restate{stupidOprop}

\begin{proof}
If we include in the sum those (nonnegative) terms for which the $h_i$
are not necessarily pairwise coprime, and use the trivial inequality
$\sigma^*(f;n)\le\sigma(f;n)$, we see that
\begin{multline*}
\mathop{\sum_{h_1\le\xi_1/y}\dots\sum_{h_k\le\xi_1/y}}_{(h_i,h_j)=1 \;
(1\le i<j\le k)} G(f;\hdots)\sigma^*(f_1;h_1) \dots \sigma^*(f_k;h_k) \\
\le \sum_{h_1\le\xi_1/y}\dots\sum_{h_k\le\xi_k/y}
G(f;\hdots)\sigma(f_1;h_1) \dots \sigma(f_k;h_k).
\end{multline*}
Given $\ep>0$, we see from Lemma~\ref{eplem} that
\begin{equation*}
\begin{split}
\sum_{h_1\le\xi_1/y}\dots\sum_{h_k\le\xi_k/y}
G(f;\hdots)\sigma(f_1;h_1) \dots \sigma(f_k;h_k) &\ll_{f,\ep}
\sum_{h_1\le\xi_1/y}\dots\sum_{h_k\le\xi_k/y} (\hdots)^{2\ep} \\
&\ll_{f,\ep} \big( {\xi_1\dots\xi_k\over y^k} \big)^{1+2\ep} \le
(x^{k(d-1/U)})^{1+2\ep}.
\end{split}
\end{equation*}
Since $k(d-1/U)<1$ by the upper bound~(\ref{Urange}) on $U$, we can
choose $\ep$ so small (depending on $f$ and $U$) that the right-hand
side is $\ll_{f,U} x/\log x$, which establishes the proposition.
\end{proof}

When applying Proposition~\ref{mvabprop} to a particular multiplicative
function, it is of course necessary to verify the
hypothesis~(\ref{hypothesis}) for that function. When the function in
question is $\sigma(f;n)$ for some polynomial $f$, the relevant
asymptotic formula is well-known. If $f$ is a polynomial with
integer coefficients with $k$ distinct irreducible factors, then the
values taken by $\sigma$ on primes are $k$ on average; more precisely,
Nagel \cite{Nag:GdTdT} showed that the asymptotic formula
\begin{equation}
\sum_{p\le w} \frac{\sigma(f;p)\log p}p = k\log w + O_f(1)
\label{nagel}
\end{equation}
holds for all $w\ge2$. For the purposes of establishing
Proposition~\ref{usehypprop} (and, in the next section,
Proposition~\ref{brunprop}), we need to verify this hypothesis for the
slightly different function $G(f;n)\sigma(f;n)$, while for
Proposition~\ref{lipsprop}, the appropriate function is
$G(f;n)\sigma^*(f;n)$. These verifications are the subject of the
following lemma.

\begin{lemma}
Let $f$ be a squarefree polynomial with integer coefficients, and let
$f_i$ denote any one of the irreducible factors of $f$. Then
\begin{equation*}
\sum_{p\le w} {G(f;p)\sigma(f_i;p)\log p\over p} = \log w + O_f(1),
\end{equation*}
and the same is true if $\sigma$ is replaced by $\sigma^*$.
\label{kappaverifylem}
\end{lemma}

\begin{proof}
The constants implicit in the $\ll$ and $O$-notations in this proof
may depend on the polynomial $f$ and thus on $f_i$ as well. By
Lemma~\ref{Olem} we know that
\begin{equation*}
G(f;p)\sigma(f_i;p) = \sigma(f_i;p)(1+O(1/p)) = \sigma(f_i;p) + O(1/p)
\end{equation*}
using Lemma~\ref{siglem}, and so
\begin{equation*}
\sum_{p\le w} {G(f;p)\sigma(f_i;p)\log p\over p} = \sum_{p\le w}
{\sigma(f_i;p)\log p\over p} + O\bigg( \sum_p {\log p\over p^2} \bigg)
= \log w + O(1)
\end{equation*}
by the asymptotic formula~(\ref{nagel}) and the fact that the last sum
is convergent.

Again by Lemmas~\ref{siglem} and~\ref{Olem}, we see that
$G(f;p)\sigma^*(f_i;p) = \sigma(f_i;p)+O(1/p)$, and so
\begin{equation*}
\sum_{p\le w} {G(f;p)\sigma^*(f_i;p)\log p\over p} = \sum_{p\le w}
{\sigma(f_i;p)\log p\over p} + O\bigg( \sum_p {\log p\over p^2}
\bigg) = \log w + O(1)
\end{equation*}
as before. This establishes the lemma.
\end{proof}

We may now establish an upper bound of the correct order of magnitude
for the summatory functions that will arise in the proofs of
Proposition~\ref{usehypprop} and (in the next section)
Proposition~\ref{brunprop}.

\begin{lemma}
Let $f$ be a polynomial with integer coefficients that is squarefree
and admissible, and let $f_i$ be any one of the irreducible factors of
$f$. Then
\begin{equation*}
\sum_{n\le x} {G(f;n)\sigma(f_i;n)\over n} \ll_f \log x
\quad\hbox{and}\quad \sum_{n\le x} {G(f;n)\sigma(f_i;n)\over
n^{1-\beta}} \ll_{f,\beta} x^\beta
\end{equation*}
for any $\beta>0$.
\label{oneDlem}
\end{lemma}

\begin{proof}
We want to apply Proposition~\ref{mvabprop} to establish the first claim
of the lemma. The function $G(f;n)\sigma(f_i;n)$ is nonnegative and
multiplicative, and by Lemma~\ref{eplem} we know that
$G(f;n)\sigma(f_i;n)\ll_f n^{1/4}$, say. Furthermore,
Lemma~\ref{kappaverifylem} verifies the condition~(\ref{hypothesis}) with
$\kappa=1$. Therefore we may apply Proposition~\ref{mvabprop}(a) to see
that
\begin{equation*}
\sum_{n\le x} {G(f;n)\sigma(f_i;n)\over n} =
c(G(f;n)\sigma(f_i;n))\log x + O_f(1).
\end{equation*}
This establishes the first claim of the lemma; the second claim
follows easily from the first by a simple partial summation argument.
\end{proof}

We are now ready to establish Proposition~\ref{usehypprop}, showing
that \hyp\ implies an upper bound for the expression $H(f;x,y)$, which
is a sum of error terms of the type $E(f_i;x_i)$. Again we recall from
page~\pageref{defspage} the definitions of $d$, $k$, $u$, $U$,
$\xi_i$, and~$y$.

\restate{usehypprop}

\begin{proof}
All constants implicit in the $\ll$ and $O$-notations in this proof
may depend on the polynomial $f$ and the parameter $U$. We recall the
definition~(\ref{Hfxydef}) of $H(f;x,y)$:
\begin{equation}
H(f;x,y) = \mathop{\sum_{h_1\le\xi_1/y}\dots
\sum_{h_k\le\xi_k/y}}_{(h_i,h_j)=1\, (1\le i<j\le k)} \sum_{b\in\Rhs}
E\big( f\hhb;{x-b\over\hdots} \big).
\label{needaB}
\end{equation}
Notice that \begin{equation*}
h_1\dots h_k \le \xi_1\dots\xi_k/y^k \ll x^{kd}(x^{1/u})^{-d} \le
x^{k(d-1/U)}
\end{equation*}
for $u\le U$; in particular, when $U$ satisfies the bound~(\ref{Urange}),
the exponent $k(d-1/U)$ is strictly less than~1. Consequently the
coefficients of the polynomial $f\hhb$ are certainly $\ll(\hdots)^{d-1}\ll
x^d$ in size, while we have the lower bound
\begin{equation}
{x-b\over\hdots} \gg x^{1-k(d-1/U)}.  \label{xoverhdots}
\end{equation}
We may therefore apply \hyp\ with $B=d/(1-k(d-1/U))>0$ to the error terms
$E(f\hhb;(x-b)/\hdots)$ in equation~(\ref{needaB}); this yields
\begin{equation*}
H(f;x,y) \ll \mathop{\sum_{h_1\le\xi_1/y}\dots
\sum_{h_k\le\xi_k/y}}_{(h_i,h_j)=1\, (1\le i<j\le k)} \sum_{b\in\Rhs}
\big( {C(f\hhb)(x-b)/\hdots \over (\log((x-b)/\hdots))^{k+1}} +1 \big)
\end{equation*}
using \hyp. We see from the lower bound~(\ref{xoverhdots}) that the
logarithmic term in the denominator can be replaced by $\log x$,
yielding
\begin{equation*}
\begin{split}
H(f;x,y) &\ll {x\over\log^{k+1}x} \mathop{\sum_{h_1\le\xi_1/y}\dots
\sum_{h_k\le\xi_k/y}}_{(h_i,h_j)=1\, (1\le i<j\le k)} \frac1\hdots
\sum_{b\in\Rhs} C(f\hhb) \\
&\qquad+ \mathop{\sum_{h_1\le\xi_1/y}\dots
\sum_{h_k\le\xi_k/y}}_{(h_i,h_j)=1\, (1\le i<j\le k)} \sum_{b\in\Rhs}
1;
\end{split}
\end{equation*}
and by equation~(\ref{withoutli}) and the definition of $\Rhs$, this is
the same as
\begin{equation*}
\begin{split}
H(f;x,y) &\ll {C(f)x\over\log^{k+1}x} \mathop{\sum_{h_1\le\xi_1/y}\dots
\sum_{h_k\le\xi_k/y}}_{(h_i,h_j)=1\, (1\le i<j\le k)}
\frac{G(f;\hdots) \sigma^*(f_1;h_1)\dots\sigma^*(f_k;h_k)}\hdots \\
&\qquad+ \mathop{\sum_{h_1\le\xi_1/y}\dots
\sum_{h_k\le\xi_k/y}}_{(h_i,h_j)=1\, (1\le i<j\le k)}
\sigma(f_1;h_1)\dots\sigma(f_k;h_k).
\end{split}
\end{equation*}

Since the $h_i$ are pairwise coprime, we may factor the multiplicative
function $G(f;\hdots)$ into $G(f;h_1)\dots G(f;h_k)$ in the first sum
on the right-hand side of this last equation. By then deleting the
restrictions $(h_i,h_j)=1$ from the conditions of summation in both
sums, and noting that $1\le G(f;h_i)$ and
$\sigma^*(f_i,h_i)\le\sigma(f_i;h_i)$, we see that
\begin{equation}
\begin{split}
H(f;x,y) &\ll {x\over\log^{k+1}x} \prod_{i=1}^k \sum_{h_i\le\xi_i/y}
\frac{G(f;h_i) \sigma(f_i;h_i)}{h_i} + \prod_{i=1}^k
\sum_{h_i\le\xi_i/y} G(f;h_i)\sigma(f_i;h_i) \\
&\ll {x\over\log^{k+1}x} \prod_{i=1}^k \log{\xi_i\over y} +
\prod_{i=1}^k {\xi_i\over y}
\end{split}
\label{gotthere}
\end{equation}
by Lemma~\ref{oneDlem} (applied with $\beta=1$ in the second
sum). Clearly each $\log\xi/y\ll\log x$, while
\begin{equation*}
\prod_{i=1}^k {\xi_i\over y} \ll x^{k(d-1/u)} \ll {x\over\log x},
\end{equation*}
again since $k(d-1/u)<1$ by the upper bound~(\ref{Urange}) on
$U$. Therefore the estimate~(\ref{gotthere}) becomes $H(f;x,y) \ll
x/\log x$, which establishes the proposition.
\end{proof}

\section{Prime Values of Polynomials}\setcounter{equation}{0}\label{brunsec}\noindent In this section we establish Proposition~\ref{brunprop}, by using
Brun's upper bound sieve method to estimate the numbers of prime values
of the polynomials $f\hhb$. To apply Brun's sieve, we must verify some
conditions on the number $\sigma(f\hb;p)$ of local roots of the
polynomials $f\hb$; this is the subject of the following two lemmas.

\begin{lemma}
Let $f$ be a squarefree polynomial with integer coefficients, let $h$
be a positive integer, and let $b$ be an element of $\R(f;h)$. Then
$\sigma(f\hb;p)\ll_f1$ uniformly for all primes $p$ (where the
implicit constant does not depend on $h$ or $b$).
\label{omegachecklem1}
\end{lemma}

\begin{proof}
If $p$ does not divide $h$, then $\sigma(f\hb;p)=\sigma(f;p)$ by
Lemma~\ref{mostcaseslem}. On the other hand, if $p^\nu$ exactly
divides $h$, then $\sigma(f\hb;p)=\sigma(f_{p^\nu,b};p)$ by
Lemma~\ref{pnulem}; moreover, Lemma~\ref{sigsumlem} applied with
$h=p^\nu$ and $n=p$ certainly implies that
$\sigma(f_{p^\nu,b};p)\le\sigma(f;p^{\nu+1})$. In either case, we have
$\sigma(f\hb;p)\le\sigma(f;p^\alpha)$ for some positive integer
$\alpha$, and Lemma~\ref{siglem} tells us that
$\sigma(f;p^\alpha)\ll_f1$. This establishes the lemma.
\end{proof}

\begin{lemma}
Let $f$ be a squarefree polynomial with integer coefficients, let $h$
be a positive integer, and let $b$ be an element of $\R(f;h)$ such
that the polynomial $f\hb$ is admissible. Then
\begin{equation}
\big( 1-{\sigma(f\hb;d)\over d} \big)^{-1}\ll_f1  \label{ifadmissible}
\end{equation}
uniformly for all positive integers $d$ (where the implicit constant
does not depend on $h$ or $b$).
\label{omegachecklem2}
\end{lemma}

\begin{proof}
If $p^\nu$ is any prime power, then putting $h=p$ and $n=p^{\nu-1}$ in
the inequality~(\ref{usein7}), we see that $\sigma(f\hb;p^\nu)\le
p^{\nu-1}\sigma(f\hb;p)$. This implies that
\begin{equation*}
{\sigma(f\hb;d)\over d} = \prod_{p^\nu\exdiv d}
{\sigma(f\hb;p^\nu)\over p^\nu} \le \prod_{p\mid d}
{\sigma(f\hb;p)\over p}.
\end{equation*}
Since each $\sigma(f\hb;p)/p\le1$, it follows that for any prime $p$
dividing $d$, we have $\sigma(f\hb;d)/d \le \sigma(f\hb;p)/p$ and hence
\begin{equation*}
\big( 1-{\sigma(f\hb;d)\over d} \big)^{-1} \le \big(
1-{\sigma(f\hb;p)\over p} \big)^{-1}.
\end{equation*}
Therefore to establish the lemma, it suffices to establish the upper
bound~(\ref{ifadmissible}) for primes. But $p-\sigma(f\hb;p)$ is a
nonnegative integer, and it cannot be 0 since we are assuming that $f\hb$
is admissible. Therefore $p-\sigma(f\hb;p)\ge1$, so we can write
\begin{equation*}
\big( 1-\frac{\sigma(f\hb;p)}p \big)^{-1} = 1 + {\sigma(f\hb;p)\over
p-\sigma(f\hb;p)} \le 1+\sigma(f\hb;p) \ll_f1
\end{equation*}
by Lemma~\ref{omegachecklem1}. This establishes the lemma.
\end{proof}

In the next lemma, we apply Brun's sieve to obtain the desired upper
bound for the number of prime values of the polynomials $f\hhb$, with
the dependence on the parameters $\hsdots$ and $b$ made explicit.

\begin{lemma}
Let $f$ be a polynomial with integer coefficients that is squarefree,
effective, and primitive, and let $k$ denote the number of distinct
irreducible factors of $f$. For any positive integers $\hsdots$ and
any $b\in\Rhs$, we have
\begin{equation}
\pi(f\hhb;t) \ll_f {t\over\log^kt}G(f;\hdots)  \label{brunresult}
\end{equation}
(where the implicit constant does not depend on $\hsdots$, or~$b$).
\label{brunlem}
\end{lemma}

\begin{proof}
If the polynomial $f\hhb$ is not admissible then $\pi(f\hhb;t)\ll_f1$
(see the remarks following the statement of \hyp\ on
page~\pageref{hyppage}). Clearly $G(f;\hsdots)\ge1$ from its
definition, and so the upper bound~(\ref{brunresult}) holds easily in
this case. Thus for the remainder of the proof, we may assume that
$f\hhb$ is in fact admissible.

We first bound $\pi(f\hhb;t)$ by noting that
\begin{equation}
\begin{split}
\pi(f\hhb;t) &= \3{n\le t\colon \hbox{each }f\6i\hsb(n)\hbox{ is
prime}} \\
&\le \3{n\le t\colon \hbox{for each $i$, }p\mid f\6i\hsb(n) \implies
p>(f\6i\hsb(n))^{1/2}} \\
&\le \3{n\le t\colon p\mid f\hhb(n) \implies p>t^{1/3}} + O(t^{2/3}),
\end{split}
\label{tobrun}
\end{equation}
since $f\6i\hsb(n)>n$ for each $i$. This casts the problem of
estimating $\pi(f\hhb;t)$ in terms of a sieving problem.

The version of Brun's sieve that we now employ can be found in
Halberstam--Richert \cite[Theorem~2.2]{HalRic:SM}. Given a set $\cal
A$ of integers, assume that a real number $X$ and a multiplicative
function $\omega(d)$ can be chosen such that
\begin{equation}
\3{a\in{\cal A}\colon d\mid a} = {\omega(d)\over d}X + R_d,
\quad\hbox{where }|R_d|\le\omega(d).
\label{Rdcond}
\end{equation}
Assume that $\omega$ satisfies the two conditions
\begin{equation}
0\le\omega(p)\ll1 \quad\hbox{and}\quad \big( 1-{\omega(d)\over d}
\big)^{-1} \ll1.
\label{omegaconds}
\end{equation}
Then we have the upper bound
\begin{equation}
\3{a\in{\cal A}\colon p\mid a\implies p>z} \ll X\prod_{p\le z} \big(
1-{\omega(p)\over p} \big)
\label{brunbd}
\end{equation}
uniformly for $2\le z\le X$, where the implicit constant depends only
on the implicit constants in the conditions~(\ref{omegaconds}) on
$\omega$.

For our application, we take ${\cal A}=\{f\hhb(n)\colon 1\le n\le
t\}$ and $z=t^{1/3}$; with these choices, the left-hand side of
the bound~(\ref{brunbd}) is exactly $\3{1\le n\le t \colon p\mid f\hhb(n)
\implies p>t^{1/3}}$. We also set $X=t$ and $\omega(d)=\sigma(f\hhb;d)$.
Of course the condition~(\ref{Rdcond}) is satisfied, since the polynomial
$f\hhb$ has $\sigma(f\hhb;d)$ zeros\mod d in every block of $d$
consecutive integers. On the other hand, since we are assuming that
$f\hhb$ is admissible, Lemmas~\ref{omegachecklem1}
and~\ref{omegachecklem2} show that the conditions~(\ref{omegaconds}) are
satisfied. The bound~(\ref{brunbd}) thus implies that
\begin{multline}
\3{1\le n\le t \colon p\mid f\hhb(n) \implies p>t^{1/3}} \ll_f
t\prod_{p\le t^{1/3}} \big( 1-{\sigma(f\hhb;p)\over p} \big) \\
\ll_f {t\over\log^kt}\prod_{p\le t^{1/3}} \big( 1-\frac1p \big)^{-k}
\big( 1-{\sigma(f;p)\over p} \big) \prod \begin{Sb}p\le t^{1/3} \\
p\mid\hdots\end{Sb} \big( 1-{\sigma(f;p)\over p} \big)^{-1} \big(
1-{\sigma(f\hhb;p)\over p} \big)
\label{brunfinish}
\end{multline}
using Lemma~\ref{mostcaseslem}. The first product converges to $C(f)$
as $t$ tends to infinity, while the factors $1-\sigma(f\hhb;p)/p$ in
the second product can be deleted for the purpose of finding an upper
bound. Therefore
\begin{equation*}
\begin{split}
\3{n\le t \colon p\mid f\hhb(n) \implies p>t^{1/3}} &\ll_f
{t\over\log^kt} (C(f)+o(1)) \prod \begin{Sb}p\le t^{1/3} \\
p\mid\hdots\end{Sb} \big( 1-{\sigma(f;p)\over p} \big)^{-1} \\
&\ll_f {t\over\log^kt}G(f;\hdots).
\end{split}
\end{equation*}
In light of the inequality~(\ref{tobrun}), this establishes the lemma.
\end{proof}

The last lemma we give before proceeding to the proof of
Proposition~\ref{brunprop} is an elementary lemma concerning the
behavior of the quantity $\eta_\hsdots$.

\begin{lemma}
Let $\hsdots$ be positive integers satisfying $h_i\le\xi_i/y$. Suppose
the quantity $\eta_\hsdots$ satisfies
$\eta_\hsdots\asymp(y\max\{\hsdots\})^{1/d}(\hdots)^{-1}$. Then
\begin{equation}
\eta_\hsdots \le y^{1/d}(h_1^{1/d}+\dots+h_k^{1/d})(\hdots)^{-1}
\label{etaineq}
\end{equation}
and $\log\eta_\hsdots\gg_{k,d,U}\log x$ uniformly for $x\ge1$.
\label{etalem}
\end{lemma}

\begin{proof}
The inequality~(\ref{etaineq}) is easy to see by noting that
\begin{equation*}
(\max\{\hsdots\})^{1/d} = \max\{h_1^{1/d},\dots,h_1^{1/d}\} \le
h_1^{1/d}+\dots+h_1^{1/d}.
\end{equation*}
As for the second assertion, the hypothesis $h_i\le\xi_i/y$ implies that
\begin{equation*}
\begin{split}
\eta_\hsdots \gg {(y\max\{\hsdots\})^{1/d}\over\hdots} &\ge
y^{1/d}(\max\{\hsdots\})^{1/d-k} \\
&\ge y^k(\max\{\xi_1,\dots,\xi_k\})^{1/d-k} \gg x^{k/u}x^{1-dk}.
\end{split}
\end{equation*}
Since the exponent $k/u+1-dk$ of $x$ is positive by the upper
bound~(\ref{Urange}) on $u$, this establishes the lemma.
\end{proof}

We now have the tools we need to prove:

\restate{brunprop}

\begin{proof}
All of the constants implicit in the $\ll$-notation in this proof may
depend on the polynomial $f$. We use the bound given by
Lemma~\ref{brunlem} on each term in the inner sum, yielding
\begin{equation}
\sum_{b\in\Rhs} \pi(f\hhb;\eta_\hsdots) \ll
{\eta_\hsdots\over\log^k\eta_\hsdots} G(f;\hdots) \sum_{b\in\Rhs} 1.
\label{threekings}
\end{equation}
By the second assertion of Lemma~\ref{etalem}, we may replace
the denominator $\log^k\eta_\hsdots$ by $\log^kx$. In addition, since
the $h_i$ are pairwise coprime we may factor the term
$G(f;\hdots)=G(f;h_1)\dots G(f;h_k)$; and the final sum in
equation~(\ref{threekings}) is precisely
$\sigma(f_1;h_1)\dots\sigma(f_k;h_k)$ by the definition of
$\Rhs$. Therefore
\begin{multline*}
\mathop{\sum_{h_1\le\xi_1/y} \dots \sum_{h_k\le\xi_k/y}}_{(h_i,h_j)=1
\; (1\le i<j\le k)} \sum_{b\in\Rhs} \pi(f\hhb;\eta_\hsdots) \\
\ll {1\over\log^kx} \sum_{h_1\le\xi_1/y} \dots \sum_{h_k\le\xi_k/y}
\eta_\hsdots G(f;h_1)\sigma(f_1;h_1)\dots G(f;h_k)\sigma(f_k;h_k),
\end{multline*}
where we have deleted the condition $(h_i,h_j)=1$.

By the first assertion of Lemma~\ref{etalem}, this estimate becomes
\begin{multline}
\mathop{\sum_{h_1\le\xi_1/y} \dots \sum_{h_k\le\xi_k/y}}_{(h_i,h_j)=1
\; (1\le i<j\le k)} \sum_{b\in\Rhs} \pi(f\hhb;\eta_\hsdots) \\
\ll {y^{1/d}\over\log^kx} \sum_{i=1}^k \sum_{h_1\le\xi_1/y} \dots
\sum_{h_k\le\xi_k/y} {h_i^{1/d} G(f;h_1)\sigma(f_1;h_1)\dots
G(f;h_k)\sigma(f_k;h_k)\over\hdots}.
\label{i}
\end{multline}
For the term with $i=1$, for instance, we may use the second claim of
Lemma~\ref{oneDlem} with $\beta=1/d$ to bound the sum over $h_1$, and
the first claim of Lemma~\ref{oneDlem} to bound the remaining sums
over $h_2,\dots,h_k$:
\begin{multline*}
\sum_{h_1\le\xi_1/y} {G(f;h_1)\sigma(f_1;h_1)\over h_1^{1-1/d}}
\prod_{i=2}^k \sum_{h_k\le\xi_k/y} {G(f;h_i)\sigma(f_i;h_i)\over h_i}
\ll \big( {\xi_1\over y} \big)^{1/d} \prod_{i=2}^k \log{\xi_i\over y}
\ll {x\log^{k-1}x\over y^{1/d}}
\end{multline*}
since $\log(\xi_i/y)\ll\log x$. The other terms are similar, and we
see from the estimate~(\ref{i}) that
\begin{equation*}
\mathop{\sum_{h_1\le\xi_1/y} \dots \sum_{h_k\le\xi_k/y}}_{(h_i,h_j)=1
\; (1\le i<j\le k)} \sum_{b\in\Rhs} \pi(f\hhb;\eta_\hsdots) \ll
{x\over\log x},
\end{equation*}
which establishes the proposition.
\end{proof}

%\make{2mult}{Multiplicative Functions of Several Variables}
\section{Partial Summation}\setcounter{equation}{0}\label{psappendixsec}\noindent The proof of Proposition~\ref{lipsprop} will occupy all of this section.
Our first goal is to understand an unweighted version of the sum
on the left-hand side of equation~(\ref{werenot}), where the
$\li_\hsdots(f;x)$ terms have been omitted. Proposition~\ref{mv2prop} of
Appendix~A is precisely the tool we need to evaluate asymptotically a
multivariable sum of multiplicative functions of this form. The beauty of
the outcome is that the local factors conspire to make the constant in
the asymptotic formula exactly equal to $C(f)^{-1}$.

\begin{lemma}
Let $f$ be a polynomial with integer coefficients that is squarefree,
admissible, and exclusive, and let $f_1,\dots,f_k$ be the distinct
irreducible factors of $f$. For all real numbers $1\le
t_1,\dots,t_k\le x$, we have
\begin{equation*}
\mathop{\sum_{h_1\le t_1}\dots\sum_{h_k\le t_k}}_{(h_i,h_j)=1 \; (1\le
i<j\le k)} {G(f;h_1\dots h_k) \sigma^*(f_1;h_1)\dots\sigma^*(f_k;h_k)
\over h_1\dots h_k} = C(f)^{-1} \bigg( \prod_{i=1}^k \log t_i \bigg) +
O_f(\log^{k-1}x).
\end{equation*}
\label{cancellem}
\end{lemma}

\begin{proof}
All of the constants implicit in the $O$-notation in this proof may
depend on the polynomial $f$ (and thus on $k$ and its irreducible
factors $f_i$ as well). For each $1\le i\le k$, the function
$G(f;n)\sigma^*(f_i;n)$ is a nonnegative multiplicative function
satisfying $G(f;n)\sigma^*(f_i;n)\le G(f;n)\sigma(f_i;n)\ll n^\ep$
by Lemma~\ref{eplem} (note that each factor $f_i$ of the admissible
polynomial $f$ is itself admissible). Moreover, by
Lemma~\ref{kappaverifylem} the asymptotic formula~(\ref{kaphyp}) holds
for each $G(f;n)\sigma^*(f_i;n)$ with $\kappa_i=1$. We may therefore
conclude from Proposition~\ref{mv2prop} that
\begin{multline*}
\mathop{\sum_{h_1\le t_1}\dots\sum_{h_k\le t_k}}_{(h_i,h_j)=1 \; (1\le
i<j\le k)} {G(f;h_1\dots h_k) \sigma^*(f_1;h_1)\dots\sigma^*(f_k;h_k)
\over h_1\dots h_k} \\
= c(G(f;n)\sigma^*(f_1;n),\dots,G(f;n)\sigma^*(f_k;n)) \prod_{i=1}^k
(\log t_i + O(1)).
\end{multline*}
We have $\prod_{i=1}^k (\log t_i + O(1)) = \prod_{i=1}^k \log t_i +
O(\log^{k-1}x)$ since each $t_i\le x$, and so to establish the lemma it
suffices to show that
\begin{equation}
c\big( G(f;n)\sigma^*(f_1;n),\dots,G(f;n)\sigma^*(f_k;n) \big) = C(f)^{-1}.
\label{cfsuffit}
\end{equation}

Since each $\kappa_i=1$, we see from the definition~(\ref{cggdef})
that
\begin{multline}
c\big( G(f;n)\sigma^*(f_1;n),\dots,G(f;n)\sigma^*(f_k;n) \big) \\
= \Gamma(2)^{-k} \prod_p \big( 1-\frac1p \big)^k \bigg( 1 +
\sum_{\nu=1}^\infty {G(f;p^\nu)\sigma^*(f_1;p^\nu) + \dots +
G(f;p^\nu)\sigma^*(f_k;p^\nu)\over p^\nu} \bigg).
\label{reallycf}
\end{multline}
Rewriting
\begin{equation*}
%1 + 
\sum_{\nu=1}^\infty {G(f;p^\nu)\sigma^*(f_1;p^\nu) + \dots +
G(f;p^\nu)\sigma^*(f_k;p^\nu)\over p^\nu} % \\
= %1 +
\big( 1-\frac{\sigma(f;p)}p \big)^{-1} \sum_{i=1}^k
\bigg( \sum_{\nu+1}^\infty {\sigma^*(f_i;p^\nu)\over p^\nu} \bigg)
%\big( {\sigma^*(f_i;p)\over p} + {\sigma^*(f_i;p^2)\over p^2} + \dotsb
%\big)
\end{equation*}
by the definition of $G$, we note that the inner sum on the right-hand
side is a telescoping series by the definition of $\sigma^*$:
\begin{equation*}
\sum_{\nu+1}^\infty {\sigma^*(f_i;p^\nu)\over p^\nu} =
\sum_{\nu+1}^\infty \big( {\sigma(f_i;p^\nu)\over p^\nu} -
{\sigma(f_i;p^{\nu+1})\over p^{\nu+1}} \big) = {\sigma(f_i;p)\over p}.
%{\sigma^*(f_i;p)\over p} + {\sigma^*(f_i;p^2)\over p^2} + \dotsb \\
%= \big( {\sigma(f_i;p)\over p} - {\sigma(f_i;p^2)\over p^2} \big) +
%\big( {\sigma(f_i;p^2)\over p^2} - {\sigma(f_i;p^3)\over p^3} \big) +
%\dotsb = {\sigma(f_i;p)\over p}.
\end{equation*}
Thus equation~(\ref{reallycf}) becomes
\begin{equation*}
c\big( G(f;n)\sigma^*(f_1;n),\dots,G(f;n)\sigma^*(f_k;n) \big) = \prod_p
\big( 1-\frac1p \big)^k \bigg( 1 + \big( 1-\frac{\sigma(f;p)}p \big)^{-1}
\sum_{i=1}^k {\sigma(f_i;p)\over p} \bigg).
\end{equation*}
Since $f$ is an exclusive polynomial, for any prime $p$ the number of
roots of $f\mod p$ is equal to the sum of the numbers of roots of each
$f_i\mod p$; in other words, $\sum_{i=1}^k \sigma(f_i;p) =
\sigma(f;p)$. Therefore
\begin{equation*}
\begin{split}
c\big( G(f;n)\sigma^*(f_1;n),\dots,G(f;n)\sigma^*(f_k;n) \big) &= \prod_p
\big( 1-\frac1p \big)^k \big( 1 + \big( 1-\frac{\sigma(f;p)}p \big)^{-1}
\cdot{\sigma(f;p)\over p} \big) \\
&= \prod_p \big( 1-\frac1p \big)^k \big( 1-\frac{\sigma(f;p)}p
\big)^{-1} = C(f)^{-1},
\end{split}
\end{equation*}
which establishes equation~(\ref{cfsuffit}) and hence the lemma.
\end{proof}

Next we use partial summation to convert the asymptotic
formula in Lemma~\ref{cancellem} to the asymptotic formula asserted
in Proposition~\ref{lipsprop}, completing the proof of that
proposition. The partial summation argument is not deep but is quite
messy, both because the function $\li_\hsdots(f;x)$ defined in
equation~(\ref{lihsFxdef}) is somewhat complicated, and because the
$k$-fold sum in the statement of Proposition~\ref{lipsprop} requires
an inductive argument with partial summation being employed in each
variable. Addressing the former of these difficulties, the following
lemma gives the asymptotics of the function $\li_\hsdots(f;x)$.  We
recall from page~\pageref{defspage} the definitions of the parameters
$d$, $k$, $u$, $U$, $\xi_i$, and $y$ for use in this section, and we
allow all constant implicit in the $\ll$ and $O$-notations in this
section to depend on the polynomial $f$ (and thus on $d$, $k$, and the
$\xi_i$ as well) and on $U$.

\begin{lemma}
For any effective polynomial $f$ and any positive integers $h_1$,
\dots, $h_k$ satisfying $h_i\le\xi_i/y$, we have
\begin{equation}
\li_\hsdots(f;x) = \frac x{\pri k (\log\xi_i/h_i)} + O\big( \frac
x{(\log x)^{k+1}} \big)
\label{lisizeeq}
\end{equation}
uniformly for $x\ge1$.
\label{lisizelem}
\end{lemma}

\begin{proof}
By adjusting the constant implicit in the $O$-notation if necessary,
it suffices to establish the asymptotic formula~(\ref{lisizeeq}) when
$x$ is sufficiently large (in terms of $f$ and $U$). Notice that
$f_i(x/(\log x)^{k+1}) > x/(\log x)^{k+1}$ for each $1\le i\le k$,
from the fact that $f$ is effective. Notice also that each $\xi_i/y\ll
x^{d-1/U}$, where the exponent $d-1/U$ is strictly less than
$1/k\le 1$. Therefore for each $1\le i\le k$, the expression $f_i(x/(\log
x)^{k+1}) \big/ (\xi_i/y)$ tends to infinity with $x$. We assume that $x$
is so large that
\begin{equation}
\min_{1\le i\le k} \big\{ {f_i(x/(\log x)^{k+1}) \over \xi_i/y} \big\}
\ge2.
\label{ninestar}
\end{equation}
By the monotonicity of the $f_i$ and the hypothesis bounding the $h_i$,
this certainly implies that
\begin{equation}
\min_{1\le i\le k} \{f_i(t)/h_i\} \ge 2  \label{xrealybig}
\end{equation}
for any $t\ge x/(\log x)^{k+1}$.

We define $v(t) = 1/\pri k \log(f_i(t)/h_i)$, the integrand in the
definition~(\ref{lihsFxdef}) of $\li_\hsdots(f;x)$. Splitting that
integral at the point $x/(\log x)^{k+1}$ yields
\begin{equation*}
\li_\hsdots(f;x) = \int\limits\begin{Sb}0<t\le x/(\log x)^{k+1} \\
\min\{f_1(t)/h_1,\dots,f_k(t)/h_k\}\ge2\end{Sb} v(t)\,dt +
\int\limits\begin{Sb}x/(\log x)^{k+1}<t\le x \\
\min\{f_1(t)/h_1,\dots,f_k(t)/h_k\}\ge2\end{Sb} v(t)\,dt.
\end{equation*}
We estimate the first integral trivially by noting that $v(t)\le(\log
2)^{-k}\ll1$ when $t$ lies in the range of integration. In the second
integral, the condition $\min\{f_1(t)/h_1,\dots,f_k(t)/h_k\}\ge2$
holds since we are considering only values of $x$ that are so
large that the inequality~(\ref{xrealybig}) is satisfied. Therefore
\begin{equation*}
\li_\hsdots(f;x) = \int\limits_{x/(\log x)^{k+1}}^x v(t)\, dt + O\big(
{x\over(\log x)^{k+1}} \big).
\end{equation*}
and thus to establish the lemma it suffices to show that
\begin{equation}
\int\limits_{x/(\log x)^{k+1}}^x v(t)\,dt = xv(x) + O\big( \frac
x{(\log x)^{k+1}} \big).
\label{suffiteq}
\end{equation}

We accomplish this by integrating by parts:
\begin{equation}
\begin{split}
\int\limits_{x/(\log x)^{k+1}}^x v(t)\, dt &= tv(t)
\mathop{\big|}\nolimits_{x/(\log x)^{k+1}}^x - \int\limits_{x/(\log
x)^{k+1}}^x tv'(t)\,dt \\
&= xv(x) + O\big( {x v(x/(\log x)^{k+1})\over(\log x)^{k+1}} +
x\max\{tv'(t)\colon x/(\log x)^{k+1}\le t\le x\} \big).
\end{split}
\label{byparts}
\end{equation}
For any $t\ge x/(\log x)^{k+1}$, we note that
\begin{equation}
\log{f_i(t)\over h_i} \ge \log\big( {x/(\log x)^{k+1}\over\xi_i/y}
\big) \gg\log x
\label{fiuhi}
\end{equation}
for each $1\le i\le k$ (see the remarks preceding
equation~(\ref{ninestar})), and so we have $v(t)\ll(\log x)^{-k}$ for
$t$ in this range. Also, by logarithmic differentiation we see that
\begin{equation*}
v'(t) = v(t) \bigg( \sum_{i=1}^k {-f_i'(t)\over
f_i(t)\log(f_i(t)/h_i)} \bigg) \ll \frac1{\log^k x} \sum_{i=1}^k
{f_i'(t)\over f_i(t)\log x} \ll {1\over t(\log x)^{k+1}}
\end{equation*}
from the estimate~(\ref{fiuhi}) and the fact that
$f_i'(t)/f_i(t)\ll1/t$ for any polynomial $f_i$. Given these estimates
for $v(t)$ and $v'(t)$, we see that equation~(\ref{byparts}) implies
equation~(\ref{suffiteq}) and hence the lemma.
\end{proof}

The next lemma encapsulates the $k$-fold partial summation argument
mentioned prior to the statement of Lemma~\ref{lisizelem}.

\begin{lemma}
For each integer $0\le l\le k$, the asymptotic formula
\begin{multline*}
\mathop{\sum_{n_1\le t_1}\dots\sum_{n_k\le t_k}}_{(n_i,n_j)=1 \; (1\le
i<j\le k)} {G(f;n_1\dots n_k)\sigma^*(f_1;n_1)\dots
\sigma^*(f_k;n_k)\over n_1\dots n_k \cdot
(\log\xi_1/n_1)\dots(\log\xi_l/n_l)} \\
= C(f)^{-1} \bigg( \pri l \log \big( {\log\xi_i\over\log\xi_i/t_i}
\big) \bigg) \bigg( \prlk \log t_i \bigg) + O((\log x)^{k-l-1})
\end{multline*}
holds for all real numbers $t_1,\dots,t_k$ satisfying $1\le
t_i\le\xi_i/y$.
\label{lipslem}
\end{lemma}

\begin{proof}
Set
\begin{equation*}
W_l(t_1,\dots,t_k) = \mathop{\sum_{n_1\le t_1}\dots\sum_{n_k\le
t_k}}_{(n_i,n_j)=1 \; (1\le i<j\le k)} {G(f;n_1\dots
n_k)\sigma^*_1(n_1)\dots \sigma^*_k(n_k)\over n_1\dots n_k \cdot
(\log\xi_1/n_1)\dots(\log\xi_l/n_l)}.
\end{equation*}
The key to the proof is to notice that $W_l$ can be expressed in terms
of $W_{l-1}$ by partial summation:
\begin{equation}
\begin{split}
W_l(t_1,\dots,t_k) &= \int\limits_{t=1^-}^{t_l}
{dW_{l-1}(t_1,\dots,t_{l-1},t,t_{l+1},\dots,t_k)\over\log\xi_l/t} . \\
&= {W_{l-1}(t_1,\dots,t_k)\over\log\xi_l/t_l} - \int\limits_1^{t_l}
{W_{l-1}(t_1,\dots,t_{l-1},t,t_{l+1},\dots,t_k) \,dt\over
t(\log\xi_l/t)^2}
\end{split}
\label{instep}
\end{equation}
upon integration by parts.

We proceed by induction on $l$; the base case $l=0$ follows immediately
from Lemma~\ref{cancellem}, since $\xi_i/y\le x$ for $x$ sufficiently
large. For the inductive step, suppose that the lemma holds for the case
$l-1$, so that we know the asymptotic formula for $W_{l-1}$; we can then
insert this asymptotic formula into equation~(\ref{instep}) to obtain
\begin{equation}
\begin{split}
W_l(t_1, & \dots,t_k) = \bigg\{ C(f)^{-1} \bigg( \pri{l-1} \log \big(
{\log\xi_i\over\log\xi_i/t_i} \big) \bigg) \bigg( \prod_l^k \log t_i
\bigg) + O((\log x)^{k-l}) \bigg\} {1\over\log\xi_l/t_l} \\
&\qquad- \int\limits_1^{t_l} \bigg\{ C(f)^{-1} \bigg( \pri{l-1} \log
\big( {\log\xi_i\over\log\xi_i/t_i} \big) \bigg) \log t \bigg( \prlk
\log t_i \bigg) + O((\log x)^{k-l}) \bigg\} {dt\over t(\log\xi_l/t)^2}
\\
&= C(f)^{-1} \bigg( \pri{l-1} \log \big( {\log\xi_i\over\log\xi_i/t_i}
\big) \bigg) \bigg( \prlk \log t_i \bigg) \bigg\{ {\log
t_l\over\log\xi_l/t_l} - \int\limits_1^{t_l} {\log t \,dt\over
t(\log\xi_l/t)^2} \bigg\} \\
&\qquad+ O\bigg( {(\log x)^{k-l}\over\log\xi_l/t_l} +
\int\limits_1^{t_l} {(\log x)^{k-l} \,dt\over t(\log\xi_l/t)^2}
\bigg).
\end{split}
\label{asdesired}
\end{equation}
Since $\log\xi_l/t\ge\log y\gg\log x$ for all $1\le t\le t_l$, this
error term is $\ll(\log x)^{k-l-1}$. Moreover, making the change of
variables $t=\log t$ yields
\begin{equation*}
\begin{split}
\int\limits_1^{t_l} {\log t \,dt\over t(\log\xi_l/t)^2} =
\int\limits_0^{\log t_l} {t\,dt\over(\log\xi_l-t)^2} &= \big(
\log(\log\xi_l-t) + {\log\xi_l\over\log\xi_l-t} \big)
\mathop{\big|}\nolimits_0^{\log t_l} \\
&= \big( \log(\log\xi_l/t_l) + {\log\xi_l\over\log\xi_l/t_l} \big) -
\big( \log(\log\xi_l) + 1 \big) \\
&= {\log t_l\over\log\xi_l/t_l} - \log\big (
{\log\xi_l\over\log\xi_l/t_l} \big).
\end{split}
\end{equation*}
Therefore equation~(\ref{asdesired}) becomes
\begin{equation*}
W_l(t_1,\dots,t_k) = C(f)^{-1} \bigg( \pri{l-1} \log \big(
{\log\xi_i\over\log\xi_i/t_i} \big) \bigg) \bigg( \prlk \log t_i
\bigg) \log\big( {\log\xi_l\over\log\xi_l/t_l} \big) + O((\log
x)^{k-l-1}),
\end{equation*}
which is the desired asymptotic formula for the case $l$. This
establishes the lemma.
\end{proof}

We are now ready to establish:

\restate{lipsprop}

\begin{proof}
By Lemma~\ref{lisizelem},
\begin{equation*}
\begin{split}
\mathop{\sum_{h_1\le\xi_1/y}\dots\sum_{h_k\le\xi_k/y}}_{(h_i,h_j)=1 \;
(1\le i<j\le k)} & {G(f;\hdots)\sigma^*(f_1;h_1)\dots
\sigma^*(f_k;h_k) \li_\hsdots(f;x) \over \hdots} \\
&\hskip-2em= x \mathop{\sum_{h_1\le\xi_1/y} \dots
\sum_{h_k\le\xi_k/y}}_{(h_i,h_j)=1 \; (1\le i<j\le k)} {G(f;h_1\dots
h_k)\sigma^*(f_1;h_1)\dots \sigma^*(f_k;h_k) \over h_1\dots h_k \cdot
(\log\xi_1/h_1)\dots(\log\xi_k/h_k)} \\
&+ O\bigg( {x\over(\log x)^{k+1}}
\mathop{\sum_{h_1\le\xi_1/y}\dots\sum_{h_k\le\xi_k/y}}_{(h_i,h_j)=1 \;
(1\le i<j\le k)} {G(f;h_1\dots h_k)\sigma^*(f_1;h_1)\dots
\sigma^*(f_k;h_k) \over h_1\dots h_k} \bigg).
\end{split}
\end{equation*}
The main term can be evaluated by Lemma~\ref{lipslem} with $l=k$,
while the error term is
\begin{multline*}
\le {x\over(\log x)^{k+1}} \sum_{h_1\le\xi_1/y}\dots\sum_{h_k\le\xi_k/y}
{G(f;h_1)\sigma(f_1;h_1)\dots G(f;h_k)\sigma(f_k;h_k) \over h_1\dots
h_k} \\
\ll {x\over(\log x)^{k+1}} \log\frac{\xi_1}y \dots
\log\frac{\xi_k}y \ll {x\over\log x}
\end{multline*}
by $k$ applications of Lemma~\ref{oneDlem}. We obtain
\begin{multline}
\mathop{\sum_{h_1\le\xi_1/y}\dots\sum_{h_k\le\xi_k/y}}_{(h_i,h_j)=1 \;
(1\le i<j\le k)} {G(f;\hdots)\sigma^*(f_1;h_1)\dots
\sigma^*(f_k;h_k) \li_\hsdots(f;x) \over \hdots} \\
= x \big( C(f)^{-1} \log \big( {\log\xi_1\over\log y} \big) \dots \log
\big( {\log\xi_k\over\log y} \big) + O((\log x)^{-1}) \big) + O\big(
{x\over\log x} \big).
\label{onthefly}
\end{multline}
Since $\log\xi_i=\log f_i(x)=d\log x+O(1)$ for each $1\le i\le k$, we
see that $\log\xi_i/\log y=du+O(1/\log x)$, whence the right-hand side
of equation~(\ref{onthefly}) becomes $C(f)^{-1} x\log^k(du) + O(x/\log
x)$. This establishes the proposition.
\end{proof}

\section{Smooth Shifted Primes}\setcounter{equation}{0}\label{linearsec}\noindent The purpose of this section is to establish Theorem~\ref{primethm}. The
proof of this theorem has a structure similar to that of
Theorem~\ref{mainthm}, though of course it is much simpler to describe
since we are dealing with a very concrete situation. First the number of
smooth shifted primes $q+a$ is related combinatorially to the numbers of
prime values of certain polynomials $f_{h,a}$ defined in
Lemma~\ref{Cfprimelem} below; then the conjectured asymptotic formulas
for the corresponding expressions $\pi(f_{h,a},\cdot)$ are used and the
resulting sums analyzed as in the proof of Theorem~\ref{mainthm}. Where
appropriate, therefore, we omit some of the details of the following
proofs and refer the reader to the relevant arguments in earlier sections.

To begin with, we need to understand the behavior of the Dickman
function $\rho(u)$ in the expanded range $u\le3$, which is the
subject of the first lemma.

\begin{lemma}
We have
\begin{equation}
1-\log u+\sum_{t^{1/u}<p\le t^{1/2}} p^{-1} \log\big( {\log
t/p\over\log p} \big) = \rho(u) + O\big( {1\over\log t} \big)
\label{abcde}
\end{equation}
uniformly for $1\le u\le3$ and $t>1$ (where the sum in
equation~(\ref{abcde}) is empty if $u\le2$).
\label{rhocheatlem}
\end{lemma}

\begin{proof}
We could use partial summation to asymptotically evaluate the sum in
equation~(\ref{abcde}) in terms of elementary functions and the
dilogarithm function, and then verify that the resulting expression
satisfies the differential-difference equation characterizing the
function~$\rho$. However, if we take it as known that
\begin{equation}
\Psi(t,t^{1/u}) = t\rho(u) + O\big( {t\over\log t} \big)  \label{known}
\end{equation}
uniformly for $1\le u\le3$ and $t>1$, then we can use the following
simpler argument. For $u$ in this range we can write
\begin{equation*}
\begin{split}
\Psi(t,t^{1/u}) &= \3{n\le t} - \sum_{t^{1/u}<p\le t} \3{n\le t\colon
p\mid n} + \sum \begin{Sb}t^{1/u}<p_1<p_2 \\ p_1p_2\le t\end{Sb}
\3{n\le t\colon p_1p_2\mid t} \\
&= t+O(1) - \sum_{t^{1/u}<p\le t} \big( {t\over p} +O(1) \big) +
\sum_{t^{1/u}<p_1\le t^{1/2}} \sum_{p_1<p_2\le t/p_1} \big( {t\over
p_1p_2} +O(1) \big)
\end{split}
\end{equation*}
(where the final sum on each line is empty if $u\le2$). Using Mertens'
formula this becomes
\begin{equation*}
\begin{split}
\Psi(t,t^{1/u}) &= t - t\big( \log\big( {\log t\over\log t^{1/u}} \big)
+ O\big( {1\over\log t^{1/u}} \big) \big) + O(\pi(t)) \\
&\qquad+ \sum_{t^{1/u}<p_1\le t^{1/2}} \bigg( \frac t{p_1} \big(
\log\big( {\log t/p_1\over\log p_1} \big) + O\big( {1\over\log p_1}
\big) \big) + O(\pi(t/p_1)) \bigg) \\
&= t\bigg( 1-\log u + \sum_{t^{1/u}<p\le t^{1/2}} p^{-1} \log\big(
{\log t/p\over\log p} \big) \bigg) + O\big( {t\over\log t} \big).
\end{split}
\end{equation*}
Comparing this to the known asymptotic formula~(\ref{known})
establishes the lemma.
\end{proof}

Next we must establish a result analogous to
Proposition~\ref{lipsprop}, giving an asymptotic formula for a
weighted sum of a certain multiplicative function. In the current
situation, the sum in question is a simpler one-dimensional sum; also,
the relevant constants $C(f_{h,a})$ can be evaluated explicitly,
obviating the need for a result analogous to
Proposition~\ref{cflbprop}. The following lemma provides an asymptotic
formula for the unweighted version of the appropriate sum.

\begin{lemma}
Let $a$ be a nonzero integer and let $f(t)=t(t-a)$, so that
$f_{h,a}(t)=(ht+a)t$ for any integer $h$. Then
\begin{equation*}
\sum \begin{Sb}h\le t\end{Sb} {C(f_{h,a})\over h} =
\log t + O_a(1)
\end{equation*}
uniformly for $t\ge1$.
\label{Cfprimelem}
\end{lemma}

\begin{proof}
Define the multiplicative function
\begin{equation*}
g(n)=\prod \begin{Sb}p\mid n \\ p>2\end{Sb} {p-1\over p-2}.
\end{equation*}
For the particular polynomial $f_{h,a}$, it is an easy exercise to
compute from the definition~(\ref{CFdef}) of $C(f_{h,a})$ that
\begin{equation*}
C(f_{h,a}) = \begin{cases} 0, &\rmif 2\dnd ha\hbox{ or }(h,a)>1, \\
2C_2g(ha), &\rmif 2\mid ha\hbox{ and }(h,a)=1, \end{cases}
\end{equation*}
where $C_2$ is the twin primes constant
\begin{equation*}
C_2 = \prod_{p>2} \big( 1-{1\over(p-1)^2} \big).
\end{equation*}
If we assume that $a$ is even, this gives
\begin{equation}
\sum \begin{Sb}h\le t\end{Sb} {C(f_{h,a})\over h} =
2C_2g(a) \sum \begin{Sb}h\le t \\ (h,a)=1\end{Sb} {g(h)\over h}.
\label{Ctog}
\end{equation}

Note that $g(n)\le 2^{\omega(n)} \ll_\ep n^\ep$ for any positive $\ep$
(see the proof of Lemma~\ref{eplem}), and that
\begin{equation}
\sum_{p\le w} {g(p)\log p\over p} = {\log2\over2} + \sum_{2<p\le w} {\log p\over p}
\big( 1 + {1\over p-2} \big) = \log w + O(1),
\label{govergeq}
\end{equation}
since the sum $\sum_{p>2}(\log p)/p(p-2)$ converges. We can thus
apply Proposition~\ref{mvabprop}(b) with $\kappa=1$ to see that
\begin{equation*}
\sum \begin{Sb}h\le t \\ (h,a)=1\end{Sb} {g(h)\over h} =
c_a(g) \log t + O(\delta(a)),
\end{equation*}
where $c_a(g)$ is as defined in equation~(\ref{cqgdef}). In light of
equation~(\ref{Ctog}), therefore, to establish the lemma (for even
$a$) it suffices to show that $2C_2g(a)c_a(g)=1$. But by the
definition~(\ref{cqgdef}),
\begin{equation*}
\begin{split}
c_a(g) &= \prod_{p\mid a} \big( 1-\frac1p \big) \prod_{p\dnd a} \big(
1-\frac1p \big) \big( 1 + {g_a(p)\over p} + {g_a(p^2)\over p^2} +
\dots \big) \\
&= \prod_{p\mid a} \big( 1-\frac1p \big) \prod_{p\dnd a} {(p-1)^2\over
p(p-2)} \\
&= \frac12 \prod \begin{Sb}p\mid a \\ p>2\end{Sb} \big( {p-1\over p}
\big) \cdot C_2^{-1} \prod \begin{Sb}p\mid a \\ p>2\end{Sb}
{p(p-2)\over(p-1)^2} = {1\over2C_2g(a)},
\end{split}
\end{equation*}
as desired. This establishes the lemma when $a$ is even, and the same
argument slightly modified holds when $a$ is odd.
\end{proof}

Next we analyze the weighted version of the sum in
Lemma~\ref{Cfprimelem}, analogous to (but again simpler than) the
proof of Lemma~\ref{lipslem} in Section~\ref{psappendixsec}.

\begin{lemma}
Let $a$ be a nonzero integer and let $f(t)=t(t-a)$. Given real numbers
$x\ge1$ and $2\le u\le3$, let $y=x^{1/u}$ and $\xi=x-a$. Then
\begin{equation*}
\sum \begin{Sb}h\le\xi/y\end{Sb} C(f_{h,a}) \li\big (
f_{h,a};\frac\xi h \big) = {x\log u\over\log x} + O\big
( {x\over\log^2x} \big)
\end{equation*}
and
\begin{equation*}
\sum_{y<p\le\xi^{1/2}} \sum \begin{Sb}p<h\le\xi/p^2\end{Sb}
C(f_{ph,a}) \li\big ( f_{ph,a};\frac\xi{ph} \big) = {x\over\log x}
\sum_{y<p\le\xi^{1/2}} p^{-1} \log\big( {\log\xi/p\over\log p} \big) +
O\big( {x\over\log^2x} \big)
\end{equation*}
uniformly for $x\ge1+\max\{a,0\}$.
\label{twosumslem}
\end{lemma}

\begin{proof}
Using the same techniques as in the proofs of Lemmas~\ref{lichangelem}
and~\ref{lisizelem}, we can see that
\begin{equation*}
\li\big( f_{h,a};\frac\xi h \big) = {x\over h\log x\log(\xi/h)} +
O\big ( {x\over h\log^3x} \big),
\end{equation*}
and so
\begin{equation*}
\sum \begin{Sb}h\le\xi/y\end{Sb} C(f_{h,a}) \li\big (
f_{h,a};\frac\xi h \big) = {x\over\log x} \sum \begin{Sb}h\le\xi/y \\
(h,a)=1\end{Sb} {C(f_{h,a})\over h\log(\xi/h)} + O\bigg (
{x\over\log^3x} \sum \begin{Sb}h\le\xi/y\end{Sb}
{C(f_{h,a})\over h} \bigg).
\end{equation*}
The sum in the error term is $\ll\log\xi/y$ by Lemma~\ref{Cfprimelem},
whence the error term is $\ll x/\log^2x$ in its entirety. Using
partial summation (see the proof of Lemma~\ref{lipslem}), we can show
that Lemma~\ref{Cfprimelem} implies
\begin{equation*}
\sum \begin{Sb}h\le\xi/y\end{Sb} {C(f_{h,a})\over
h\log(\xi/h)} = \log\big( {\log\xi\over\log y} \big) + O\big(
{1\over\log x} \big) = \log u + O\big( {1\over\log x} \big).
\end{equation*}
This establishes the first claim of the lemma.

Similarly,
\begin{multline*}
\sum_{y<p\le\xi^{1/2}} \sum \begin{Sb}p<h\le\xi/p^2\end{Sb}
C(f_{ph,a}) \li\big ( f_{ph,a};\frac\xi{ph} \big) = {x\over\log x}
\sum_{y<p\le\xi^{1/2}} \frac1p \sum_{p<h\le\xi/p^2} {C(f_{ph,a})\over
h\log(\xi/ph)} \\
+ O\bigg( {x\over\log^3x} \sum_{y<p\le\xi^{1/2}} \frac1p \sum
\begin{Sb}p<h\le\xi/p^2\end{Sb} {C(f_{ph,a})\over h}
\bigg).
\end{multline*}
Again the error term can be shown to be $\ll x/\log^2x$, while the
inner sum in the main term can be evaluated by a similar partial
summation argument:
\begin{equation*}
\sum \begin{Sb}h\le\xi/p^2\end{Sb} {C(f_{ph,a})\over
h\log(\xi/ph)} = \log\big( {\log x/p\over\log p} \big) + O\big(
{1\over\log x} \big).
\end{equation*}
This establishes the second part of the lemma.
\end{proof}

We are now prepared to establish Theorem~\ref{primethm}.

\begin{pflike}{Proof of Theorem~\ref{primethm}:}
Let $a$ be a nonzero integer and define $f(t)=t(t-a)$, so that
$f_{h,a}(t)=(ht+a)t$. All constants implicit in the $\ll$ and
$O$-notations in this proof may depend on the nonzero integer $a$ and
thus on the polynomial $f$ as well. Let $x\ge1+\max\{a,0\}$ be a real
number, let $u$ and $U$ be real numbers satisfying $1\le u\le U<3$
(since the theorem is trivially true for $u<1$), and define $\xi=x-a$
and $y=x^{1/u}$. Reserving the letters $p$ and $q$ to denote primes
always, we have
\begin{equation*}
\begin{split}
\Phi_a(x,y) &= \3{q\le x\colon q-a\hbox{ is $y$-smooth}} \\
&= \pi(x) - \sum_{y<p\le\xi} \3{q\le x\colon p\mid(q-a)} + \sum
\begin{Sb}y<p_1<p_2 \\ p_1p_2\le\xi\end{Sb} \3{q\le x\colon
p_1p_2\mid(q-a)},
\end{split}
\end{equation*}
where this last sum (and similar ones to follow) is empty if
$y^2>\xi$, which holds for $x$ sufficiently large if and only if
$u<2$. If we make the substitutions $q=ph+a$ in the first sum and
$q=p_1p_2h+a$ in the second, and transform the resulting expressions
in a manner similar to the proof of Proposition~\ref{alterMprop}, the
end result is
\begin{equation}
\begin{split}
\Phi_a(x,y) &= \pi(x) - \sum_{h\le\xi/y} \3{y<p\le\xi/h\colon
p,\, hp+a\hbox{ are both prime}} \\
&\qquad + \sum_{y<p_1\le\xi^{1/2}} \sum_{h\le\xi/p_1^2}
\3{p_1<p_2\le\xi/p_1h\colon p_2,\, p_1p_2h+a\hbox{ are both
prime}}. \\
&= \pi(x) - \sum_{h\le\xi/y} \big( \pi\big( f_{h,a}; \frac\xi h \big)
- \pi(f_{h,a};y) \big) + \sum_{y<p\le\xi^{1/2}} \sum_{h\le\xi/p^2}
\big( \pi\big( f_{ph,a}; \frac\xi{ph} \big) - \pi(f_{ph,a};p) \big).
\end{split}
\label{likeprop3}
\end{equation}

Now by Lemma~\ref{brunlem},
\begin{equation*}
\sum_{h\le\xi/y} \pi(f_{h,a};y) \ll {y\over\log^2y} \sum_{h\le\xi/y}
G(f;h).
\end{equation*}
Since $G(f;p)=1+O(1/p)$ by Lemma~\ref{Olem}, we can apply the
estimate~(\ref{mvcor}) following from Proposition~\ref{mvprop}, with
$\beta=1$, to see that
\begin{equation*}
\sum_{h\le\xi/y} \pi(f_{h,a};y) \ll {y\over\log^2y} \cdot \frac\xi y
\ll {x\over\log^2x}.
\end{equation*}
Similarly,
\begin{equation*}
\begin{split}
\sum_{y<p\le\xi^{1/2}} \sum_{h\le\xi/p^2} \pi(f_{ph,a};p) &\ll
\sum_{y<p\le\xi^{1/2}} {p\over\log^2p} \sum_{h\le\xi/p^2} G(f;ph) \\
&\ll \log\xi^{1/2} \sum_{y<p\le\xi^{1/2}} {pG(f;p)\over\log^3p}
\sum_{h\le\xi/p^2} G(f;h) \\
&\ll \xi\log x \sum_{y<p\le\xi^{1/2}} {1\over p\log^3p}
\end{split}
\end{equation*}
since $G(f;p)\ll1$. The resulting expression is $\ll\xi\log
x/\log^3y\ll x/\log^2x$. Therefore from equation~(\ref{likeprop3}),
\begin{equation*}
\Phi_a(x,y) = {x\over\log x} - \sum_{h\le\xi/y} \pi\big( f_{h,a};
\frac\xi h \big) + \sum_{y<p\le\xi^{1/2}} \sum_{h\le\xi/p^2} \pi\big(
f_{ph,a}; \frac\xi{ph} \big) + O\big( {x\over\log^2x} \big)
\end{equation*}
since $\pi(x)=x/\log x+O(x/\log^2x)$ by the prime number theorem.

We now write this as
\begin{multline*}
\Phi_a(x,y) = {x\over\log x} - \sum \begin{Sb}h\le\xi/y\end{Sb}
C(f_{h,a}) \li\big( f_{h,a}; \frac\xi h \big) \\
+ \sum_{y<p\le\xi^{1/2}} \sum \begin{Sb}h\le\xi/p^2\end{Sb}
C(f_{ph,a}) \li\big( f_{ph,a}; \frac\xi{ph} \big) + H(x,y) + O\big(
{x\over\log^2x} \big),
\end{multline*}
where
\begin{equation}
H(x,y) = - \sum_{h\le\xi/y} E\big( f_{h,a}; \frac\xi h \big) +
\sum_{y<p\le\xi^{1/2}} \sum_{h\le\xi/p^2} E\big( f_{ph,a};
\frac\xi{ph} \big).
\end{equation}
We see from Lemma~\ref{twosumslem} that this is the same as
\begin{equation}
\begin{split}
\Phi_a(x,y) &= {x\over\log x} \bigg( 1 - \log u +
\sum_{y<p\le\xi^{1/2}} p^{-1} \log\big( {\log x/p\over\log p} \big)
\bigg) + H(x,y) + O\big( {x\over\log^2x} \big) \\
&= \pi(x)\rho(u) + H(x,y) + O\big( {\pi(x)\over\log x} \big)
\end{split}
\end{equation}
by Lemma~\ref{abcde} and the prime number theorem. Moreover, assuming
\hyp\ we can show by the method of the proof of
Proposition~\ref{usehypprop} that $H(x,y)\ll\pi(x)/\log x$. This
establishes Theorem~\ref{primethm}.
\qed
\end{pflike}

\setcounter{section}{0}
\renewcommand{\thesection}{Appendix \Alph{section}}
\section{Sums of Multiplicative Functions}\setcounter{equation}{0}\label{additionsec}\noindent\renewcommand{\thesection}{\Alph{section}}%
In this appendix we establish asymptotic formul\ae\ for summatory functions
associated with multiplicative functions $g(n)$, typified by the complete
one-dimensional sum
\begin{equation*}
M_g(x) = \sum_{n\le x} \gov n.  \label{Mgdef}
\end{equation*}
We shall also consider the modified sum $M_g(x,q)$, where the sum is taken
over only those integers coprime to $q$, as well as multidimensional
analogues, where the several variables of summation are restricted to be
coprime to one another. We are interested in such asymptotic formul\ae\ when
each multiplicative function $g$ is constant on average over primes, as is
usually the case for multiplicative functions that arise in sieve
problems, for example. Specifically, we impose the condition on $g$ that there
is a constant $\kappa=\kappa(g)$ such that
\begin{equation}
\sum_{p\le w} \glov p = \kappa\log w + O_g(1)  \label{consteq}
\end{equation}
for all $w\ge2$.

Although the ideas used in establishing the following proposition have
been part of the ``folklore'' for some time, the literature does not
seem to contain a result in precisely this form. Wirsing's pioneering
work \cite{Wir:DAVvSuMF}, for instance, requires $g$ to be a nonnegative
function and implies an asymptotic formula for $M_g(x)$ without a
quantitative error term; while Halberstam and Richert \cite[Lemma
5.4]{HalRic:SM} give an analogous result with a quantitative error
term, but one that requires $g$ to be supported on squarefree integers
in addition to being nonnegative. Both results are slightly too
restrictive for our purposes as stated.

Consequently we provide a self-contained proof of an asymptotic
formula for $M_g(x)$ with a quantitative error term, for
multiplicative functions $g$ that are not necessarily supported on
squarefree integers. The proof below, which is based on unpublished
work of Iwaniec (used with his kind permission) that stems from ideas
of Wirsing and Chebyshev, has the advantage that $g$ is freed from the
requirement of being nonnegative. We state the result in a more
general form than is required for our present purposes, with a mind
towards other applications and because the proof is exactly the same
in the more general setting.

\begin{proposition}
Suppose that $g(n)$ is a complex-valued multiplicative function such
that the asymptotic formula {\rm(\ref{consteq})} holds for some
complex number $\kappa=\xi+i\eta$ satisfying $\eta^2<2\xi+1$ (so that
$\xi>-1/2$ in particular). Suppose also that
\begin{equation}
\sum_p \aglov p \sum_{r=1}^\infty \agov{p^r} + \sum_p
\sum_{r=2}^\infty \aglov{p^r} < \infty,
\label{minorconverge}
\end{equation}
and that there exists a nonnegative real number $\beta=\beta(g)<\xi+1$
such that
\begin{equation}
\prod_{p\le x} \big( 1 + \agov{p} \big) \ll_g \log^{\beta} x
\label{minorprod}
\end{equation}
for all $x\ge2$. Then the asymptotic formula
\begin{equation}
M_g(x) = c(g)\log^\kappa x + O_g((\log x)^{\beta-1})  \label{mvform}
\end{equation}
holds for all $x\ge2$, where $\log^\kappa x$ denotes the principal
branch of $t^\kappa$, and $c(g)$ is defined by the convergent product
\begin{equation}
c(g) = \Gamma(\kappa+1)^{-1} \prod_p \big( 1-\frac1p \big)^\kappa \big(
1 + \gov p + \gov{p^2} + \dotsb \big).
\label{cgdef}
\end{equation}
\label{mvprop}
\end{proposition}

\noindent The conditions (\ref{minorconverge}) and (\ref{minorprod}) are
usually very easily verified in practice. We remark that the condition
(\ref{minorprod}) cannot hold with any $\beta<\abs\kappa$ if $g$ satisfies the
asymptotic formula (\ref{consteq}). The necessity that $\beta$ be less than
$\xi+1$, so that the formula (\ref{mvform}) is truly an asymptotic formula,
requires us to consider only those $\kappa$ for which $\abs\kappa<\xi+1$;
this is the source of the condition $\eta^2<2\xi+1$ on $\kappa$. We also
remark that from equation (\ref{mvform}), it follows easily by partial
summation that
\begin{equation}
\sum_{n<x}g(n)\ll_g x\log^{\beta-1}x  \label{mvcor}
\end{equation}
under the hypotheses of the proposition.

\vskip12pt

\begin{proof}
All of the constants implicit in the $O$- and $\ll$ symbols in this
proof may depend on the multiplicative function $g$, and thus on
$\kappa$ and $\beta$ as well. We begin by examining an analogue of
$M_g(x)$ weighted by a logarithmic factor. We have
\begin{equation}
\begin{split}
\sum_{n\le x} \glov n &= \sum_{n\le x} \gov n \sum_{p^r\mid\mid n}
\log p^r \\
&= \sum_{r=1}^\infty \sum_{p\le x^{1/r}} \glov{p^r} \sum
\begin{Sb}m\le x/p^r \\ p\dnd m\end{Sb} \gov m \\
&= \sum_{p\le x} \glov p \sum_{m\le x/p} \gov m - \sum_{p\le x} \glov
p \sum \begin{Sb}m\le x/p \\ p\mid m\end{Sb} \gov m \\
& \qquad {}+ \sum_{r=2}^\infty \sum_{p\le x^{1/r}} \glov{p^r} \sum
\begin{Sb}m\le x/p^r \\ p\dnd m\end{Sb} \gov m \\
&= \Sigma_1 - \Sigma_2 + \Sigma_3,
\end{split}
\label{sigmaseq}
\end{equation}
say. If we define the function $\Delta(x)$ by
\begin{equation}
\Delta(x) = \sum_{p\le x} \glov p - \kappa\log x,  \label{deltadef}
\end{equation}
then $\Sigma_1$ becomes
\begin{equation}
\Sigma_1 = \sum_{m\le x} \gov m \sum_{p\le x/m} \glov p = \kappa
\sum_{m\le x} \gov m \log\frac xm + \sum_{m\le x} \gov m \Delta\big(
\frac xm \big).
\label{using}
\end{equation}
Since $M_g(x)=1$ for $1\le x<2$ and
\begin{equation*}
M_g(x)\log x - \sum_{m\le x} \frac{g(m)\log m}m = \sum_{m\le x} \gov
m\log\frac xm = \int_1^x M_g(t) \frac{dt}t
\end{equation*}
by partial summation, we can rewrite equation (\ref{sigmaseq}) using
equation (\ref{using}) as
\begin{equation}
M_g(x)\log x - (\kappa+1) \int_2^x M_g(t)t^{-1} dt = E_g(x),
\label{usetwice}
\end{equation}
where we have defined
\begin{equation}
E_g(x) = (\kappa+1)\log 2 + \sum_{m\le x} \gov m \Delta\big( \frac xm \big)
- \Sigma_2 + \Sigma_3.  \label{Egdef}
\end{equation}

We integrate both sides of equation (\ref{usetwice}) against
$x^{-1}(\log x)^{-\kappa-2}$, obtaining
\begin{multline}
\int_2^x M_g(u)u^{-1}(\log u)^{-\kappa-1} du - (\kappa+1)\int_2^x u^{-1}(\log
u)^{-\kappa-2} \int_2^u M_g(t) t^{-1} dt\,du \\
= \int_2^x E_g(u) u^{-1}(\log u)^{-\kappa-2} du.
\label{okthen}
\end{multline}
Some cancellation can be obtained on the left-hand side by switching
the order of integration in the double integral and evaluating the new
inner integral; equation (\ref{okthen}) becomes simply
\begin{equation*}
(\log x)^{-\kappa-1} \int_2^x M_g(u)u^{-1} du = \int_2^x E_g(u) u^{-1}(\log
u)^{-\kappa-2} du.
\end{equation*}
We can substitute this into equation (\ref{usetwice}), divide by $\log
x$, and rearrange terms to get
\begin{equation}
M_g(x) = (\kappa+1) \log^\kappa x \int_2^x E_g(u) u^{-1}(\log
u)^{-\kappa-2} du + E_g(x)\log^{-1}x.
\label{rearrange}
\end{equation}

An upper bound for $E_g(x)$ is now needed. Since $\Delta(x)$ is
bounded from its definition (\ref{deltadef}) and the asymptotic
formula (\ref{consteq}), we have
\begin{equation}
\sum_{m\le x} \gov m \Delta\big( \frac xm \big) \ll \sum_{m\le x}
\agov m.
\label{modpf}
\end{equation}
We also have
\begin{equation}
\sum_{m\le x} \agov m \le \prod_{p\le x} \bigg( 1 + \sum_{r=1}^\infty
\agov{p^r} \bigg) \le \prod_{p\le x} \big( 1+\agov p \big) \prod_{p\le
x} \bigg( 1 + \sum_{r=2}^\infty \agov{p^r} \bigg).
\label{protoMabsg}
\end{equation}
Because the sum $\sum_p \sum_{r=2}^\infty \abs{g(p^r)}/p^r$ converges
by the hypothesis (\ref{minorconverge}), the last product in equation
(\ref{protoMabsg}) is bounded as $x$ tends to infinity. Therefore the
hypothesis (\ref{minorprod}) implies that
\begin{equation}
\sum_{m\le x} \agov m \ll \log^{\beta}x.
\label{Mabsg}
\end{equation}
The terms $\Sigma_2$ and $\Sigma_3$ can be estimated by
\begin{equation*}
\Sigma_2 = \sum_{p\le x} \glov p \sum_{r=1}^\infty \gov{p^r} \sum
\begin{Sb}l\le x/p^{r+1} \\ p\dnd l\end{Sb} \gov l \ll \sum_{p\le x}
\aglov p \sum_{r=1}^\infty \agov{p^r} \sum_{l\le x} \agov l
\end{equation*}
and
\begin{equation*}
\Sigma_3 \ll \sum_{p\le x} \sum_{r=2}^\infty \aglov{p^r} \sum_{m\le
x} \agov m,
\end{equation*}
and so both $\Sigma_2$ and $\Sigma_3$ are $\ll \log^{\beta}x$ by the
estimate (\ref{Mabsg}) and the hypothesis (\ref{minorconverge}).
Therefore, by the definition (\ref{Egdef}) of $E_g(x)$, we see that
\begin{equation}
E_g(x) \ll \log^{\beta}x.  \label{Egxbd}
\end{equation}

In particular, since $\beta<\xi+1$, we have
\begin{equation}
\int_x^\infty E_g(u) u^{-1}(\log u)^{-\kappa-2} du \ll \int_x^\infty
u^{-1}(\log u)^{\beta-\xi-2} du \ll (\log x)^{\beta-\xi-1},
\label{tenw}
\end{equation}
and so equation (\ref{rearrange}) and the bound (\ref{Egxbd}) give us
the asymptotic formula
\begin{equation}
M_g(x) = c(g) \log^\kappa x + O((\log x)^{\beta-1})  \label{allbutcg}
\end{equation}
for $x\ge2$, where
\begin{equation}
c(g) = (\kappa+1) \int_2^\infty E_g(u) u^{-1}(\log u)^{-\kappa-2} du.
\end{equation}

To complete the proof of the proposition, we need to show that $c(g)$
can be written in the form given by (\ref{cgdef}); we accomplish this
indirectly, using the asymptotic formula (\ref{allbutcg}). Consider
the zeta-function $\zeta_g(s)$ formed from $g$, defined by
\begin{equation*}
\zeta_g(s) = \sum_{n=1}^\infty {g(n)\over n^s}.
\end{equation*}
From the estimate (\ref{Mabsg}) and partial summation, we see that
$\zeta_g(s)$ converges absolutely for $s>1$ (we shall only need to
consider real values of $s$), and thus has an Euler product
representation
\begin{equation}
\zeta_g(s) = \prod_p \big( 1 + {g(p)\over p^s} + {g(p^2)\over p^{2s}}
+ \dotsb \big)
\label{eulerg}
\end{equation}
for $s>1$.

We can also use partial summation to write
\begin{equation}
\zeta_g(s+1) = s \int_1^\infty M_g(t) t^{-s-1}dt
\label{homestretch}
\end{equation}
for $s>0$. Since $M_g(x)=1$ for $1\le x<2$, it is certainly true that
\begin{equation*}
M_g(x) = c(g) \log^\kappa x + O(1+\log^\xi x)
\end{equation*}
in that range; using this together with the asymptotic formula
(\ref{allbutcg}), equation (\ref{homestretch}) becomes
\begin{multline*}
\zeta_g(s+1) = s \int_1^\infty c(g) \log^\kappa t \cdot t^{-s-1}
dt \\
{}+ O\bigg( s \int_1^2 (1+\log^\xi t) t^{-s-1}dt + s \int_2^\infty
(\log t)^{\beta-1} t^{-s-1}dt \bigg),
\end{multline*}
valid uniformly for $s>0$. Making the change of variables $t=e^{u/s}$
in all three integrals and multiplying through by $s^\kappa$ yields
\begin{equation}
\begin{split}
s^\kappa \zeta_g(s+1) &= c(g) \int_0^\infty u^\kappa e^{-u}du +
O\bigg( \int_0^{s\log2} (s^\xi+u^\xi)e^{-u}du + s^{\xi-\beta+1}
\int_{s\log2}^\infty u^{\beta-1}e^{-u}du \bigg) \\
&= c(g)\Gamma(\kappa+1) + O(s^{\xi-\beta+1}\log s^{-1})
\end{split}
\label{sto}
\end{equation}
as $s\to0^+$, where the exponent $\xi-\beta+1$ is positive and at most
1 (since $\beta\ge\abs\kappa\ge\xi$). Because the Riemann
$\zeta$-function satisfies $s\zeta(s+1)=1+O(s)$ as $s\to0^+$, equation
(\ref{sto}) implies
\begin{equation}
\zeta(s+1)^{-\kappa} \zeta_g(s+1) = c(g)\Gamma(\kappa+1) +
O(s^{\xi-\beta+1}\log s^{-1}).
\label{recent}
\end{equation}

On the other hand, from equation (\ref{eulerg}) we certainly have the
Euler product representation
\begin{equation*}
\zeta(s+1)^{-\kappa} \zeta_g(s+1) = \prod_p \big( 1-\frac1{p^{s+1}}
\big)^\kappa \big( 1 + {g(p)\over p^{s+1}} + {g(p^2)\over p^{2(s+1)}}
+ \dotsb \big)
\end{equation*}
for $s>0$, and one can show that in fact this Euler product converges
uniformly for $s\ge0$. The important contribution comes from the sum
$\sum_p (g(p)-\kappa)/p^{s+1}$, and we see from the hypothesis
(\ref{consteq}) and partial summation that
\begin{equation*}
\sum_{p>x} {g(p)-\kappa\over p^{s+1}} \ll \frac1{x^s\log x}
\end{equation*}
uniformly for $s\ge0$ and $x\ge2$. The remaining contributions can be
controlled using the hypothesis (\ref{minorconverge}).

Consequently, taking the limit of both sides of equation
(\ref{recent}) as $s\to0^+$ gives us
\begin{equation*}
\prod_p \big( 1-\frac1{p} \big)^\kappa \big( 1 + {g(p)\over p} +
{g(p^2)\over p^{2}} + \dotsb \big) = c(g) \Gamma(\kappa+1)
\end{equation*}
(where we have just shown that the product on the left-hand side
converges), which is equivalent to (\ref{cgdef}). This establishes the
proposition.
\end{proof}

From Proposition~\ref{mvprop} we can quickly derive a similar asymptotic
formula for the restricted sum
\begin{equation*}
M_g(x,q) = \sum \begin{Sb}n\le x \\ (n,q)=1\end{Sb} \gov n.
\end{equation*}

\begin{proposition}
Suppose that $g(n)$ satisfies the hypotheses of Proposition
\ref{mvprop}. Then the asymptotic formula
\begin{equation*}
M_g(x,q) = c_q(g)\log^\kappa x + O_g(\delta(q)(\log x)^{\beta-1})
\end{equation*}
holds uniformly for all $x\ge2$ and all nonzero integers $q$, where
\begin{equation}
c_q(g) = \Gamma(\kappa+1)^{-1} \big( \frac{\phi(q)}q \big)^\kappa
\prod_{p\dnd q} \big( 1-\frac1p \big)^\kappa \big( 1 + \gov p +
\gov{p^2} + \dotsb \big)
\label{firstcqgdef}
\end{equation}
and $\delta(q) = 1 + \sum_{p\mid q} \abs{g(p)}(\log p)/p$.
\label{coprimeprop}
\end{proposition}

\begin{proof}
We would like to apply Proposition \ref{mvprop} to the multiplicative
function $g_q(n)$ defined by
\begin{equation*}
g_q(n) = \begin{cases}g(n), &\rmif (n,q)=1, \\ 0, &\rmif
(n,q)>1.\end{cases}
\end{equation*}
Certainly $\abs{g_q(n)}\le\abs{g(n)}$, and so the estimates
(\ref{minorconverge}) and (\ref{minorprod}) for $g_q$ follow from the
same estimates for $g$. We also have
\begin{equation*}
\begin{split}
\sum_{p\le x} \frac{g_q(p)\log p}p = \sum \begin{Sb}p\le x \\ p\dnd
q\end{Sb} \glov p &= \sum_{p\le x} \glov p - \sum \begin{Sb}p\le x \\
p\mid q\end{Sb} \glov p \\
&= \kappa\log x + O_g(1) + O\bigg( \sum_{p\mid q} \aglov p \bigg)
\end{split}
\end{equation*}
from the assumption that $g$ satisfies equation (\ref{consteq}).
Therefore $g_q$ satisfies equation (\ref{consteq}) as well, with the
error term being $\ll_g \delta(q)$ uniformly in $x$.

If we keep this dependence on $q$ explicit throughout the proof of
Proposition \ref{mvprop}, the only modification necessary is to
include a factor of $\delta(q)$ on the right-hand sides of the
estimates (\ref{modpf}), (\ref{Egxbd}), and (\ref{tenw}) and in the
error term in equation (\ref{allbutcg}). Therefore, the application of
Proposition \ref{mvprop} to $g_q$ yields
\begin{equation*}
M_g(x,q) = M_{g_q}(x) = c(g_q) \log^\kappa x + O_g(\delta(q)(\log
x)^{\beta-1}),
\end{equation*}
where the implicit constant is independent of $q$. Because
\begin{equation*}
\begin{split}
c(g_q) &= \Gamma(\kappa+1)^{-1} \prod_p \big( 1-\frac1p \big)^\kappa \big(
1 + \gqov p + \gqov{p^2} + \dotsb \big) \\
&= \Gamma(\kappa+1)^{-1} \prod_{p\mid q} \big( 1-\frac1p \big)^\kappa
\prod_{p\dnd q} \big( 1-\frac1p \big)^\kappa \big( 1 + \gov p +
\gov{p^2} + \dotsb \big) = c_q(g),
\end{split}
\end{equation*}
the proposition is established.
\end{proof}

For convenience, we state below a particular case of Propositions~\ref{mvprop}
and~\ref{coprimeprop}, where the hypotheses have been somewhat simplified.

\begin{proposition}
Let $g(n)$ be a nonnegative multiplicative function satisfying $g(n)\ll
n^\alpha$ for some constant $\alpha<1/2$. Suppose that there is a real
number $\kappa$ (necessarily nonnegative) such that
\begin{equation}
\sum_{p\le w} \glov p = \kappa\log w + O_g(1)  \label{hypothesis}
\end{equation}
for all $w\ge2$. Then:
\begin{enumerate}
\item[\rm(a)] the asymptotic formula
\begin{equation}
\sum_{n\le x} \gov n = c(g)\log^\kappa x + O_g(\log^{\kappa-1}x)
\label{MVT}
\end{equation}
holds for all $x\ge2$, where $c(g)$ is defined by the convergent
product~(\ref{cgdef});
\item[\rm(b)] the asymptotic formula
\begin{equation*}
\sum \begin{Sb}n\le x \\ (n,q)=1\end{Sb} \gov n = c_q(g)\log^\kappa x
+ O_g(\delta(q)\log^{\kappa-1}x)
\end{equation*}
holds uniformly for all $x\ge2$ and all positive integers $q$, where
\begin{equation}
c_q(g) = c(g) \prod_{p\mid q} \big( 1 + \gov p + \gov{p^2} + \dotsb
\big)^{-1} \quad\hbox{and}\quad \delta(q) = 1 + \sum_{p\mid q} \glov
p.
\label{cqgdef}
\end{equation}
\end{enumerate}
\label{mvabprop}
\end{proposition}

\begin{proof}
We have
\begin{equation*}
\begin{split}
\sum_p \aglov p \sum_{r=1}^\infty \agov{p^r} + \sum_p \sum_{r=2}^\infty
&\aglov{p^r} \\
&\ll \sum_p {p^\alpha\log p\over p} \sum_{r=1}^\infty {p^{r\alpha}\over p^r} +
\sum_p \sum_{r=2}^\infty {p^{r\alpha}\log p^r\over p^r} \\
&= \sum_p {\log p\over p^{2-2\alpha}(1-p^{\alpha-1})} + \sum_p
{(2-p^{\alpha-1})\log p\over p^{2-2\alpha}(1-p^{\alpha-1})^2} \\
&\ll \sum_p {\log p\over p^{2-2\alpha}} < \infty,
\end{split}
\end{equation*}
since $2-2\alpha>1$. Hence the condition~(\ref{minorconverge}) is
satisfied. Also, since $g$ is a nonnegative function, we have
\begin{equation*}
\prod_{p\le x} \big( 1 + \agov{p} \big) \le \prod_{p\le x}
\exp\big( {|g(p)|\over p} \big) = \exp\bigg( \sum_{p\le x} {g(p)\over p}
\bigg) = \exp(\kappa\log x+O_g(1)) \ll_g \log^\kappa x
\end{equation*}
by the asymptotic formula~(\ref{consteq}), which shows that the
condition~(\ref{minorprod}) is also satisfied with $\beta=\kappa$. Therefore
the two parts of this proposition are just special cases of
Propositions~\ref{mvprop} and~\ref{coprimeprop}, respectively. The form given
in equation~(\ref{cqgdef}) for $c_q(g)$ is equivalent to the form given in
equation~(\ref{firstcqgdef}), since the assumption that $g$ is nonnegative
implies that the sum $(1+g(p)/p+g(p^2)/p^2+\dotsb)$ is nonzero.
\end{proof}

Finally, the following proposition analyzes the behavior of a general sum of
several multiplicative functions, where the variables of summation are not
permitted to have common factors.

\begin{proposition}
Let $k$ be a positive integer, and let $g_1(n)$, \dots, $g_k(n)$ be
nonnegative multiplicative functions. Suppose that there exists a real number
$\alpha<1/2$ such that, for each $1\le i\le k$, the estimate $g_i(n)\ll
n^\alpha$ holds for all positive integers $n$. Suppose further that
$\kappa_1$, \dots, $\kappa_k$ are real numbers such that, for each
$1\le i\le k$,
\begin{equation}
\sum_{p\le w} \gilov ip = \kappa_i\log w + O_{g_i}(1)  \label{kaphyp}
\end{equation}
for all $w\ge2$. Then the asymptotic
formula
\begin{equation}
\manysum {g_1(n_1)\dots g_k(n_k)\over n_1\dots n_k} = c(g_1,\dots,g_k)
\prod_{i=1}^k \big( \log^{\kappa_i}x_i + O_{g_1,\dots,g_k}\big( (\log
x_i)^{\kappa_i-1} \big) \big)
\label{forl}
\end{equation}
holds for all $x_1,\dots,x_k\ge1$, where $c(g_1,\dots,g_k)$ is defined
by the convergent product
\begin{multline}
c(g_1,\dots,g_k) = \big( \Gamma(\kappa_1+1)\dotsb\Gamma(\kappa_k+1)
\big)^{-1} \\
{}\times \prod_p \big( 1-\frac1p \big)^{\kappa_1+\dots+\kappa_k} \big(
1 + \frac{g_1(p)+\dots+g_k(p)}p + {g_1(p^2)+\dots+g_k(p^2)\over p^2} +
\dotsb \big).
\label{cggdef}
\end{multline}
\label{mv2prop}
\end{proposition}

\begin{proof}
All of the constants implicit in the $O$- and $\ll$ symbols in this
proof may depend on the multiplicative functions $g_i$, and thus on
$k$, $\alpha$, and the $\kappa_i$ as well. Define
\begin{equation*}
S_k = S_k(g_1,\dots,g_k;x_1,\dots,x_k) = \manysum {g_1(n_1)\dots
g_k(n_k)\over n_1\dots n_k}.
\end{equation*}
We establish the desired asymptotic formula~(\ref{forl}) for $S_k$ by
induction on $k$. The base case $k=1$ of the induction is exactly the
statement of Proposition~\ref{mvabprop}(a). Supposing now that we know that
the asymptotic formula~(\ref{forl}) holds for sums of the form $S_k$, we
wish to show that it holds for $S_{k+1}$.

We rewrite $S_{k+1}$ as
\begin{equation}
S_{k+1} = \manysum {g_1(n_1)\dots g_k(n_k)\over n_1\dots n_k} \sum
\begin{Sb}m\le x_{k+1} \\ (m,n_1\dots n_k)=1\end{Sb} \giov{k+1}m,
\label{firststep}
\end{equation}
and we can use Proposition~\ref{mvabprop}(b) to obtain an asymptotic formula
for this inner sum. For any index $i$, let $\gamma_i(n)$ be the multiplicative
function defined on prime powers by
\begin{equation*}
\gamma_i(p^\nu) = 1 + \giov ip + \giov i{p^2} + \dotsb,
\end{equation*}
which satisfies
\begin{equation}
\gamma_i(p) = 1 + O\big( \frac{p^\alpha}p + {p^{2\alpha}\over p^2} +
\dotsb \big) = 1 + O(p^{-1+\alpha})
\label{onepluso}
\end{equation}
by the hypothesized estimate on the size of $g_i$, and set
\begin{equation*}
\delta(n) = 1 + \sum_{p\mid n} {g_{k+1}(p)\log p\over p}.
\end{equation*}
Then by Proposition~\ref{mvabprop}(b), equation~(\ref{firststep}) becomes
\begin{equation}
\begin{split}
S_{k+1} &= \manysum \Big\{ {g_1(n_1)\dots g_k(n_k)\over n_1\dots n_k} \\
&\qquad\times \big( c(g_{k+1})\gamma_{k+1}(n_1\dots n_k)^{-1} (\log
x_{k+1})^{\kappa_{k+1}} + O\big( \delta(n_1\dots n_k) (\log
x_{k+1})^{\kappa_{k+1}-1} \big) \big) \Big\} \\
&= c(g_{k+1}) (\log x_{k+1})^{\kappa_{k+1}} T_k + O\big( (\log
x_{k+1})^{\kappa_{k+1}-1} U_k \big),
\end{split}
\label{halfway}
\end{equation}
where $c(g_{k+1})$ is defined in equation~(\ref{cgdef}), and where we
have defined
\begin{equation*}
T_k = \manysum {g_1(n_1)\gamma_{k+1}(n_1)^{-1}\dots
g_k(n_k)\gamma_{k+1}(n_k)^{-1}\over n_1\dots n_k}
\end{equation*}
(using the multiplicativity of $\gamma_{k+1}$ and the fact that the $n_i$
are pairwise coprime) and
\begin{equation*}
U_k = \manysum {g_1(n_1)\dots g_k(n_k) \delta(n_1\dots n_k) \over
n_1\dots n_k}.
\end{equation*}

First we obtain an asymptotic formula for $T_k$. The sum $T_k$ is
precisely of the form $S_k$ with each $g_i$ replaced by
$g_i\gamma_{k+1}^{-1}$. By equation~(\ref{onepluso}), for each $1\le i\le k$
the multiplicative function $g_i\gamma_{k+1}^{-1}$ satisfies
\begin{equation*}
\begin{split}
\sum_{p\le x} \frac{g_i(p)\gamma_{k+1}(p)^{-1}\log p}p &= \sum_{p\le x}
\gilov ip (1 + O(p^{-1+\alpha})) \\
&= \sum_{p\le x} \gilov ip + O\bigg( \sum_{p\le x} p^{-2+2\alpha}\log
p \bigg) = \kappa_i\log x + O(1)
\end{split}
\end{equation*}
since $\alpha<1/2$. We may therefore invoke the induction hypothesis
to show that the asymptotic formula
\begin{multline}
T_k = \manysum {g_1(n_1) \gamma_{k+1}(n_1)^{-1} \dots g_k(n_k)
\gamma_{k+1}(n_k)^{-1} \over n_1\dots n_k} \\
= c(g_1\gamma_{k+1}^{-1},\dots,g_k\gamma_{k+1}^{-1}) \prod_{i=1}^k \big(
\log^{\kappa_i}x_i + O\big( (\log x_i)^{\kappa_i-1}
\big) \big)
\label{invoke}
\end{multline}
holds for $T_k$.

Next we obtain an estimate for $U_k$. By the definition of $\delta(n)$
we have
\begin{equation}
\begin{split}
U_k &= \manysum {g_1(n_1)\dots g_k(n_k) \over n_1\dots n_k} \bigg( 1 +
\sum_{p\mid n_1\dots n_k} \gilov{k+1}p \bigg) \\
&= \manysum {g_1(n_1)\dots g_k(n_k) \over n_1\dots n_k} \\
&\qquad+ \sum_{p\le x_1\dots x_k} \gilov{k+1}p \mathop{\sum_{n_1\le
x_1}\dots\sum_{n_k\le x_k}} \begin{Sb}(n_i,n_j)=1 \; (1\le i<j\le k)
\\ p\mid n_1\dots n_k\end{Sb} {g_1(n_1)\dots g_k(n_k) \over n_1\dots
n_k} \\
&= S_k + \sum_{p\le x_1\dots x_k} \gilov{k+1}p \sum_{i=1}^k V_{i}(p),
\end{split}
\label{second}
\end{equation}
where we have defined
\begin{equation*}
V_{i}(p)=\mathop{\sum_{n_1\le x_1}\dots\sum_{n_k\le x_k}}
\begin{Sb}(n_i,n_j)=1 \; (1\le i<j\le k) \\ p\mid n_i\end{Sb}
{g_1(n_1)\dots g_k(n_k) \over n_1\dots n_k}.
\end{equation*}
We can rewrite $V_{k}(p)$, for example, as
\begin{equation*}
\begin{split}
V_{k}(p) &= \mathop{\sum_{n_1\le x_1}\dots\sum_{n_{k-1}\le
x_{k-1}}}_{(n_i,n_j)=1 \; (1\le i<j\le k-1)} {g_1(n_1)\dots
g_{k-1}(n_{k-1}) \over n_1\dots n_{k-1}} \sum \begin{Sb}n_k\le x_k \\
(n_k,n_1\dots n_{k-1})=1 \\ p\mid n_k\end{Sb} {g_k(n_k)\over n_k} \\
&= \mathop{\sum_{n_1\le x_1}\dots\sum_{n_{k-1}\le
x_{k-1}}}_{(n_i,n_j)=1 \; (1\le i<j\le k-1)} {g_1(n_1)\dots
g_{k-1}(n_{k-1}) \over n_1\dots n_{k-1}} \bigg( \sum_{\nu=1}^\infty
{g_k(p^\nu)\over p^\nu} \sum \begin{Sb}m\le x_k/p^\nu \\ (m,pn_1\dots
n_{k-1})=1\end{Sb} {g_k(m)\over m} \bigg).
\end{split}
\end{equation*}
In this innermost sum of nonnegative terms, we can delete the
restriction that $n_k$ be coprime to $p$ and extend the range of
summation from $m\le x_k/p^\nu$ to $m\le x_k$; this yields (after
renaming $m$ to $n_k$) the upper bound
\begin{equation*}
V_{k}(p) \le \manysum {g_1(n_1)\dots g_k(n_k) \over n_1\dots n_k}
\sum_{\nu=1}^\infty {g_k(p^\nu)\over p^\nu} = (\gamma_k(p)-1) S_k.
\end{equation*}
The same analysis shows that $(\gamma_i(p)-1)S_k$ is an
upper bound for $V_{i}(p)$ for each $1\le i\le k$.

With this upper bound, equation~(\ref{second}) becomes
\begin{equation*}
U_k \le S_k \bigg( 1+ \sum_{i=1}^k \sum_{p\le x_1\dots x_k}
{g_{k+1}(p)(\gamma_i(p)-1)\log p\over p} \bigg).
\end{equation*}
By the induction hypothesis for $S_k$, we see that
$S_k\ll(\log^{\kappa_1}x_1 \dots \log^{\kappa_k}x_k)$. Furthermore,
both $g_{k+1}(p)/p$ and each $\gamma_i(p)-1$ are $\ll p^{-1+\alpha}$ (the
latter by equation~(\ref{onepluso})). Therefore
\begin{equation*}
U_k \ll (\log^{\kappa_1}x_1 \dots \log^{\kappa_k}x_k) \bigg( 1 +
\sum_p p^{2\alpha-2}\log p \bigg) \ll \log^{\kappa_1}x_1 \dots
\log^{\kappa_k}x_k
\end{equation*}
since $\alpha<1/2$.

Using this estimate for $U_k$ and the asymptotic
formula~(\ref{invoke}) for $T_k$, we see that equation~(\ref{halfway})
becomes
\begin{equation}
\begin{split}
S_k &= c(g_{k+1}) \Big\{ (\log x_{k+1})^{\kappa_{k+1}}
c(g_1\gamma_{k+1}^{-1},\dots,g_k\gamma_{k+1}^{-1}) \prod_{i=1}^k \big(
\log^{\kappa_i}x_i + O\big( (\log x_i)^{\kappa_i-1} \big) \big) \Big\}
\\
&\qquad{}+ O\big( (\log x_{k+1})^{\kappa_{k+1}-1} \log^{\kappa_1}x_1
\dots \log^{\kappa_k}x_k \big) \\
&= c(g_{k+1}) c(g_1\gamma_{k+1}^{-1},\dots,g_k\gamma_{k+1}^{-1})
\prod_{i=1}^{k+1} \big( \log^{\kappa_i}x_i + O\big( (\log
x_i)^{\kappa_i-1} \big) \big).
\end{split}
\label{wouldbut}
\end{equation}
This would establish the lemma if only we had $c(g_1,\dots,g_{k+1})$
in place of the product $c(g_{k+1}) \* c(g_1\gamma_{k+1}^{-1},\dots,
g_k\gamma_{k+1}^{-1})$. However, the $\Gamma$-factors of these two
expressions are certainly equal by inspection. For each prime $p$, moreover,
the power of $(1-1/p)$ in the infinite products of the two expressions equals
$\kappa_1+\dots+\kappa_{k+1}$ in both cases, and we also have
\begin{equation*}
\begin{split}
\gamma_{k+1}(p) \times \bigg( 1 +{}& \frac{g_1(p)\gamma_{k+1}(p)^{-1} + \dots
+ g_k(p)\gamma_{k+1}(p)^{-1}}p \\
&\qquad{}+ {g_1(p^2)\gamma_{k+1}(p^2)^{-1} +
\dots+g_k(p^2)\gamma_{k+1}(p^2)^{-1}\over p^2} + \dotsb \bigg) \\
&= \gamma_{k+1}(p) + \frac{g_1(p)+\dots+g_k(p)}p +
{g_1(p^2)+\dots+g_k(p^2)\over p^2} + \dotsb \\
&= 1 + \frac{g_1(p)+\dots+g_{k+1}(p)}p +
{g_1(p^2)+\dots+g_{k+1}(p^2)\over p^2} + \dotsb.
\end{split}
\end{equation*}
Therefore the local factors in the infinite products of $c(g_{k+1})
\times c(g_1\gamma_{k+1}^{-1},\dots,g_k\gamma_{k+1}^{-1})$ and
$c(g_1,\dots,g_{k+1})$ are also equal, and so the asymptotic
formula~(\ref{wouldbut}) is equivalent to the statement of the lemma.
\end{proof}

\renewcommand{\thesection}{Appendix \Alph{section}}
\section{Prime Values of Linear Polynomials}\setcounter{equation}{0}\label{zultimatesec}\noindent\renewcommand{\thesection}{\Alph{section}}%
In this appendix, we show that in the case of linear
polynomials, \hyp\ is equivalent to a well-believed statement about
the number of primes in short segments of arithmetic progressions,
which for purposes of reference we shall call \hypap:

\begin{pflike}{\hypap.}
Given real numbers $0<\ep<1$ and $C>1$, the asymptotic formula
\begin{equation}
\pi(x;q,a) - \pi(x-y;q,a) = {y\over\phi(q)\log x} + O\big(
{y\over\phi(q)\log^2x} \big)
\label{hypapeq}
\end{equation}
holds uniformly for all real numbers $x$ and $y$ satisfying $1\le y\le
x\le y^C$ and all coprime integers $1\le a\le q\le y^{1-\ep}$.
\end{pflike}

Of course, it is equivalent to ask only that the asymptotic
formula~(\ref{hypapeq}) hold when $x$ and $y$ are sufficiently large,
by adjusting the constant implicit in the $O$-notation if
necessary. The conditions $x\le y^C$ and $q\le y^{1-\ep}$ mean that
the primes being counted are only polynomially large as a function of
the number of terms $y/q$ in the segment of the arithmetic
progression. This restriction is not made merely for simplicity:
Friedlander and Granville~\cite{FriGra:LttEoP}, expanding on the
ground-breaking ideas of Maier, showed that even in the case $y=x$, the
asymptotic formula~(\ref{hypapeq}) can fail when the size of $q$ is
$x/\log^D x$ for arbitrarily large~$D$. Certainly one can construct by
elementary methods, for any given $y$, an integer $x$ so that the
interval $[x-y,x]$ contains no primes whatsoever, so
that~(\ref{hypapeq}) cannot hold without some restriction on~$x$.

As remarked in Section~\ref{introsec}, \hyp\ holds automatically for
non-admissible polynomials. Note that a linear polynomial $qt+b$
(where by multiplying by $-1$ if necessary, we may assume that $q$ is
positive) is admissible if and only if $(b,q)=1$, in which case
$C(qt+b)$ is easily seen to equal $q/\phi(q)$. So for linear
polynomials, \hyp\ can be stated as follows:

\begin{pflike}{\hypo.}
Given a constant $B>0$, we have
\begin{equation}
\pi(qt+b;T) = {q\over\phi(q)} \li(qt+b;T) + O\big( {q\over\phi(q)}
{T\over\log^2 T} \big)
\label{hypoeq}
\end{equation}
uniformly for all real numbers $T\ge1$ and all coprime integers
$q$ and $b$ satisfying $1\le q\le T^B$ and $|b|\le T^B$.
\end{pflike}

Again it is clearly equivalent to ask the the asymptotic
formula~(\ref{hypoeq}) hold for sufficiently large integer values
of~$T$. Before demonstrating the equivalence between \hypo\ and
\hypap, we remark that $\li(F;T)$ can be expressed in terms of the
ordinary logarithmic integral $\li(x)$ when $F(t)=qt+b$ is a linear
polynomial. Assuming that $q$ and $b$ are positive, we can make the
change of variables $v=qt+b$ in the integral in the
definition~(\ref{liFxdef}) of $\li(F;T)$ to see that
\begin{equation}
\li(F;T) = \int\limits \begin{Sb}0<t<T \\ qt+b\ge2\end{Sb}
{dt\over\log(qt+b)} = \frac1q \int\limits_{\max\{b,2\}<v<qT+b}
{dv\over v} = {\li(qT+b)-\li(b)\over q} + O(1).
\label{asacons}
\end{equation}
In fact, this formula holds without the assumption that $q$ and $b$
are positive, if we make the conventions that $\li(x)=-\li(|x|)$ if
$x\le-2$ and $\li(x)=0$ if $|x|\le 2$.

We also need a lemma on the behavior of the logarithmic integral
$\li(x)$.

\begin{lemma}
We have
\begin{equation*}
\li(x)-\li(x-y) = {y\over\log x} + O\big( {y\over\log^2x} \big)
\end{equation*}
uniformly for $2\le y\le x-2$.
\label{calclem}
\end{lemma}

\begin{proof}
It is easily seen by integration by parts that
\begin{equation}
\li(x) = {x\over\log x} + O\big( {x\over\log^2x} \big).  \label{listand}
\end{equation}
First we consider the case where $x-x/\log x\le y\le x-2$. In this case,
we have
\begin{equation*}
\li(x-y) \le \li\big( {x\over\log x} \big) \ll {x\over\log^2x}
\end{equation*}
by equation~(\ref{listand}). Also by~(\ref{listand}),
\begin{equation*}
\li(x) = {x\over\log x} + O\big( {x\over\log^2x} \big) = {y\over\log
x} + O\big( {x-y\over\log x} + {x\over\log^2x} \big).
\end{equation*}
But since $x-y\le x/\log x$ we have $x\ll y$, and so we see that the
error term is $\ll y/\log^2x$, which establishes the lemma in this
case.

In the remaining case, where $2\le y\le x-x/\log x$, we note that
\begin{equation*}
\li(x) - \li(x-y) = \int_{x-y}^x {dt\over\log t} = {x\over\log x} -
{x-y\over\log(x-y)} + O\bigg( \int_{x-y}^x {dt\over\log^2t} \bigg)
\end{equation*}
by integration by parts. The integral in the error term is over an
interval of length $y$, and the integrand never exceeds
$(\log(x-y))^{-2} \ll (\log x)^{-2}$, and so the error term is $\ll
y/\log^2x$. As for the main term, the fact that $y\le x-x/\log x$ implies
that $(1-y/x)^{-1}\le\log x$, and so we can write
\begin{equation*}
\begin{split}
{x\over\log x} - {x-y\over\log(x-y)} &= {x\over\log x} - {x-y\over\log
x + \log(1-y/x)} \\
&= {x\over\log x} - {x-y\over\log x}\big( 1+O\big (
{\log(1-y/x)^{-1}\over\log x} \big) \big) \\
&= {y\over\log x} + O\big ( {y\over\log^2x} A\big( 1-\frac yx
\big) \big),
\end{split}
\end{equation*}
where we have defined the function $A(t)=t\log t^{-1}/(1-t)$. One can
check that this function $A$ is bounded on the interval $(0,1)$, and
so this error term is simply $O(y/\log^2x)$, which establishes the
lemma.
\end{proof}

As a consequence of Lemma~\ref{calclem}, we see that the asymptotic
formula~(\ref{hypapeq}) in \hypap\ can be restated as
\begin{equation}
\pi(x;q,a) - \pi(x-y;q,a) = {\li(x)-\li(x-y)\over\phi(q)} + O\big(
{y\over\phi(q)\log^2x} \big).
\label{hypapeqre}
\end{equation}
On the other hand, as a consequence of equation~(\ref{asacons}), the
asymptotic formula~(\ref{hypoeq}) in \hypo\ can be restated as
\begin{equation}
\pi(qt+b;T) = {\li(qT+b)-\li(b)\over\phi(q)} + O\big( {q\over\phi(q)}
{T\over\log^2 T} \big).
\label{hypoeqre}
\end{equation}

We are now able to show that \hypo\ and \hypap\ are equivalent statements.

\begin{pflike}{Proof that \hypo\ implies \hypap:}
Let $0<\ep<1$ and $C>1$ be real numbers, let $x$ and $y$ be
sufficiently large real numbers satisfying $y\le x\le y^C$, and let
$a$ and $q$ be coprime integers satisfying $1\le a\le q\le
y^{1-\ep}$. We want to show that the asymptotic
formula~(\ref{hypapeqre}) holds for $\pi(x;q,a)-\pi(x-y;q,a)$.

Suppose first that $x$ and $y$ are integer multiples of $q$. Then
\begin{equation}
\begin{split}
\pi(x;q,a)-\pi(x-y;q,a) &= \3{x-y<p\le x\colon p\equiv a\mod q} \\
&= \3{0<m\le y/q\colon qm+(x-y)+a\hbox{ is prime}} \\
&= \pi(qt+b;T) + O(1),
\end{split}
\label{usehyp}
\end{equation}
where we have defined $b=x-y+a$ and $T=y/q$. For these values of $b$
and $T$ we have
\begin{equation*}
1\le\max\{b,q\} = \max\{x-y+a,q\} \le x+q \le y^C + y^{1-\ep}.
\end{equation*}
If we choose $B>C/\ep$, then this implies for $y$ sufficiently large
\begin{equation*}
1\le\max\{b,q\} < (y^\ep)^B \le \big( \frac yq \big)^B = T^B.
\end{equation*}
Therefore, we may apply \hypo\ with this choice of $B$ to the
expression $\pi(qt+b,T)$. Using the equivalent
formulation~(\ref{hypoeqre}) of \hypo, equation~(\ref{usehyp}) becomes
\begin{equation*}
\begin{split}
\pi(x;q,a)-\pi(x-y;q,a) &= {\li(qT+b)-\li(b)\over\phi(q)} + O\big
( {q\over\phi(q)} {T\over\log^2T} \big) + O(1) \\
&= {\li(x+a)-\li(x-y+a)\over\phi(q)} + O\big( {q\over\phi(q)}
{y/q\over\log^2(y/q)} \big) \\
&= {\li(x)-\li(x-y)\over\phi(q)} + O\big( {a\over\phi(q)} \big) +
O\big( {y\over\phi(q)\log^2(y/q)} \big).
\end{split}
\end{equation*}
Using the assumptions on the sizes of $x$, $y$, $q$, and $a$, the
error terms can be replaced by $O(y/(\phi(q)\log^2x))$, and so we have
derived the desired asymptotic formula~(\ref{hypapeqre}).

This shows that \hypo\ implies \hypap\ in the case where $x$ and $y$
are integer multiples of $q$. However, if we let $x'$ and $y'$ be the
integer multiples of $q$ closest to $x$ and $y$, respectively, then
$\pi(x';q,a)=\pi(x;q,a)+O(1)$ and similarly for
$\pi(x'-y';q,a)$. Therefore \hypo\ implies \hypap\ for any values of
$x$ and $y$ in the appropriate range.
\end{pflike}

\begin{pflike}{Proof that \hypap\ implies \hypo:}
Let $B$ be a positive real number, let $T$ be a sufficiently large
real number, and let $q$ and $b$ be coprime integers satisfying $1\le
q\le T^B$ and $|b|\le T^B$. We want to show that the asymptotic
formula~(\ref{hypoeqre}) holds for $\pi(qt+b;T)$.

Suppose first that $b$ is positive. If we let $a$ denote the smallest
positive integer congruent to $b\mod q$, then
\begin{equation}
\begin{split}
\pi(qt+b;T) &= \3{1\le m\le T\colon qm+b\hbox{ is prime}} \\
&= \3{b<n\le qT+b\colon n\hbox{ is prime, }n\equiv a\mod q} \\
&= \pi(x;q,a) - \pi(x-y;q,a),
\end{split}
\label{useotherhyp}
\end{equation}
where we have defined $x=qT+b$ and $y=qT$. Clearly we have $1\le a\le
q$ and $1\le y\le x$. Moreover, if we choose $C>B+1$, then
\begin{equation*}
x = qT+b \le T^{B+1}+T^B < T^C \le (qT)^C = y^C
\end{equation*}
since $T$ is sufficiently large; and if we also let $\ep=(B+1)^{-1}$,
then
\begin{equation}
q = q^{1-\ep} q^{\ep} \le q^{1-\ep} (T^B)^{\ep} =
(qT)^{1-\ep} \le y^{1-\ep}.
\label{ineqs}
\end{equation}
Therefore we can apply \hypap\ with these values of $C$ and $\ep$ to the
difference $\pi(x;q,a)-\pi(x-y;q,a)$. Using the equivalent
formulation~(\ref{hypapeqre}) of \hypap, equation~(\ref{useotherhyp})
becomes
\begin{equation*}
\begin{split}
\pi(qt+b;T) &= {\li(x)-\li(x-y)\over\phi(q)} + O\big
( {y\over\phi(q)\log^2x} \big) \\
&= {\li(qT+b)-\li(b)\over\phi(q)} + O\big( {qT\over\phi(q)\log^2T}
\big)
\end{split}
\end{equation*}
(since $\log x\ge\log T$), which is the desired asymptotic
formula~(\ref{hypoeqre}).

This shows that \hypap\ implies \hypo\ in the case where $b$ is
positive. Notice that $\pi(qt+b;T)$ counts the number of primes in the
set $\{q+b,2q+b,\dots,qT-q+b,qT+b\}$; on the other hand,
$\pi(qt-(qT+b);T)$ counts the number of primes in the set
$\{-(qT-q+b),\dots,-(q+b),-b\}$, which differs from the aforementioned
set only by the negative signs on each element and a difference of one
element at each end. Consequently, $\pi(qt+b;T) =
\pi(qt-(qT+b);T)+O(1)$, and so if $b$ is so negative that $qT+b$ is
also negative, we can replace $b$ by $-(qT+b)$ and reduce to the case
already considered.

Finally, consider the case where $b$ is negative but $qT+b$ is
positive. Replacing $b$ by $qT+b$ as in the previous paragraph if
necessary, we may assume that $qT+b\ge|b|$. In this case the analogous
equation to~(\ref{useotherhyp}) is
\begin{equation}
\begin{split}
\pi(qt+b;T) &= \3{1\le m\le T\colon qm+b\hbox{ is prime}} \\
&= \3{1\le n\le qT+b\colon n\hbox{ is prime, }n\equiv a\mod q} \\
&\qquad+ \3{b<n\le-1\colon |n|\hbox{ is prime, }n\equiv a\mod q} \\
&= \pi(x_1;q,a) + \pi(x_2;q,q-a),
\end{split}
\label{usesomething}
\end{equation}
with $a$ defined (as above) to be the smallest positive integer
congruent to $b\mod q$, and where we have defined $x_1=qT+b=qT-|b|$
and $x_2=|b|$. Notice that
\begin{equation*}
qT+b\ge|b| \quad\implies\quad qT/2 \ge |b| \quad\implies\quad
x_1=qT-|b|\ge qT/2.
\end{equation*}
Notice also that $q\le(qT)^{1-\ep}$ was shown in
equation~(\ref{ineqs}) (where $\ep=(B+1)^{-1}$ as before). If we
choose a real number $\ep'$ satisfying
$0<\ep'<\min\{(B+1)^{-1},(2B)^{-1},1/3\}$, we see that
\begin{equation*}
q\le(qT)^{1-\ep} < (2x_1)^{1-\ep} < x_1^{1-\ep'}
\end{equation*}
since $T$ is sufficiently large. Consequently we may apply \hypap\ to
$\pi(x_1;q,a)$ with $x=y=x_1$; the equivalent
formulation~(\ref{hypapeqre}) gives us
\begin{equation}
\pi(x_1;q,a) = {\li(x_1)\over\phi(q)} + O\big(
{x_1\over\phi(q)\log^2x_1} \big).
\label{pix1}
\end{equation}
The idea is now to apply a similar argument to the other term
$\pi(x_2;q,q-a)$ in equation~(\ref{usesomething}) when $x_2$ is
reasonably large, and to bound this expression trivially when $x_2$ is
rather small.

In this vein, assume first that
\begin{equation}
x_2>\sqrt T\hbox{ and }x_2^{1-\ep'}\ge q.  \label{x2conds}
\end{equation}
In this case, we can apply equation~(\ref{hypapeqre}) to
$\pi(x_2;q,q-a)$ with $x=y=x_2$, resulting in
\begin{equation*}
\pi(x_2;q,q-a) = {\li(x_2)\over\phi(q)} + O\big(
{x_2\over\phi(q)\log^2x_2} \big).
\end{equation*}
This, together with equation~(\ref{pix1}), means that
equation~(\ref{usesomething}) becomes
\begin{equation*}
\begin{split}
\pi(qt+b;T) &= {\li(x_1)\over\phi(q)} + O\big(
{x_1\over\phi(q)\log^2x_1} \big) + {\li(x_2)\over\phi(q)} + O\big(
{x_2\over\phi(q)\log^2x_2} \big) \\
&= {\li(qT-|b|)+\li(|b|)\over\phi(q)} + O\big(
{qT-|b|\over\phi(q)\log^2x_1} + {|b|\over\phi(q)\log^2x_1} \big) \\
&= {\li(qT+b)-\li(b)\over\phi(q)} + O\big( {q\over\phi(q)}
{T\over\log^2T} \big)
\end{split}
\end{equation*}
(since $\log x_1\ge\log x_2\ge\log\sqrt T$), using the convention
about $\li(b)$ mentioned after equation~(\ref{asacons}).

On the other hand, assume that one of the two
conditions~(\ref{x2conds}) fails for $x_2$. If $x_2\le\sqrt T$ then
certainly $\pi(x_2;q,q-a)\le\sqrt T$. Also, if $x_2^{1-\ep'}<q$ then
\begin{equation*}
\pi(x_2;q,q-a) \le 1+{x_2\over q} \ll x_2^{\ep'} = |b|^{\ep'} \le
(T^B)^{\ep'} < \sqrt T
\end{equation*}
as well. Using this and the asymptotic formula~(\ref{pix1}),
equation~(\ref{usesomething}) now becomes
\begin{equation}
\begin{split}
\pi(qt+b;T) &= {\li(x_1)\over\phi(q)} + O\big(
{x_1\over\phi(q)\log^2x_1} \big) + O(\sqrt T) \\
&= {\li(qT-|b|)\over\phi(q)} + O\big( {qT-|b|\over\phi(q)\log^2x_1} +
\sqrt T \big) \\
&= {\li(qT+b)\over\phi(q)} + O\big( {q\over\phi(q)}
{T\over\log^2T} \big),
\end{split}
\label{insertme}
\end{equation}
since $\log x_1\gg\log T$. 

Now
\begin{equation*}
{-\li(b)\over\phi(q)} = {\li(|b|)\over\phi(q)} \le {x_2\over
q/\log\log q} \le {\max\{\sqrt T,q^{(1-\ep')^{-1}}\}\over q/\log\log q}
\end{equation*}
since $x_2$ fails at least one of the
conditions~(\ref{x2conds}). Because $(1-\ep')^{-1}<3\ep'/2$ by the
restriction $\ep'<1/3$, we have
\begin{equation*}
{-\li(b)\over\phi(q)} \ll \max\{\sqrt T,q^{3\ep'/2}\} \le \max\{\sqrt
T,(T^B)^{3/4B}\} \le T^{3/4} \ll {q\over\phi(q)} {T\over\log^2T}
\end{equation*}
since $\ep'<(2B)^{-1}$. Therefore the term $-\li(b)/\phi(q)$ may be
inserted into the last line of equation~(\ref{insertme}), which
establishes the asymptotic formula~(\ref{hypoeqre}) in this last case.
\end{pflike}

{\smaller\smaller\baselineskip=12pt
\begin{pflike}{Acknowledgements.}
The author would like to acknowledge the support of National Science
Foundation grant number DMS 9304580 and Natural Sciences and
Engineering Research Council grant number A5123. The author would
also like to thank the anonymous referee, whose thorough
reading of the manuscript greatly reduced the number of
typographical errors herein.\par
\end{pflike}
}
\bibliography{asympoly}

\providecommand{\bysame}{\leavevmode\hbox to3em{\hrulefill}\thinspace}
\begin{thebibliography}{10}

\bibitem{AGP}
W.~R. Alford, A.~Granville, and C.~Pomerance, \emph{There are infinitely many
  {C}armichael numbers}, Ann. of Math. (2) \textbf{139} (1994), no.~3,
  703--722.

\bibitem{BatHor}
P.~T. Bateman and R.~A. Horn, \emph{A heuristic asymptotic formula concerning
  the distribution of prime numbers}, Math. Comp. \textbf{16} (1962), 363--367.

\bibitem{DarMarTen}
C.~Dartyge, G.~Martin, and G.~Tenenbaum, \emph{Polynomial values free of large
  prime factors}, preprint.

\bibitem{Dic}
L.~E. Dickson, \emph{A new extension of {D}irichlet's theorem on prime
  numbers}, Messenger of Math. \textbf{33} (1903), 155--161.

\bibitem{FriGra:LttEoP}
J.~Friedlander and A.~Granville, \emph{Limitations to the equi-distribution of
  primes. {I}}, Ann. of Math. (2) \textbf{129} (1989), no.~2, 363--382.

\bibitem{HalRic:SM}
H.~Halberstam and H.-E. Richert, \emph{Sieve methods}, Academic Press, London,
  1974.

\bibitem{HarLit}
G.~H. Hardy and J.~E. Littlewood, \emph{Some problems of ``partitio
  numerorum''; {III}: On the expression of a number as a sum of primes}, Acta
  Math. \textbf{44} (1923), 1--70.

\bibitem{HilTen}
A.~Hildebrand and G.~Tenenbaum, \emph{Integers without large prime factors}, J.
  Th\'eor. Nombres Bordeaux \textbf{5} (1993), no.~2, 411--484.

\bibitem{Hux:ANoPC}
M.~N. Huxley, \emph{A note on polynomial congruences}, Recent Progress in
  Analytic Number Theory, vol.~1 (Durham, 1979), Academic Press, London--New
  York, 1981, pp.~193--196.

\bibitem{Nag:GdTdT}
T.~Nagel, \emph{G\'en\'eralisation d'un th\'eor\`eme de {T}chebycheff}, J.
  Math. Pures Appl. (8) \textbf{4} (1921), 343--356.

\bibitem{NairTen}
M.~Nair and G.~Tenenbaum, \emph{Short sums of certain arithmetic functions},
  Acta Math. \textbf{180} (1998), no.~1, 119--144.

\bibitem{SchSie}
A.~Schinzel and W.~Sierpi\'nski, \emph{Sur certaines hypoth\`eses concernant
  les nombres premiers}, Acta Arith. \textbf{4} (1958), 185--208; {\it
  erratum\/} {\bf 5} (1958), 259.

\bibitem{Ste:OtNoSoPCaTE}
C.~L. Stewart, \emph{On the number of solutions of polynomial congruences and
  {T}hue equations}, J. Amer. Math. Soc. \textbf{4} (1991), no.~4, 793--835.

\bibitem{Wir:DAVvSuMF}
E.~Wirsing, \emph{Das asymptotische {V}erhalten von {S}ummen \"uber
  multiplikative {F}unktionen}, Math. Ann. \textbf{143} (1961), 75--102.

\end{thebibliography}
\bibliographystyle{../amsplain}
\end{document}